\documentclass[aap,MSNbibl,citesort,seceqn,dvips]{arximspdf}
\usepackage{graphicx}

\doi{10.1214/11-AAP784}
\volume{22}
\issue{2}
\pubyear{2012}
\firstpage{734}
\lastpage{791}

\makeatletter
\newcommand{\bbref}[1]{(\ref{#1})}
\newtheorem{theorem}{Theorem}[section]
\newtheorem{lemma}[theorem]{Lemma}
\newtheorem{lemm}{Lemma}
\newtheorem{cor}[theorem]{Corollary}
\newtheorem{prop}[theorem]{Proposition}

\newproclaim{assu}[theorem]{Assumption}
\newproclaim{remark}[theorem]{Remark}
\newproclaim{definition}[theorem]{Definition}
\newproclaim{defn}[theorem]{Definition}
\newproclaim{cond}[theorem]{Condition}
\newproclaim{pr}[theorem]{Property}
\newproclaim{example}{Example}
\newcommand{\eqref}[1]{(\ref{#1})}
\DeclareMathAlphabet{\mathpzc}{OT1}{pzc}{m}{it}
\newcommand{\curvp}{\mathsf{p}}
\newcommand{\curvm}{\bolds{\mathsf{m}}}
\newcommand{\curvk}{\bolds{\mathsf{k}}}
\newcommand{\III}{\mathbb{I}}
\newcommand{\curvi}{\mathpzc{i}}
\newcommand{\curvq}{\mathbf{\mathpzc{q}}}
\newcommand{\curvz}{\mathbf{\mathpzc{z}}}
\newcommand{\IE}{\mathbb{E}}
\newcommand{\IZ}{\mathbb{Z}}
\newcommand{\IP}{\mathbb{P}}
\newcommand{\E}{\mathbb{E}}
\newcommand{\R}{\mathbb{R}}
\newcommand{\N}{\mathbb{N}}
\newcommand{\Df}{\doteq}
\newcommand{\II}{\mathbf{I}}
\newcommand{\JJ}{\mathbf{J}}
\newcommand{\KK}{\mathbf{K}}
\newcommand{\BB}{\mathbf{B}}
\newcommand{\LL}{\mathbf{L}}
\newcommand{\bfi}{{\bolds{\Phi}}^r}
\newcommand{\bfij}[1]{{\bolds{\Phi}}^{#1,r}}
\newcommand{\bfijl}[2]{{\bolds{\Phi}}^{#2,r}_{#1}}
\makeatother

\begin{document}
\begin{frontmatter}

\title{Controlled stochastic networks in heavy traffic: Convergence of value functions}
\runtitle{Convergence of value functions}

\begin{aug}
\author[A]{\fnms{Amarjit} \snm{Budhiraja}\thanksref{t1}\ead[label=e1]{budhiraj@email.unc.edu}}
\and
\author[B]{\fnms{Arka P.} \snm{Ghosh}\corref{}\thanksref{t2}\ead[label=e2]{apghosh@iastate.edu}}
\thankstext{t1}{Supported in part by the NSF
(DMS-10-04418),
Army Research Office (W911NF-0-1-0080, W911NF-10-1-0158) and the
US-Israel Binational Science Foundation (2008466).}
\thankstext{t2}{Supported in part by NSF Grant DMS-06-08634.}
\runauthor{A. Budhiraja and A. P. Ghosh}
\affiliation{University of North Carolina and Iowa State University}
\address[A]{Department of Statistics\\
\quad and Operations Research\\
University of North Carolina\\
Chapel Hill, North Carolina 27599-3260\\
USA\\
\printead{e1}} 
\address[B]{Department of Statistics\\
3216 Snedecor Hall\\
Iowa State University\\
Ames, Iowa 50011-1210\\
USA\\
\printead{e2}}
\end{aug}

\received{\smonth{3} \syear{2010}}
\revised{\smonth{12} \syear{2010}}

%
\begin{abstract}
Scheduling control problems for a family of unitary networks
under heavy traffic with general interarrival and service times,
probabilistic routing and an infinite horizon
discounted linear holding cost are studied. Diffusion control problems,
that have been proposed as approximate models for the
study of these critically loaded controlled stochastic networks, can be
regarded as formal scaling limits of such stochastic systems. However,
to date, a rigorous limit theory that justifies the use of such
approximations for a general family of controlled networks has been
lacking. It is shown that, under broad conditions, the value function
of the suitably scaled network control
problem converges to that of the associated diffusion control problem.
This scaling limit result, in addition to giving a precise mathematical
basis for the above approximation approach, suggests a general strategy
for constructing near optimal controls for the physical stochastic
networks by solving the associated diffusion control problem.
\end{abstract}

%

\end{frontmatter}

\section{Introduction}\label{sec:intro}
As an approximation to control problems for critically-loaded
stochastic networks,
Harrison (in \cite{harri2}, see also \cite{harri1,harri-canon}) has
formulated a stochastic control problem
in which the state process is driven by a multidimensional Brownian
motion along with an additive control that
satisfies certain feasibility and nonnegativity constraints. This
control problem, that is, usually referred to as the
\textit{Brownian Control Problem} (BCP) has been one of the key
developments in the heavy traffic theory of controlled stochastic
processing networks (SPN). BCPs can be regarded as formal scaling
limits for a broad range of
scheduling and sequencing control problems for multiclass queuing
networks. Finding optimal
(or even near-optimal)
control policies for such networks---which may have quite general
non-Markovian primitives, multiple server capabilities and
rather complex routing geometry---is in general prohibitive. In that
regard, BCPs that provide significantly more tractable approximate models
are very useful. In this diffusion approximation approach to policy
synthesis, one first finds an optimal
(or near-optimal) control for the BCP which is then suitably
interpreted to construct a
scheduling policy for the underlying physical network. In recent years
there have been many works
\cite{ata-kumar,bellwill,bellwill2,BudGho,meyn,ward-kumar,chen-yao,dai-lin}
that consider specific
network models for which the associated BCP is explicitly solvable
(i.e., an optimal control process can be written
as a known function of the driving Brownian motions) and, by suitably
adapting the solution to the underlying network,
construct control policies that are asymptotically (in the heavy
traffic limit) optimal. The paper \cite{KuMa} also carries
out a~similar program for the crisscross network where the state--space
is three dimensional, although an explicit solution for the BCP here is not
available.

Although now there are several papers which establish a rigorous
connection between a network control problem
and its associated BCP by exploiting the explicit form of the solution
of the latter, a systematic theory which justifies the use of BCPs as
approximate models has been missing.
In a recent work \cite{BudGho2} it was shown that for a large family of
\textit{Unitary Networks}
(following terminology of \cite{Will-Bram-2work}, these are networks
with a structure as described in Section \ref{secsetup}), %
with general interarrival and service times, probabilistic routing and an
infinite horizon discounted linear holding cost, the cost associated
with any admissible control policy for the network
is asymptotically, in the heavy traffic limit, bounded below by the
value function of the BCP. This inequality, which
provides a useful bound on the best achievable asymptotic performance
for an admissible control policy, was a key step
in developing a rigorous general theory relating BCPs with SPN in heavy traffic.

The current paper is devoted to the proof of the reverse inequality.
The network model is required to satisfy
assumptions made in \cite{BudGho2} (these are summarized above Theorem
\ref{ab937}). In addition, we impose a nondegeneracy
condition (Assumption~\ref{non-deg}), a condition on the underlying
renewal processes regarding probabilities of deviations
from the mean (Assumption~\ref{ldp}) and regularity of a certain
Skorohod map (Assumption~\ref{ab1012}) (see next paragraph for a discussion
of these conditions).
Under these assumptions we prove that the value function of the BCP
is bounded below by the heavy traffic limit (limsup) of the value
functions of the network control problem (Theorem~\ref{main1016}). Combining
this with the result obtained in \cite{BudGho2} (see Theorem~\ref
{ab937}), we obtain
the main result of the paper (Theorem \ref{maincorr}). This theorem
says that, under broad conditions, the value function of
the network control problem converges to that of the BCP. This result
provides, under general conditions, a rigorous basis for regarding
BCPs as approximate models for critically loaded stochastic networks.

Conditions imposed in this paper allow for a wide range of SPN
models. Some such models, whose description is taken from
\cite{Will-Bram-2work}, are discussed in detail in Examples
\ref{exammodel}{(a)--(c)}. We note that our approach does not
require the BCP to be explicitly solvable and the result covers many
settings where explicit solutions are unavailable. Most of the
conditions that we impose are quite standard and we only comment
here on three of them: Assumptions \ref{ab148}, \ref{G-column} and
\ref{ab1012}. Assumption \ref{ab148} says that each buffer is
processed by at least one basic activity (see Remark \ref{rem-HT}).
This condition, which was introduced in~\cite{Will-Bram-2work}, is
fundamental for our analysis. In fact, \cite{Will-Bram-2work}~has
shown that without this assumption even the existence of a
nonnegative workload matrix may fail. Assumption \ref{G-column} is
a natural condition on the geometry of the underlying network.
Roughly speaking, it says that a nonzero control action leads to
a~nonzero state displacement. Assumption \ref{ab1012} is the third
key
requirement in this work. It says that the Skorohod problem associated
with a certain reflection matrix $D$ [see equation \eqref{refmat} for the
definition of $D$] is well posed and the associated Skorohod map is
Lipschitz continuous.
As Example \ref{exammodel} discusses, this condition holds for a broad
family of networks (including \textit{all } multiclass open queuing
networks, as well as a large family of parallel server networks and
job-shop networks).

The papers \cite{ata-kumar,bellwill,bellwill2,BudGho,ward-kumar,dai-lin} noted earlier, that treat the setting of explicitly solvable
BCP, do much more
than establish convergence of value functions. In particular, these
works give an explicit implementable control policy for the underlying
network that
is asymptotically optimal in the heavy traffic limit. In the
generality treated in the current work, giving explicit recipes (e.g.,
threshold type policies)
is unfeasible, however, the policy sequence constructed in Section \ref
{subconstruct} suggests a general approach for building near
asymptotically optimal
policies for the network \textit{given} a near optimal control for the
BCP. Obtaining near optimal controls for the BCP in general requires
numerical approaches (see, e.g., \cite{DuKu,Kushbook,meyn-book}),
discussion of which is beyond the scope of the current work.

We now briefly describe some of the ideas in the proof of the main
result---Theorem \ref{main1016}. We begin by choosing, for an
arbitrary $\varepsilon> 0$, a
suitable $\varepsilon$-optimal control $\tilde Y$ for the BCP and then, using
$\tilde Y$, construct a sequence of control policies
$\{T^r\}_{r\ge1}$ for the network model such that the (suitably
scaled) cost associated with $T^r$ converges to that associated with
$\tilde Y$,
as $r\to\infty$. This yields the desired reverse inequality.
One of the
key difficulties is in the translation of a given control for the BCP
to that for the physical network. Indeed, a (near) optimal control for
the BCP can be a very general adapted process with RCLL paths. Without
additional information on such
a stochastic process, it is not at all clear how one adapts and applies
it to a given network model. A control policy for the network
needs to specify how each server distributes its effort among various
job classes at any given time instant. By a series of approximations
we show that one can find a~rather simple $\varepsilon$-optimal control for
the BCP, that is, easy to interpret and implement
on a network control problem. As a first step, using PDE
characterization results for general singular control problems with
state constraints from \cite{AtBu} (these, in particular, make use of
the nondegeneracy
assumption---Assumption \ref{non-deg}), one can argue that a
near-optimal control can be taken to be adapted to the driving
Brownian motion and be further assumed to have moments that are
subexponential in the time variable
(see Lemma \ref{cont-bcp}). Using results from \cite{BuRo}, one can
perturb this control so that it
has continuous sample paths without significantly affecting the cost.
Next, using ideas developed by Kushner and Martins \cite{KuMa} in the
context of a two-dimensional BCP,
one can further approximate such a control by a process with a fixed
(nonrandom) finite
number of jumps that take values in a finite set. Two main requirements
(in addition to the usual adaptedness condition)
for such a process to be an admissible control of a BCP (see Definition
\ref{BCP}) are the nonnegativity constraints \eqref{ab934}
and state constraints \eqref{ab440}.
It is relatively easy to construct a pure jump process that satisfies
the first requirement of admissibility, namely, the nonnegativity constraints,
however, the nondegenerate Brownian motion in the dynamics rules out
the satisfaction of the second requirement, that is, state constraints,
without additional modifications.
This is where the regularity assumption on a certain Skorohod map
(Assumption \ref{ab1012}) is used.
The pure jump control is modified in a manner such that in between
successive jumps one uses the Skorohod map
to employ minimal control needed in order to respect state
constraints. Regularity of the Skorohod problem
ensures that this modification does not change the associated cost
much. The Skorohod map also plays a key
role in the weak convergence arguments used to prove convergence of
costs. 
The above construction is the essential content of Theorem \ref{main-bcp}.
The $\varepsilon$-optimal control that we use for the construction of the
policy sequence requires two additional modifications [see part (iii)
of Theorem \ref{main-bcp} and below \eqref{111}] which facilitate adapting
such a control for the underlying physical network and in some weak
convergence proofs, but we leave that discussion
for later in Section \ref{nearopt} (see Remark \ref{on3} and above
Theorem \ref{newmain518}).

Using a near-optimal control $\tilde Y$ of the form given in Section
\ref
{nearopt} (cf. Theorem~\ref{newmain518}), we then proceed to construct
a sequence
of policies $\{T^r\}$ for the underlying network. The key relation that
enables translation of $\tilde Y$ into
$\{T^r\}$ is \eqref{Y-defn} using which one can loosely interpret
$\tilde
Y(t)$ as the asymptotic
deviation, with suitable scaling, of $T^r(t)$ from the nominal
allocation $x^*t$ (see Definition \ref{defn-ht} for the definition of nominal
allocation vector). Recall that $\tilde Y$ is constructed by modifying,
through a Skorohod constraining mechanism,
a pure jump process (say, $\tilde Y_0$). In particular, $\tilde Y$~has
sample paths that are, in general, discontinuous.
On the other hand, note that an admissible policy~$T^r$ is required to
be a Lipschitz function (see Remark \ref{rem921}). This suggests the
following construction for $T^r$. Over time periods
(say, $\Delta t$)
of constancy of~$\tilde Y_0$ one should use the nominal allocation (i.e.,
$x^*\Delta t$), while jump-instants should be
stretched into periods of length of order $r$ (note that in the scaled
network, time is accelerated by a factor of $r^2$
and so such periods translate to intervals of length $1/r$ in the
scaled evolution and thus are negligible)
over which a nontrivial control action is employed as dictated by the
jump vector (see Figure~\ref{policyfig} for a more complete description).
This is analogous to the idea of a discrete review policy proposed by
Harrison \cite{bigstep} (see also \cite{ata-kumar} and references therein).
There are
some obvious difficulties with the above prescription, for example, a nominal
allocation corresponds to the average behavior of the system and for a
given realization is feasible only when the buffers are nonempty. Thus,
one needs to modify the above construction to incorporate
idleness, that is, caused due to empty buffers. The effect of such a
modification is, of course, very similar to that of a Skorohod
constraining mechanism and it is tempting to hope that the deviation
process corresponding to this modified policy converges to
$\tilde Y$ (in an appropriate sense), as $r\to\infty$. However, without
further modifications, it is not obvious that the reflection terms
that are produced from the idling periods under this policy are
asymptotically consistent with those obtained from the Skorohod
constraining mechanism applied to (the state process corresponding to)
$\tilde Y_0$. The additional modification [see \eqref{ab543}] that we make
roughly says that jobs are processed from a given buffer over a small
interval $\Delta$, only if at the beginning of this interval there are
a ``sufficient'' number of jobs in the buffer. This idea of safety
stocks is not new and has been used in previous works (see, e.g.,
\cite{bellwill,bellwill2,ata-kumar,BudGho,meyn-book}). The
modification, of course, introduces a~somewhat nonintuitive idleness
even when there are jobs that
require processing. However, the analysis of Section \ref{secproof}
shows that this idleness does not significantly affect the asymptotic
cost. The above very rough sketch of construction of~$T^r$ is made
precise in Section \ref{subconstruct}.

The rest of the paper is devoted to showing that the cost associated
with~$T^r$ converges to that associated with $\tilde Y$.
It is unreasonable to expect convergence of controls (e.g., with
the usual Skorohod topology)---in particular, note that $T^r$ has Lipschitz paths for every $r$ while
$\tilde Y$ is a (modification of) a~pure jump
process -- however, one finds that the convergence of costs holds. This
convergence proof, and the related weak convergence
analysis, is carried out in Sections \ref{secconvprf} and \ref{weakcgce}.

The paper is organized as follows. Section \ref{secsetup} describes
the network structure, all the associated stochastic processes and the
heavy-traffic assumptions as well as the other assumptions of the
paper. The section also presents the SPN control problem, that is,
considered here, along with the main result of the paper (Theorem \ref
{maincorr}). Section \ref{nearopt} constructs (see Theorem \ref
{newmain518}) a~near-optimal control policy for the BCP which can be
suitably adapted to the network control problem. In Section \ref
{secproof} the near-optimal control policy from
Section \ref{nearopt} is used to obtain a sequence of admissible
control policies for the scaled SPN. The main result of the section is
Theorem~\ref{mainweak}, which establishes weak convergence of various
scaled processes. Convergence of costs (i.e., Theorem \ref{jrtoj}) is
an immediate
consequence of this weak convergence result. Theorem~\ref{maincorr}
then follows on combining Theorem~\ref{jrtoj} with results of \cite
{BudGho2} (stated as Theorem \ref{ab937} in the current work). Finally,
the \hyperref[appen]{Appendix} collects proofs of some auxiliary results.

The following notation will be used.
The space of reals (nonnegative reals),
positive (nonnegative) integers will be denoted by $\R$ ($\R_+$),
$\N$ ($\N_0$), respectively.
For $m\geq1$ and $\theta\in(0,\infty),$ $\mathcal{C}^m [\mathcal{C}
^m_{\theta}]$ will denote the space of
continuous functions from $[0, \infty)$ (resp. $[0, \theta]$) to
$\R^m$ with the topology of uniform convergence on compacts (resp.
uniform convergence). Also, $\mathcal{D}^m [\mathcal{D}^m_{\theta}]$
will denote
the space of right continuous functions with left limits, from $[0,
\infty)$ (resp. $[0, \theta]$) to $\R^m$ with the usual Skorohod
topology. For $y \in\mathcal{D}^m$ and $t, \delta> 0$, we write
$\sup_{0\le s \le t} |y(s)| = |y|_{\infty, t}$ and $\sup_{0\le s_1
\le s_2 \le t, |s_1 - s_2| \le\delta} |y(s_1) - y(s_2)| =
\varpi_y^t(\delta)$, where for $z = (z_1, \ldots, z_m)' \in\R^m$, $|z|^2
= \sum_{i=1}^m |z_i|^2.$
All vector inequalities are to be interpreted component-wise. We will
call a function $f\in\mathcal{D}^m$ nonnegative if $f(t)\ge0$ for all
$t\in{\R}_{+}$. A function $f\in\mathcal{D}^m$ is called
nondecreasing if it is nondecreasing in each component. All
(stochastic) processes in this work will have sample
paths that are right continuous and have left limits, and thus can be
regarded as $\mathcal{D}^m$-valued random variables with a~suitable $m$.
For a Polish space $\mathcal{E}$,
$\mathcal{B}(\mathcal{E})$ will denote the corresponding Borel
sigma-field. Weak
convergence of $(\mathcal{E},
\mathcal{B}(\mathcal{E}))$ valued random variables $Z_n$ to $Z$
will be denoted as $Z_n\Rightarrow Z$.
Sequence of processes $\{Z_n\}$ is tight if and only if the measures
induced by
$Z_n$'s on $(\mathcal{D}^m, \mathcal{B}(\mathcal{D}^m))$ form a
tight sequence. A~sequence of processes with paths in $\mathcal{D}^m$ ($m\geq1$) is
called $\mathcal{C}$-tight if it is tight in $\mathcal{D}^m$ and any weak
limit point of the sequence has paths in ${\mathcal{C}}^m$ almost
surely (a.s.).
For processes $\{Z_n\}$, $Z$ defined on a common probability
space, we say that $Z_n$ converge to $Z$, uniformly on compact
time intervals (u.o.c.), in probability (a.s.) if for all $t
> 0$, $\sup_{0\leq s\leq t} |Z_n(s) - Z(s)|$ converges to zero in
probability (resp. a.s.). To ease the notational burden, standard
notation (that follow \cite{bram-will-1,Will-Bram-2work}) for
different processes are
used (e.g., $Q$ for queue-length,~$I$ for idle time, $W$ for workload
process etc.). We also use standard notation, for example, $\bar W,
\hat
W$, to denote fluid scaled, respectively, diffusion scaled, versions of
various processes of interest [see \eqref{fl-scaled} and \eqref{diff-scaled}].
All vectors will be column vectors.
An $m$-dimensional vector with all entries $1$
will be denoted by
${\mathbf{1}}_{m}$. For a
vector $a$, $\operatorname{diag}(a)$ will denote the diagonal matrix such that
the vector of its diagonal entries is $a$. $M'$ will denote the
transpose of a matrix $M$. Also, $C_i, i=0, 1, 2, \ldots,$ will denote generic constants whose values may change from one proof
to the next.

\section{Multiclass queueing networks and the control problem}
\label{secsetup}
Let $(\Omega, \mathcal{F}, \IP)$ be a probability space. All the random
variables associated with the network model described below
are assumed to be defined on this probability space.
The expectation operation under $\IP$ will be denoted by $\IE$.

\subsection*{Network structure}

\begin{figure}

\includegraphics{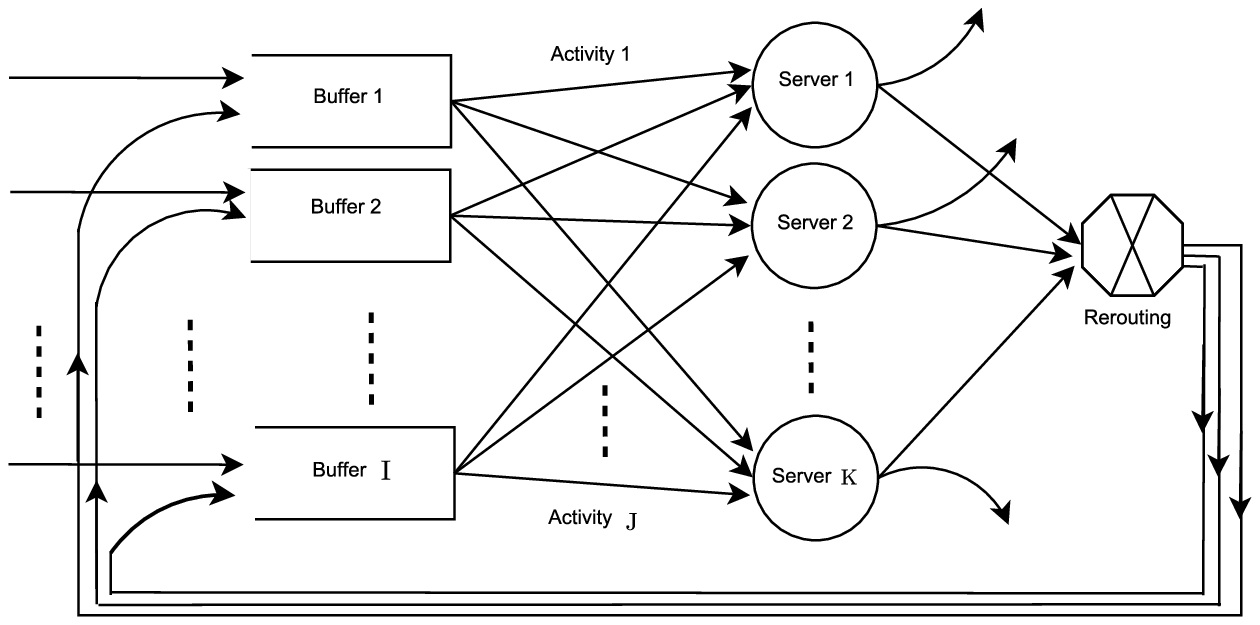}

\caption{A network with $\II$ buffers, $\JJ$ activities, $\KK$ servers
and probabilistic routing (given by the matrix $P$).} \label{fig1-general}
\vspace*{6pt}
\end{figure}

We begin by introducing the family
of stochastic processing network models that will be considered
in this work. We closely follow the terminology and notation used in
\cite{harri2,harri1,harri-canon,bellwill2,Will-Bram-2work,bram-will-1}. The
network has~$\II$ infinite capacity buffers (to store $\II$ many
different classes of jobs) and $\KK$ nonidentical servers for
processing jobs. Arrivals of jobs, given in terms of suitable renewal
processes, can be from outside the system
and/or from the internal rerouting of jobs that have already been
processed by some server. Several different servers may process jobs
from a particular buffer. Service from a given buffer $i$ by a given
server $k$ is called an \textit{activity}. Once a job starts being
processed by an
activity, it must complete its service with that activity, even if
its service is interrupted for some time (e.g., for preemption by a
job from another buffer). When service of a partially completed job
is resumed, it is resumed from the point of preemption---that is, the job
needs only the remaining service time from the server to get
completed (preemptive-resume policy). Also, an activity must
complete service of any job that it started before starting another
job from the same buffer. An activity always selects the oldest job
in the buffer that has not yet been served, when starting a new service
[i.e., First In First Out (FIFO) within class].
There are $\JJ$
activities [at most one activity for a server-buffer pair $(i,k)$,
so that $\JJ\le\II\cdot\KK$]. Here the integers $\II,\JJ,\KK$
are strictly positive. Figure \ref{fig1-general} gives a schematic for
such a model.
%

Let $\mathbb{I} = \{1, \ldots, \II\}$, $\mathbb{J} = \{1, \ldots,
\JJ\}$
and $\mathbb{K} = \{1, \ldots, \KK\}$. The correspondence between the
activities and buffers, and
activities and servers are described by two matrices $C$ and $A$
respectively. $C$ is an $\II\times\JJ$ matrix with $C_{i j}=1$ if
the $j$th activity processes jobs from buffer $i$, and $C_{i j}=0$
otherwise. The matrix~$A$ is $\KK\times\JJ$ with $A_{k j}=1$ if the
$k$th server is associated with the $j$th activity, and $A_{k
j}=0$ otherwise. Each activity associates one buffer and one
server, and so each column of $C$ has exactly one 1 (and
similarly, every column of~$A$ has exactly one~1).
We will
further assume that each
row of $C$ (and $A$) has at least one~1, that is, each buffer is
processed by (server is processing,
resp.) at least one activity.
For $j \in\mathbb{J}$, let $\sigma(j) \equiv(\sigma_1(j), \sigma
_2(j)) = (i,k)$,
if activity $j$ corresponds to the $k$th server processing class $i$
jobs. Let, for $k \in\mathbb{K}$,
$\mathbb{J}(k) \Df\{j \in\mathbb{J}\dvtx  \sigma_2(j) = k\}$ and
$\mathbb{I}(k) \Df\{\sigma_1(j)\dvtx  j \in\mathbb{J}(k)\}$.
Thus, for the $k$th server, $\mathbb{J}(k)$ denotes the set of activities
that the server can perform, and $\mathbb{I}(k)$ represents the
corresponding buffers from which the jobs can be processed.

\subsection*{Stochastic primitives}

We are interested in the study of networks that are nearly critically
loaded. Mathematically, this is modeled by considering a
sequence of networks $\{\mathcal{N}^r\}$ that ``approach heavy
traffic,'' as $r
\to\infty$, in the sense of
Definition~\ref{defn-ht} below. Each network in the sequence has
identical structure, except for the rate
parameters that may depend on $r$.
Here $r\in\mathbb{S}\subseteq
{\R}^{+}$, where $\mathbb{S}$ is a countable set: $\{r_1, r_2,\ldots
\}$
with $1\le r_1<r_2<\cdots$ and $r_n \rightarrow\infty$, as
$n\rightarrow\infty$. One thinks of the physical
network of interest as the $r$th network
embedded in this sequence, for a fixed large value of $r$. For
notational simplicity, throughout the paper,
we will write the limit along the sequence $r_n$ as $n\rightarrow
\infty$ simply as ``$r\rightarrow\infty$.'' Also, $r$ will always be
taken to be an element of $\mathbb{S}$ and, thus, hereafter the
qualifier $r \in\mathbb{S}$ will not be stated explicitly.

The $r$th network $\mathcal{N}^r$ is described as follows.
If the $i$th
class ($i\in\mathbb{I}$) has exogenous job arrivals,
the interarrival times of such jobs are given by a~sequence of
nonnegative random variables $\{u_{i}^r(n)\dvtx  n\ge1\}$ that are i.i.d with
mean and standard deviation $1/{\alpha_{i}^r}, \sigma_{i}^{u,r}
\in(0,\infty)$ respectively. Let, by relabeling if needed, the
buffers with exogenous arrivals correspond to $i \in\{ 1, \ldots,
\II'\} :=\mathbb{I}'$, where $\II' \le\II$. We set ${\alpha
_{i}^r}, \sigma_{i}^{u,r}
= 0$ and $u^r_i(n)=\infty, n\ge1$, for $i \in\mathbb{I}\setminus
\mathbb{I}'$.
Service times for the $j$th type of activity (for
$j \in\mathbb{J}$) are given by a sequence of nonnegative random
variables $\{v_{j}^r(n)\dvtx  n\ge1\}$ that are i.i.d. with mean and standard
deviation $1/{\beta_{j}^r}, \sigma_{j}^{v,r} \in(0,\infty)$
respectively. We will assume that the above random variables are
in fact strictly positive, that is,
%
\begin{equation}\label{stpos1029} \mbox{for
all } i \in\mathbb{I}, j \in\mathbb{J},\qquad \IP\bigl(u^r_i(1)
>0\bigr)=\IP\bigl(v^r_j(1) >0\bigr)=1.
\end{equation}
We will further impose the following uniform integrability condition:
%
\begin{equation}\label{ui1032}
\begin{tabular}{p{240pt}@{}}
the collection $\{(u_i^r(1))^2, (v_j^r(1))^2; r \ge1, j \in
\mathbb{J}, i \in\mathbb{I}' \}$ is uniformly integrable.
\end{tabular}
\end{equation}
%
Rerouting of jobs completed by the $j$th activity is specified
by a sequence of $(\II+1)$-dimensional vector
$\{(\phi^{j,r}_0(n),\phi^{j,r}(n))',n\ge1\}$, where
$\phi^{j,r}(n)=({\phi}^{j,r}_i(n)\dvtx  i \in\mathbb{I})$. For each $j\in
\mathbb{J}$ and $i\in
\mathbb{I}\cup\{0\}$, 
${\phi}^{j,r}_i(n) =1$ if the $n$th completed job by activity $j$
gets rerouted to buffer $i$, and takes the value zero otherwise,
where $i=0$ represents jobs leaving the system. It is assumed that
for each fixed $r$, $\{(\phi^{j,r}_0(n),\phi^{j,r}(n)),n\ge1\}$, $j
\in
\mathbb{J}$,
are (mutually) independent sequences of i.i.d
$\operatorname{Multinomial}_{(\II+1)}(1, (p^j_0,p^{j}))$, where
$p^{j}=(p^{j}_i\dvtx i=1,\ldots,\II)$. That, in particular, means, for
$j \in\mathbb{J}, n \ge1$, $ \sum_{i=0}^{\II}
{\phi}^{j,r}_i(n)=\sum_{i=0}^{\II} p^{j}_i=1$. Furthermore, for fixed
$j \in\mathbb{J}, i_1, i_2 \in\mathbb{I}$,
%
\begin{equation}
\label{sigma-phi}
\operatorname{Cov}(\phi_{i_1}^{j,r}(n),\phi
_{i_2}^{j,r}(n))=\sigma^{\phi_{j}}_{i_1 i_2}=
- p_{i_1}^jp_{i_2}^j + p_{i_1}^{j}\delta_{i_1,i_2},
\end{equation}
where $\delta_{i_1,i_2}$ is $1$ if $i_1=i_2$ and $0$ otherwise.
We also assume that, for each $r$, the random variables
%
\begin{equation}\label{ass1120}
\begin{tabular}{p{230pt}@{}}
$\{u_{i}^r(n),v_{j}^r(n),\phi^{j,r}_0(n), \phi^{j,r}(n), n \ge1, i
\in\mathbb{I}, j \in\mathbb{J}\}$
are mutually independent.
\end{tabular}
\end{equation}

Next we introduce the primitive renewal processes, $(E^r, S^r
)$, that describe the state dynamics. The process
$(E^r_1, \ldots, E^r_{\II'})$ is the $\II'$-dimensional exogenous
arrival process, that is, for each $i \in\mathbb{I}'$, $E_i^r(t)$ is a
renewal process which denotes the number of jobs that have arrived
to buffer $i$ from outside the system over the interval $[0, t]$. For
class $i$
to which there are no exogenous arrivals (i.e., $i\in\mathbb
{I}\setminus
\mathbb{I}'$), we set
$E_i^r(t)=0$ for all $t\ge0$. We will denote the process
$(E^r_1, \ldots, E^r_{\II})'$ by $E^r$.
For each
activity $j \in\mathbb{J}$, $S_j^r(t)$ denotes the number of complete
jobs that
could be processed by activity $j$ in $[0,t]$ \textit{if the
associated server worked continuously and exclusively on jobs from
the associated buffer in $[0, t]$ and the buffer had an infinite
reservoir of jobs}. The vector $(S_1^r, \ldots,
S_{\JJ}^r)'$ is denoted by $S^r$. More precisely,
for $i \in\mathbb{I}, j \in\mathbb{J}, m \ge1$, let
%
\begin{equation}\xi_i^r(m) \Df\sum_{n=1}^{m} u_{i}^r(n),\qquad \eta
_j^r(m) \Df
\sum_{n=1}^{m} v_{j}^r(n).
\end{equation}
We set
$\xi_i^r(0) = 0, \eta_j^r(0) = 0$.
Then $E_i^r$, $S_j^r$ are renewal processes given as follows. For $t
\ge0$,
%
\begin{equation}\label{prim150}
\qquad E_i^r(t)= \max\{m\ge0\dvtx  \xi_i^r(m) \le t\},\quad S_j^r(t) =
\max\{m \ge1\dvtx  \eta_j^r(m) \le t\}.\hspace*{-25pt}
\end{equation}
Finally, we introduce the routing sequences. Let
$\bfij{j}_i(n)$ denote the number of jobs that
are routed to the $i$th buffer, among the first $n$ jobs completed by
activity~$j$. Thus, for $i \in\mathbb{I}, j \in\mathbb{J}$,
%
\begin{equation}
\bfij{j}_{i}(n)=\sum_{m=1}^{n} {\phi}_i^{j,r}(m),
\qquad n=1,2,\ldots.\vadjust{\goodbreak}
\end{equation}
We will denote the $\II$-dimensional
sequence $\{(\bfij{j}_1(n), \ldots, \bfij{j}_{\II}(n))'\}$
corresponding to routing of jobs completed by the $j$th activity by
$\{\bfij{j}(n)\}$. Also, $ \bfi(n) $ will denote the $\II\times\JJ$
matrix $(\bfij{1}(n), \bfij{2}(n)
,\ldots, \bfij{\JJ}(n) )$.

\subsection*{Control}

A \textit{Scheduling policy} or \textit{control} for the $r$th SPN is
specified by a~nonnegative, nondecreasing $\JJ$-dimensional process
$T^r=\{(T_1^r(t),\ldots, T_{\JJ}^r(t))',\allowbreak t\ge0\}$. For any
$j \in\mathbb{J}, t \ge0$, $T_j^r(t)$ represents the cumulative
amount of time spent on the $j$th activity up to time $t$.
For a control $T^r$ to be admissible, it must satisfy
additional properties which are specified below in Definition~\ref{t-adm-defn}.

\subsection*{State processes}

For a given scheduling policy $T^r$, the state processes of the
network are the associated $\II$-dimensional queue length process
$Q^r$ and the $\KK$-dimensional idle time process $I^r$. For each
$t\ge0$, $i \in\mathbb{I}$, $Q_i^r(t)$ represents the queue-length
at the $i$th buffer at time $t$ (including the jobs that are in
service at that time), and for $k=1,\ldots, \KK$, $I_k^r(t)$ is the
total amount of time the $k$th server has idled up to time $t$.
Let $q^r=Q^r(0) \in\mathbb{N}^{\II}$
be the
$\II$-dimensional vector of queue-lengths at time $0$. Note that, for
$j \in\mathbb{J}, t \ge0$, $S_j^r(T_j^r(t))$ is the total number
of services completed by the $j$th activity up to time $t$.
The total number of completed jobs (by activity $j$) up to time
$t$ that get rerouted to buffer $i$ equals $\bfijl
{i}{j}(S_j^r(T_j^r(t)))$. Recalling the definition of matrices $C$ and
$A$, the state of the
system at time $t \ge0$ can be described by the following equations:
%
\begin{eqnarray}
\hspace*{30pt} Q_i^r(t)&=& q^r+E_i^r(t)-\sum
_{j=1}^\JJ{C}_{i
j} S_j^r(T_j^r(t)) + \sum_{j=1}^\JJ{\bfijl{i}{j}}(S_j^r(T_j^r(t))),
\qquad i \in\mathbb{I}, \label{q-defn1}\\
I_k^r(t)&=&t-\sum_{j=1}^\JJ{A}_{k j}
T_j^r(t),\qquad k \in\mathbb{K}.\label{i-defn}
\end{eqnarray}
\subsection*{Heavy traffic}

We now describe the main heavy traffic assumption\cite{harri1,harri-canon}. We begin with a condition on the convergence of various parameters
in the sequence of networks $\{\mathcal{N}^r\}$.

\begin{assu}\label{assum-limit-param}
There are $q, \alpha, \sigma^{u} \in
{\R}^{\II}_+, \beta, \sigma^{v} \in{\R}^{\JJ}_+$, $ \theta_1
\in\R
^{\II}$, $\theta_2 \in\R^{\JJ}$
such that $ \beta>0, \sigma^v > 0$, $\alpha_i, \sigma^u_i=0$
if and only if $i \in\mathbb{I}\setminus\mathbb{I}'$, and, as $r
\to\infty$,
%
\begin{eqnarray}\label{assum-limit-param-1}
\theta_1^r &\Df& r({\alpha}^r-{\alpha})\rightarrow
\theta_1,\qquad \theta_2^r \Df r({\beta^r}-\beta) \rightarrow\theta_2,\nonumber
\\[-8pt]
\\[-8pt]
\nonumber
\sigma
^{u,r}&\rightarrow&\sigma^{u},\qquad
\sigma^{v,r}\rightarrow\sigma^{v}, \qquad\hat q^r \Df
\frac{q^r}{r} \rightarrow q.
\end{eqnarray}
\end{assu}

The definition of \textit{heavy traffic}, for the sequence
$\{\mathcal{N}^r\}$, as introduced in~\cite{harri1} (also see \cite
{Will-Bram-2work,bram-will-1,harri-canon}), is as follows.

\begin{defn}[{[Heavy traffic]}] \label{defn-ht} Define
$\II\times\JJ$ matrices $P', R$, such that $P'_{i j}\Df p^{j}_i$,
for $i \in\mathbb{I}, j \in\mathbb{J}$, and
%
\begin{equation}\label{R-defn} R
\Df(C-P')\operatorname{diag}(\beta).
\end{equation}
We say that the sequence
$\{\mathcal{N}^r\}$ approaches heavy traffic as $r \to\infty$ if, in
addition to Assumption \ref{assum-limit-param}, the following two
conditions hold:
\begin{longlist}[(ii)]
\item[(i)] There is a unique optimal solution $(x^*, \rho^*)$ to the
following linear program (LP):
%
\begin{equation}\label{LP}\qquad
\mbox{minimize } \rho\mbox{ such that } Rx = {\alpha}\mbox{ and }
Ax \le
\rho{\mathbf{1}}_{\KK} \qquad\mbox{for all }
x\ge0.
\end{equation}
\item[(ii)] The pair $(x^*, \rho^*)$
satisfies
%
\begin{equation}
\rho^*=1\quad \mbox{and}\quad Ax^*={\mathbf{1}}_{\KK}. \label{ab336}
\end{equation}
\end{longlist}
\end{defn}

\begin{assu}\label{assum-HT}
The sequence of networks $\{\mathcal{N}^r\}$ approaches heavy traffic
as $r \to
\infty$.
\end{assu}

\begin{remark}\label{rem-HT} From Assumption \ref{assum-HT},
$x^*$ given in (i) of Definition \ref{defn-ht} is the unique
$\JJ$-dimensional nonnegative vector satisfying
%
\begin{equation}\label{xstar-defn}
Rx^*= {\alpha}, \qquad Ax^*= {\mathbf{1}}_{\KK}.
\end{equation}
Following \cite{harri1}, assume without loss of generality
(by relabeling activities, if
necessary), that the first $\BB$ components of $x^*$ are
strictly positive (corresponding activities are referred to as
\textit{basic}) and the rest are zero (\textit{nonbasic} activities).
For later use, we partition the following matrices and vectors in
terms of \textit{basic} and \textit{nonbasic} components:
%
\begin{equation}\label{partition}
x^*= \left[
\matrix{
x^*_b \cr
\mathbf{0} }
\right],\quad
T^r= \left[
\matrix{
T^r_b\cr
T^r_n
}
\right],\quad A = [B:N], \quad R=[H:M],
\end{equation}
where $T^r$ is some control
policy, ${\mathbf{0}}$ is a
$(\JJ-\BB)$-dimensional vector of zeros, $B, N, H,M$ are $\KK\times
\BB$,
$\KK\times(\JJ-\BB)$, $\II\times\BB$ and $\II\times(\JJ- \BB)$
matrices, respectively.
\end{remark}
The following assumption (see \cite{Will-Bram-2work}) says that for
each buffer there is an associated basic activity.

\begin{assu}
\label{ab148}
For every $i \in\mathbb{I}$, there is a $j \in\mathbb{J}$ such that
$R_{ij} >0
$ and $x^*_j > 0$.
\end{assu}
\subsection*{Other processes}

Components of the vector $x^*$ defined above can be interpreted as the
nominal allocation rates
for the $\JJ$ activities.
Given a control policy~$T^r$, define the \textit{deviation process} $Y^r$
as the difference
between $T^r$ and the nominal allocation:
%
\begin{equation}\label{Y-defn}
Y^r(t) \Df x^*t - T^r(t),\qquad t\ge0.\vadjust{\goodbreak}
\end{equation}
It follows from (\ref{i-defn}) and (\ref{xstar-defn}) that
the idle-time process $I^r$ has the following representation:
\[
I^r(t) = A Y^r(t),\qquad t\ge0.
\]
Let $\mathbf{N}\Df\KK+\JJ-\BB$. Next we define a $\mathbf
{N}\times\JJ$ matrix
$K$ and
$\mathbf{N}$-dimensional process $U^r$ as follows:
%
\begin{equation}\label{U-K-defn}
K \Df\left[
\matrix{
B &N\cr
0 & -\III
}
\right],\qquad  U^r(t) \Df K Y^r(t),\qquad  t \ge0,
\end{equation}
where $\III$ denotes a $(\JJ-\BB)\times(\JJ-\BB)$ identity
matrix. Note that, with $T^r_n$ as in \eqref{partition},
%
\begin{equation}U^r(t) = \left[
\matrix{
I^r(t)\cr
T^r_n(t)}
\right],\qquad t \ge0. \label{flu1104}
\end{equation}
Finally, we introduce
the \textit{workload} process $W^r$ which is defined as a certain linear
transformation of the queue-length process and
is of dimension no greater than of the latter.
More precisely, $W^r$ is an $\LL$-dimensional process ($\LL= \II+\KK
-\BB$, see~\cite{Will-Bram-2work})
defined
as
%
\begin{equation}\label{workload-defn}
W^r(t) = \Lambda Q^r(t),\qquad  t\ge0,
\end{equation}
where $\Lambda$ is a $\LL\times\II$-dimensional matrix with rank
$\LL$
and nonnegative entries, called the
workload matrix. We will not give a complete description of $\Lambda$
since that requires additional notation; and we refer the reader
to \cite{Will-Bram-2work,harri-canon} for details. The key fact that
will be used in our analysis is that there is a $\LL\times\mathbf
{N}$ matrix
$G$ with nonnegative entries (see (3.11) and (3.12) in \cite
{harri-canon}) such that
%
\begin{equation}\label{mrgk}
\Lambda R = G K.
\end{equation}
We will impose the following additional assumption on $G$ which says that
each of its columns has at least one strictly positive entry. The
assumption is needed in the proof of Lemma
\ref{cont-bcp} [see \eqref{cond612}].
\begin{assu}\label{G-column}
There exists a $c>0$ such that for every $u \in
{\R}_+ ^{\mathbf{N}}$,
$|Gu| \ge c |u|$.
\end{assu}
\subsection*{Rescaled processes}

We now introduce two types of scalings. The first is the so-called
fluid scaling, corresponding to a law of large numbers, and the second
is the
standard diffusion
scaling, corresponding to a central limit theorem.

Fluid Scaled Process: This is obtained from the
original process by accelerating time by a factor of $r^2$ and
scaling down space by the same factor. The following fluid scaled
processes will
play a role in our analysis. For $t\ge0$,
%
\begin{eqnarray}\label{fl-scaled}
\bar E^r(t)& \Df& r^{-2}E ^r(r^2 t),\qquad \bar
S^r(t)\Df r^{-2}S^r(r^2 t), \nonumber\\
\bar{\bolds{\Phi}}^{r}(t) &\Df& r^{-2}{\bolds{\Phi}}^{r}(\lfloor r^2 t\rfloor
),\qquad
\bar T^r(t) \Df
r^{-2}T^r(r^2 t),\\
\bar I^r(t) &\Df&
r^{-2}I^r(r^2 t),\qquad \bar Q^r(t) \Df r^{-2}Q^r(r^2
t).\nonumber
\end{eqnarray}
Here for $x\in\R_+$, $\lfloor x\rfloor$ denotes its integer part, that
is, the greatest integer bounded by $x$.

Diffusion Scaled Process: This is obtained from the original
process by accelerating time by a factor of $r^2$ and, after
appropriate centering, scaling down
space by $r$. Some
diffusion scaled processes that will be used are as follows. For $t\ge0$,
%
\begin{eqnarray}\label{diff-scaled}
 \hat E^r(t) &\Df&\frac{(E^r(r^2
t)-{\alpha}^r
r^2 t
)}{r},\qquad \hat S^r(t) \Df\frac{(S^r(r^2 t)-{\beta^r}r^2 t
)}{r},\nonumber\\
\hat{\bolds{\Phi}}^{r}(t) &\Df&
\frac{({\bolds{\Phi}}^{r} (\lfloor r^2 t\rfloor])-\lfloor r^2 t\rfloor
P' )}{r},
\nonumber
\\[-8pt]
\\[-8pt]
\nonumber
\hat U^r(t)& \Df& r^{-1}U^r(r^2 t),\qquad
\hat Q^r(t)\Df r^{-1}Q^r(r^2
t),\\
\hat W^r(t) &\Df& r^{-1}W^r(r^2 t),\qquad \hat Y^r(t) \Df
r^{-1}Y^r(r^2 t).\nonumber
\end{eqnarray}
The processes $U^r, Q^r, W^r$ are not centered, as one finds (see Lemma
3.3 of \cite{BudGho2}) that, with any reasonable control policy, their
fluid scaled versions
converge to zero as $r\rightarrow\infty$. Define for $t \ge0$,
%
\begin{eqnarray}\label{xr-defn} \hat X^r_i(t) &\Df&\hat E^r_i(t) -
\sum_{j=1}^{\JJ} (C_{i j}-p_i^j) \hat S^r_j(\bar T^r_j(t))
\nonumber
\\[-8pt]
\\[-8pt]
\nonumber
&&{}-
\sum_{j=1}^{\JJ} \hat{\bolds{\Phi}}^{j,r}_i(\bar S^r_j(\bar T^r_j(t))),\qquad
i \in\mathbb{I}.
\end{eqnarray}
Recall $\theta_i^r$ and $\hat q^r$
from Assumption \ref{assum-limit-param}. Using (\ref{q-defn1}), (\ref
{i-defn}), \eqref{xstar-defn} and
(\ref{U-K-defn}), one has the following
relationships between the various scaled quantities defined above.
For all $t \ge0$,
\[
\hat Q^r(t)
= \hat\zeta^r(t) + R \hat Y^r(t), \qquad \hat U^r(t) = K
\hat Y^r(t),
\]
where
\begin{equation}\label{q-relation}
\hat\zeta^r(t)= \hat q^r + \hat X^r(t) + [\theta_1^r t - (C-P')
\operatorname{diag}(\theta_2^r) \bar
T^r(t)].
\end{equation}
Also, using
(\ref{workload-defn}), \eqref{mrgk} and \eqref{q-relation}, for all $t
\ge0$,
%
\begin{equation}
\quad\hat W^r(t) = \Lambda\hat q^r + \Lambda\hat X^r(t) + \Lambda
[\theta_1^r t - (C-P') \operatorname{diag}(\theta_2^r) \bar T^r(t)] +
G\hat U^r(t).\hspace*{-30pt} \label{w-relation}
\end{equation}

\subsection*{Admissibility of control policies}

The definition of admissible policies (Definition
\ref{t-adm-defn}), given below, incorporates appropriate
nonanticipativity requirements
and ensures feasibility by requiring that the
associated queue-length and idle-time processes ($Q^r, I^r$) are nonnegative.

For $m = (m_1,\ldots,m_\II) \in{\N}^{\II}, n = (n_1,\ldots,n_\JJ)
\in{\N}^{\JJ}$ we define the multiparameter filtration generated
by interarrival and\vadjust{\goodbreak} service times and routing variables as
%
\begin{eqnarray}
&&\bar{\mathcal{F}}^r((m,n))
\nonumber
\\[-8pt]
\\[-8pt]
\nonumber
&&\qquad = \sigma\{u_i^r(m'_i), v_j^r(n'_j),
\phi_{i}^{j,r}(n'_{j})\dvtx  m'_i\le m_i, n'_j\le n_j; i \in\mathbb{I}, j
\in\mathbb{J}\}.\hspace*{-20pt}
\end{eqnarray}
Then $\{ \bar{\mathcal{F}}^r((m,n))\dvtx  m \in{\N}^{\II}, n \in{\N
}^{\JJ}
\}$ is a multiparameter filtration with the following (partial)
ordering:
\[
(m^1,n^1) \le(m^2,n^2)\quad \mbox{if and only if}\quad
m^1_i \le m^2_i, n^1_j \le n^2_j; i \in\mathbb{I}, j \in\mathbb{J}.
\]
We refer the reader to Section 2.8 of \cite{Kurtz-redbook} for
basic definitions and properties of multiparameter filtrations,
stopping times and martingales.
Let
%
\begin{equation}\label{fr-sigmafield-defn} \bar{\mathcal{F}}^r \Df
\bigvee_{(m,n)\in{\N}^{\II+\JJ}} \bar{\mathcal{F}}^r((m,n)).
\end{equation}
For all $(m,n)\in\{0,1\}^{\II+\JJ}$, we define $
\bar{\mathcal{F}}^r((m,n)) = \bar{\mathcal{F}}^r(({\mathbf{1}},{\mathbf{1}}))$ where ${\mathbf{1}}$
denotes the vector of 1's.
It will be convenient to allow for extra randomness, than that captured
by $\bar{\mathcal{F}}^r$, in formulating the class of admissible policies.
Let $\mathcal{G}$ be a $\sigma$-field independent of $\bar{\mathcal
{F}}^r$. For $m
\in\mathbb{N}^{\II},
n \in\mathbb{N}^{\JJ}$, let
$ \mathcal{F}^r((m,n)) \equiv{\mathcal{F}}^r_{\mathcal{G}}(m,n) \Df
\bar{\mathcal{F}}^r((m,n))
\vee\mathcal{G}$.

\begin{defn}\label{t-adm-defn}
For a fixed $r$ and
$q^r \in\R_+^{\II}$, a scheduling policy $T^r=\{(T_1^r(t),
\ldots, T_\JJ^r(t))\dvtx  t\ge0 \}$ is called admissible for $\mathcal{N}^r$
with initial condition~$q^r$ if for some $\mathcal{G}$ independent of
$\bar
{\mathcal{F}}^r$, the following conditions hold:
\begin{longlist}[(iii)]
\item[(i)] $T^r_j$ is nondecreasing, nonnegative and satisfies
$T^r_j(0) = 0$ for $j \in\mathbb{J}$.\vspace*{1pt}
\item[(ii)] $I_k^r$ defined by (\ref{i-defn}) is nondecreasing,
nonnegative and satisfies $I_k^r(0) = 0$ for $k=1,\ldots, \KK$.
\item[(iii)] $Q_i^r$ defined in (\ref{q-defn1}) is nonnegative
for $i \in\mathbb{I}$.
\item[(iv)] Define for each $r, t\ge0$,
%
\begin{eqnarray}
\label{sigma0defn}{\sigma}^r_0(t)&= &(\sigma_0^{r, E}(t),
\sigma_0^{r, S}(t))
\nonumber
\\[-8pt]
\\[-8pt]
\nonumber
&\Df&\bigl(E_i^r(r^2 t)+1\dvtx  i \in\mathbb{I};
S_j^r(T_j^r(r^2t))+1\dvtx  j \in\mathbb{J}\bigr).
\end{eqnarray}
Then, for each $t\ge0$,
%
\begin{equation}\label{23025}{\sigma}^r_0(t)\mbox{ is a
} \{\mathcal{F}^r((m,n))\dvtx  m \in{\N}^{\II}, n \in{\N}^{\JJ} \}
\mbox{ stopping time}.
\end{equation}
Define the filtration $\{\mathcal{F}^r_1(t)\dvtx t\ge0\}$ as
%
\begin{eqnarray}\mathcal{F}^r_1(t) &\!\Df\!&
\mathcal{F}^r({\sigma}^r_0(t))
\nonumber
\\[-8pt]
\\[-8pt]
\nonumber
&\!=\!&\sigma\bigl\{A \in{\mathcal{F}}^r\dvtx  A
\cap\{{\sigma
}^r_0(t)\le
(m,n)\} \in\mathcal{F}^r((m,n)),
m \in{\N}^{\II}, n \in{\N}^{\JJ}\bigr\}.\hspace*{-35pt}
\end{eqnarray}
Then
%
\begin{equation}\label{23050} \hat U^r \mbox{ is }
\{\mathcal{F}^r_1(t)\}\mbox{-adapted}.
\end{equation}
\end{longlist}
Denote by $\mathcal{A}^r(q^r)$ the collection of all admissible
policies for
$\mathcal{N}^r$ with
initial condition $q^r$.
\end{defn}
\begin{remark} \label{rem921} (i) and (ii) in Definition \ref{t-adm-defn}
imply, in view of (\ref{i-defn})
and properties of the matrix $A$, that
%
\begin{equation}\qquad
0 \le T_j^r(t) -
T_j^r(s) \le t - s,\qquad j \in\mathbb{J}\mbox{ for all } 0 \le s
\le t < \infty. \label{lip1056}
\end{equation}
In particular, $T_j^r$ is a
process with Lipschitz continuous paths. Condition~(iv) in
Definition \ref{t-adm-defn} can be interpreted as a
nonanticipativity condition. Proposition 2.8 and Theorem 5.4 of \cite
{BudGho2} give
general sufficient conditions under which this property holds (see also
Proposition \ref{tisadmisnew} of the current work).
\end{remark}

\subsection*{Cost function}

For the network $\mathcal{N}^r$, we consider an expected infinite horizon
discounted (linear) holding cost associated with a scheduling
policy~$T^r$ and initial queue length vector $q^r$:
%
\begin{equation} \quad J^{r}(q^r,
T^r) \Df\IE\biggl(\int_0^{\infty} e^{-\gamma t} h\cdot\hat
Q^r(t) \,dt\biggr) + \IE\biggl(\int_0^{\infty} e^{-\gamma t} p \cdot
d\hat U^r(t)\biggr).\label{cost-r-th}
\end{equation}
Here, $\gamma\in(0,
\infty) $ is the ``discount factor'' and $h$, an $\II$-dimensional
vector with each component $h_i \in(0, \infty), i \in\mathbb{I}$, is
the vector of ``holding costs'' for the $\II$ buffers. In the second
term, $p\ge0$ is an $\mathbf{N}$-dimensional vector. The first~$\KK$
block of $U$ corresponds to the idleness process $I$, and, thus, the
second term in the cost, in particular, captures the idleness cost.
The last $\JJ-\BB$ components of $U$ correspond to the time spent
on \textit{nonbasic} activities. Thus, this formulation of the cost
allows, in addition to the idleness cost, the user to put a penalty
for using nonbasic activities.


The formulation of the cost function considered in our work goes back
to the original work of Harrison et al. \cite{harri1,harri2}.

The \textit{scheduling control problem} for $\mathcal{N}^r$ is to find an
admissible control policy~$T^r$ that minimizes
the cost $J^r$. The value function $V^r$ for this control problem is
defined as
%
\begin{equation}
\label{ab101}
V^r(q^r) \Df\inf_{T^r \in\mathcal{A}^r(q^r)} J^r(q^r, T^r),\qquad q^r \in
\mathbb{N}_0^{\II}.
\end{equation}

\subsection*{Brownian control problem}

The goal of this work is to characterize the limit of value functions
$V^r$ as $r\to\infty$, as the value function of a suitable diffusion
control problem.
In order to see the form of the diffusion control problem, we will
like to send $r\to\infty$ in \eqref{q-relation}.
Using the functional central limit theorem for renewal processes, it is
easily seen that, for all reasonable control policies (see again Lemma
3.3 of \cite{BudGho2}), when $\hat q^r$ converges to some $q \in\R
_+^{\II}$,
$\hat\zeta^r$ defined in (\ref{q-relation}) converges
weakly to
%
\begin{equation}\label{234b}\tilde\zeta= q + \tilde X + \theta
\curvi,
\end{equation}
where
%
\begin{equation}\label{theta-defn}
\theta\Df\theta_1-(C-P')
\operatorname{diag}(\theta_2) x^*.
\end{equation}
Here $\curvi(s) = s, s \ge0$ is the identity map and $\tilde X$ is a
Brownian motion with drift 0 and covariance
matrix
%
\begin{equation}\Sigma\Df\Sigma^u + (C-P')\Sigma^v \operatorname{diag}(x^*)
(C-P')' + \sum_{j=1}^{\JJ} {\beta_{j}}
x^*_j\Sigma^{\phi^j},\label{sigmadefn}
\end{equation}
where $\Sigma^u$ is a $\II\times\II$ diagonal matrix with diagonal
entries $(\sigma^u_i)^2, i \in\mathbb{I}$, $\Sigma^v$ is a
$\JJ\times\JJ$ diagonal matrix with diagonal entries
$(\sigma^v_j)^2, j \in\mathbb{J}$ and $\Sigma^{\phi^j}$s are
$\II\times\II$ matrices with entries $\sigma^{\phi^j}_{i_1 i_2},
i_1,i_2 \in\mathbb{I}$ [see (\ref{sigma-phi})]. Although the process~$\hat Y^r$
in~(\ref{q-relation}), for a general policy sequence $\{T^r\}$,
need not converge, upon
formally taking limit as $r \to\infty$, one is led to the following
diffusion control problem.

\begin{defn}[{[Brownian Control Problem (BCP)]}]\label{BCP}
A $\JJ$-dimensional adapted process $\tilde Y$, defined on some
filtered probability space $(\tilde\Omega, \tilde{\mathcal{F}},
\tilde
\IP,
\{\tilde\mathcal{F}(t)\})$ which supports an $\II$-dimensional
$\{\tilde\mathcal{F}(t)\}$-Brownian motion $\tilde X$ with drift 0
and covariance
matrix $\Sigma$ given by (\ref{sigmadefn}), is called an
admissible control for the Brownian control problem with the
initial condition $q \in\R_+^\II$ iff the following two properties
hold $\tilde
\IP$-a.s.:
%
\begin{eqnarray}
\tilde Q(t) &\Df& \tilde\zeta(t) + R \tilde Y(t) \ge
0 \qquad\mbox{where } \tilde\zeta
(t)=q + \tilde X(t) + \theta t, t \ge0, \label{ab440}\hspace*{-35pt}\\
\tilde U &\Df& K \tilde Y \mbox{ is
nondecreasing and } \tilde U(0) \ge0,\label{ab934}
\end{eqnarray}
where $\tilde\zeta$ and $\theta$ are as in (\ref{234b}) and (\ref
{theta-defn}) respectively.
We refer to $\Phi= (\tilde
\Omega, \tilde{\mathcal{F}},\allowbreak \tilde\IP, \{\tilde\mathcal{F}(t)\}
,\tilde X)$
as a
system.
We denote
the class of all such admissible controls by $\tilde\mathcal{A}(q)$. The
Brownian control problem is to
%
\begin{equation}\label{cost-BCP} \mbox{infimize } \tilde{J}(q,
\tilde Y) \Df\tilde
\E\Biggl[\int_0^{\infty} e^{-\gamma t} h\cdot\tilde Q(t) \,dt +
\int_{[0,{\infty})} e^{-\gamma t} p \cdot d \tilde U(t)\Biggr],\hspace*{-30pt}
\end{equation}
over all admissible controls $\tilde Y \in\tilde\mathcal{A}(q)$. Define
the value function
%
\begin{equation}\tilde J^*(q) = 
\inf_{\tilde Y \in\tilde\mathcal{A}(q)} \tilde{J}(q, \tilde
Y).\label{mincost-BCP}
\end{equation}
\end{defn}

Recall our standing assumptions \eqref{stpos1029}, \eqref{ui1032},
\eqref{ass1120}, Assumptions \ref{assum-limit-param}, \ref{assum-HT},
\ref{ab148} and \ref{G-column}. The following is the main result of
\cite{BudGho2}.

\begin{theorem}[(Budhiraja and Ghosh \cite{BudGho2}, Theorem 3.1,
Corollary 3.2)]
\label{ab937}
Fix $q \in\R_+^{\II}$ and for $r > 0$, $q^r \in\mathbb{N}^{\II}$
such that
$\hat q^r \to q$ as $r \to\infty$. Then
\[
\liminf_{r \to\infty}
V^r(q^r) \ge\tilde J^*(q).
\]
\end{theorem}

\begin{remark}
The proof in \cite{BudGho2} is presented for the case where in the definition
of $V^r(q^r)$ [see \eqref{ab101}], $\mathcal{A}^r(q^r)$ is replaced
by the
smaller family~$\bar\mathcal{A}^r(q^r)$ which consists of all $T^r
\in\mathcal{A}
^r(q^r)$ that satisfy
(iv) of Definition \ref{t-adm-defn} with $\mathcal{F}^r((m,n))$
replaced by
$\bar{\mathcal{F}}^r((m,n))$. Proof for the slightly more general setting
considered in the
current paper requires only minor modifications and, thus, we omit the details.
\end{remark}

For the main result of this work, we will need additional assumptions.

\begin{assu}
\label{non-deg}
The matrix $\Sigma$ is positive definite.
\end{assu}
We will make the following assumption on the probabilities of
deviations from the mean for the underlying
renewal processes. Similar conditions have been used in previous works
on construction of asymptotically optimal control policies
\cite{bellwill,bellwill2,BudGho,ata-kumar,dai-lin}.

\begin{assu}\label{ldp}
There exists $\curvm> 2$ and, for each $\delta> 0$, some $\varsigma
(\delta) \in(0, \infty)$ such that,
for $j \in\mathbb{J}$, $i \in\mathbb{I}$, $r \ge1$, $t \in
(1,\infty)$,
\begin{eqnarray*}
\IP\bigl(|S_j^r(t)-\beta_j^rt| \ge\delta t\bigr) &\le&\frac{\varsigma(\delta
)}{t^{\curvm}},
\\
\IP\bigl(|E_i^r(t)- \alpha_i^r t| \ge\delta t\bigr) &\le&\frac{\varsigma
(\delta)}{t^{\curvm}},
\\
\IP\bigl(|\Phi_i^{j,r}(S_j^r(t))- p^j_i\beta_j^r t| \ge\delta t\bigr) &\le&
\frac{\varsigma(\delta)}{t^{\curvm}}.
\end{eqnarray*}
\end{assu}

The third inequality above is a consequence of the first two, but we
note it explicitly here for future use.
The assumption is clearly satisfied when $E_i$ and~$S_j$ are Poisson
processes. For general renewal processes, such inequalities
hold under suitable moment conditions on the interarrival and service
time distributions. Indeed, if for some $\curvm> 1$
%
\begin{eqnarray}
\E\Bigl[ \sup_{r} [u_i^r (1)]^{2\curvm} \Bigr]&<& \infty,
\nonumber
\\[-8pt]
\\[-8pt]
\nonumber
\E
\Bigl[ \sup_{r} [v_j^r(1)]^{2\curvm}\Bigr] &<& \infty \qquad\mbox{for all
} i \in\mathbb{I}, j \in\mathbb{J},
\end{eqnarray}
then, from Theorem 4 of \cite{K-T}, Assumption \ref{ldp} is
satisfied.

We now introduce an assumption on the regularity properties of a
certain Skorohod map. This map plays a crucial role in our analysis;
see proofs
of Theorems~\ref{main-bcp} and~\ref{mainweak} [see in particular,
\eqref
{110}, \eqref{ab302}, discussion below~\eqref{ab419} and the proof of
Theorem \ref{newmain518}].
Let $\mathcal{D}_+^{\II} \Df\{x \in\mathcal{D}^{\II}\dvtx  x(0) \ge0
\}$ and
%
\begin{equation}
\label{refmat}
D= (C - P') \operatorname{diag} (\beta) \operatorname{diag}(x^*)C'.
\end{equation}

\begin{definition}
\label{ab948}
Given $x \in\mathcal{D}_+^{\II}$, we say $(z,y) \in\mathcal
{D}^{\II} \times
\mathcal{D}
^{\II}$ solve the \mbox{Skorohod} Problem (SP)
for $(x, D)$ if: (i) $z(0) = x(0)$, (ii) $z = x + Dy$,
(iii) $y$ is nondecreasing and $y(0) \ge0$,
(iv) $z(t) \ge0$ for all $t \ge0$, and\break (v) $\int_{[0, \infty)} 1_{\{
x_i(t) > 0\}}\, dy_i(t) = 0$ for all $i \in\II$.
\end{definition}

Denoting by $\mathcal{D}_0$ the set of $x \in\mathcal{D}_+^{\II}$
such that
there is a unique solution to the SP for $(x,D)$,
we define maps $\Gamma\dvtx  \mathcal{D}_0 \to\mathcal{D}^{\II}$,
$\hat\Gamma\dvtx  \mathcal{D}_0 \to\mathcal{D}^{\II}$ as
$\Gamma(x) = z$, $\hat\Gamma(x) = y$ if $(z,y)$ solve the
$SP$ for $(x, D)$.

\begin{assu}\label{ab1012}
$\mathcal{D}_0 = \mathcal{D}_+^{\II}$ and the maps $\Gamma, \hat
\Gamma$
are Lipschitz, namely, there exists $L \in(0, \infty)$
such that for all $x_1, x_2 \in\mathcal{D}_+^{\II}$,
\[
\sup_{0 \le t < \infty}
\{ |\Gamma(x_1)(t) - \Gamma(x_2)(t)| +
|\hat\Gamma(x_1)(t) - \hat\Gamma(x_2)(t)| \}
\le L \sup_{0 \le t < \infty}|x_1(t) - x_2(t)|.
\]
\end{assu}
We refer the reader to \cite{DupuisIshii1991,Dup-Ram} and
\cite{har-rei} for sufficient conditions under which the
above regularity property of the Skorohod map holds. See also Example~\ref{exammodel} below.
For later use we introduce the notation $\bar\Gamma(x) =
\operatorname{diag}(x^*)C'\hat\Gamma(x)$ for $x \in\mathcal
{D}_+^{\II}$. Since~$A$ has
nonnegative entries and $x^*_j = 0$ for
$j = \BB+1, \ldots,\JJ$, we see from the definition of $K$ [see
\eqref
{U-K-defn}] that
%
\begin{equation}
\label{abn756}
 \mbox{if } y = \bar\Gamma(x), \mbox{ then } Ky \ge
0.
\end{equation}
%
For rest of the paper, in addition to the assumptions listed above
Theorem~\ref{ab937}, Assumptions \ref{non-deg}, \ref{ldp} and
\ref{ab1012} will be in force.
The
main result of the paper is the following.

\begin{theorem}
\label{main1016}
Fix $q \in\R_+^{\II}$. Let for $r > 0$, $q^r \in\mathbb{N}_0^{\II}$
be such that $\hat q^r \to q$ as $r \to\infty$.
Then
\[
\limsup_{r \to\infty} V^r(q^r) \le\tilde J^*(q).
\]
\end{theorem}

The theorem is an immediate consequence of Theorem \ref{jrtoj} below
which is proved in Section \ref{secproof}.
For $\varepsilon> 0$, we say
$Y \in\tilde\mathcal{A}(q)$ is $\varepsilon$-optimal for the BCP with
initial value
$q$ if
\[
\tilde{J}(q, Y)\leq\tilde{J}^*(q)+\varepsilon.
\]
When clear from the context, we will omit the phrase ``for the BCP with
initial value $q$'' and merely say that
$Y$ is $\varepsilon$-optimal.

\begin{theorem}
\label{jrtoj}
Fix $q \in\R_+^{\II}$. For $r > 0$, let $q^r \in\mathbb{N}_0^{\II}$
be such that $\hat q^r \to q$ as $r \to\infty$. For every $\varepsilon> 0$,
there exists
$\tilde Y \in\tilde\mathcal{A}(q)$, which is $\varepsilon$-optimal,
and a~sequence $T^r
\in\mathcal{A}^r(q^r)$,
$r \ge1$ such that
\[
J^r(q^r, T^r) \to\tilde J(q, \tilde Y)\qquad \mbox{as } r \to\infty.
\]
\end{theorem}


Combining Theorems \ref{ab937} and \ref{main1016}, the following
is\vadjust{\goodbreak}
immediate.
\begin{theorem}
\label{maincorr}
Fix $q \in\R_+^{\II}$. For $r > 0$, let $q^r \in\mathbb{N}_0^{\II}$
be such that $\hat q^r \to q$ as $r \to\infty$.
Then, as $r \to\infty$, $V^r(q^r) \to\tilde J^*(q)$.
\end{theorem}

Assumptions made in this work can loosely be divided into two
categories: Assumptions on the underlying stochastic primitives, which
include, in particular, the heavy traffic conditions
(Assumptions \ref{assum-limit-param}, \ref{assum-HT}, \ref{non-deg}
and~\ref{ldp}), and assumptions made on the network structure (Assumptions
\ref{ab148}, \ref{G-column} and~\ref{ab1012}).
Below we discuss the validity of these structural assumptions for some
basic families of SPN models.

\begin{example} \label{exammodel}
The following examples have been described in detail in~\cite
{Will-Bram-2work}. We will assume here, without loss of generality,
that $\beta_j >0$ for all $j \in\mathbb{J}$ (an activity~$j$ for which
$\beta_j=0$ can simply be
deleted from the network description). Furthermore, for all three
settings considered below, Assumption~\ref{ab148} can be made without
loss of generality, since otherwise one can consider a~reduced system
obtained by omitting the buffers that are not processed by any basic activity.
Assumption \ref{G-column} states that the matrix $G$ can be chosen
in a manner such that it has no columns that are identically zero.
Roughly speaking, it says that a nonzero control action leads to a
nonzero state displacement.
Although this appears to be a very natural geometric condition and is
trivially satisfied for networks in part (a) below,
it is not clear that it holds always for examples in parts (b) and
(c) below. We will assume this condition to hold without further
comment.

Thus, in discussion below, we will focus only on Assumption \ref{ab1012}.

%
\begin{figure}[b]

\includegraphics{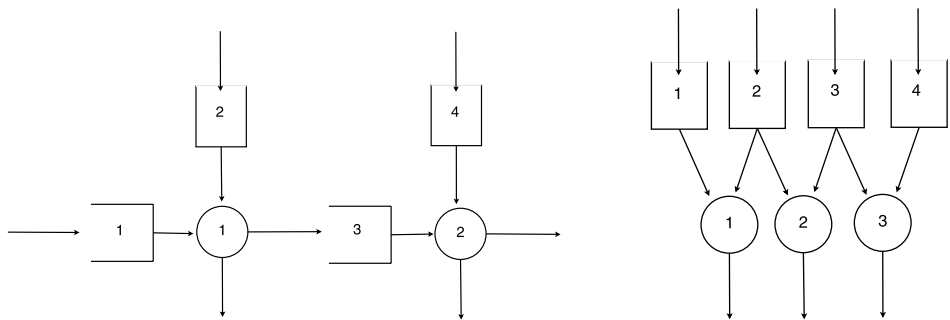}

\caption{Open multiclass network (left) and parallel-server system
(right).} \label{fig-mult}
\end{figure}

\begin{longlist}
\item[(a)] \textit{Open multiclass queueing networks}: These correspond
to a setting where each buffer is processed by exactly one
activity and, consequently, there is a one-to-one correspondence
between activities and buffers, that is, $\JJ=\II$
(see left figure in Figure \ref{fig-mult} for an example).
For such networks, $R$ is an $\II\times\II$-matrix of the form $R =
(\III-P')\operatorname{diag}(\beta)$ where $P$ is a nonnegative matrix
with spectral radius less than 1. In particular, $R$ is nonsingular,
$K$ is a $\KK\times\JJ$ matrix with full row rank and one can take
$\Lambda= KR^{-1}$ and $G= \III$. Here $D=(\III-P')\operatorname
{diag}(\beta
)\operatorname{diag}
(x^*)$ and from \cite{har-rei} it is known that for such~$D$ Assumption
\ref{ab1012} is
satisfied.

\item[(b)] \textit{Parallel server networks}: For such SPN, a buffer
can be served by more than one activity, however, each job gets
processed exactly once before leaving the system (i.e., there is no
rerouting). See right figure in Figure \ref{fig-mult} for an
example. In particular, $P=0$ and, hence, $R=C \operatorname
{diag}(\beta)$.
%
In this case, $D= C \operatorname{diag}(\beta)\operatorname
{diag}(x^*) C' \equiv\operatorname{diag}(\gamma^*)$,
where $\gamma^*_i = \sum_{j\dvtx  \sigma_1(j)=i} \beta_j x^*_j$, for $i
\in
\mathbb{I}$.
From Assumption \ref{ab148} (which, as was noted above, can be made
without loss of generality) we have that $\gamma^*
> 0$ and, thus, $D$ is a diagonal matrix with strictly positive
diagonal entries. Assumption \ref{ab1012} is clearly satisfied for
such matrices.

%
\begin{figure}

\includegraphics{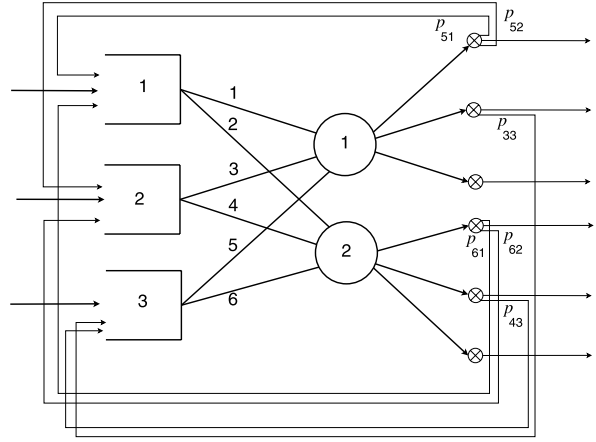}

\caption{A Job-shop network. Here $p_{61}= p_{51}$, $p_{62}=p_{52}$,
$p_{33}=p_{43}$.} \label{fig-jshop}
\vspace*{3pt}
\end{figure}


\item[(c)] \textit{Job-shop networks}: This subclass of networks
combines features of both (a) and (b): A buffer can be processed
by more than one activity and jobs, once served, can get rerouted to
another buffer for additional processing.
See Figure \ref{fig-jshop} for an example.
Following specific examples
considered in~\cite{Kumar} (see also \cite{Will-Bram-2work}), we define
job-shop networks as those which satisfy the following property: If for
some $j,j' \in\mathbb{J}$, and $i \in\mathbb{I}$, $\sigma_1(j)
=\sigma_1(j') =
i$, then $p_{j i'} = p_{j' i'}$ for all $i' \in\mathbb{I}$. Namely, jobs
corresponding to any two activities that process the same buffer $i$
have an identical (probabilistic) routing structure, following their
completion by the respective servers.
It is easily checked that in this case $D=(\III-\tilde
{P'})\operatorname{diag}
(\gamma
^*)$, where $\gamma^*$ is as introduced in (b) and $\tilde{P}$ is an
$\II\times\II$-dimensional matrix with entries $\tilde p_{i, i'}= p_{j,
i'}$ where $j \in\mathbb{J}$ is such that $\sigma_1(j)=i$.
Under the condition that $\tilde P$ has spectral radius less than 1, it
follows from \cite{har-rei} that Assumption \ref{ab1012} is satisfied.
\end{longlist}
\end{example}

\section{Near-optimal controls for BCP}\label{nearopt}

The rest of the paper is devoted to the proof of Theorem \ref{jrtoj}.
Toward that goal, in this section we construct
near-optimal
controls for the BCP with certain desirable features. This construction
is achieved in Theorem \ref{newmain518},
which is the main result of this section.

Since an admissible control is not required to be of bounded
variation, the BCP is a somewhat nonstandard diffusion control problem
and is
difficult to analyze directly. However, as shown in \cite{harri1},
under assumptions made in this paper, one can replace this control
problem by
an equivalent problem of Singular
Control with State Constraints (SCSC). This control problem, also
referred to as the Equivalent Workload
Formulation (EWF) of the BCP, is given below. We begin by introducing
the cost function, that is, optimized in this equivalent control problem.

\textit{Effective cost function}: Recall the
definition of the workload matrix $\Lambda$ introduced in
(\ref{workload-defn}). Let $\mathcal{W} \Df\{\Lambda z\dvtx  z \in
{\R}_{+}^{\II}\}$. For each $w \in\mathcal{W}$, define
%
\begin{equation}
\hat h(w) \Df\inf\{h\cdot q \dvtx  \Lambda q = w, q \ge
0\}.\label{h-hat-defn}
\end{equation}
Since $h > 0$, the infimum is attained for
all $w \in\mathcal{W}$. It is well known (see Theorem~2 of \cite{bohm})
that one can take a continuous selection of the minimizer in the
above linear program. That is, there is a continuous map $\tilde q^*
\dvtx  \mathcal{W} \to\R_+^{\II}$ such that
%
\begin{equation}
\tilde q^*(w) \in
\mathop{\operatorname{argmin}}_q \{ h\cdot q\dvtx  \Lambda q = w, q \ge0
\}.\label{z-star-defn}
\end{equation}
Thus, in particular, $\hat h$ is continuous. One can check that $\hat
h$ satisfies linear lower and upper bounds.
In order to see this,
define 
%
\begin{equation}\label{32half} q^*(w) = \tilde q^*(w)1_{\{|w|\le1\}}
+ |w| \tilde
q^*\biggl(\frac{w}{|w|}\biggr)1_{\{|w| > 1\}}.
\end{equation}
Then \eqref{z-star-defn} holds with $\tilde q^*$ replaced by $q^*$.
Since $\hat h(w) = h\cdot q^*(w)$ and $h > 0$, we have from the above
display that
%
\begin{equation}\label{cond511} b_1 |w| - b_2 \le|\hat h(w)| \le
b_3(1 + |w|),\qquad w \in\mathcal{W}
\end{equation}
for some $b_1, b_2, b_3 \in(0, \infty)$. Also, uniform continuity of $q^*$
on \mbox{$\{w \in\mathcal{W}\dvtx  |w| \le1\}$} shows that
%
\begin{eqnarray}
\label{cond512}
|\hat h(w_1) - \hat h(w_2)|
\le\hat m (\delta) (1 + |w_1| + |w_2|),
\nonumber
\\[-8pt]
\\[-8pt]
\eqntext{w_1, w_2 \in\mathcal{W},
|w_1 - w_2| \le\delta.}
\end{eqnarray}
Here $\hat m$ is a modulus, that is, a nondecreasing function from $[0,
\infty) \to[0, \infty)$
satisfying $\hat m(0+) = 0$.
Inequalities \eqref{cond511} and \eqref{cond512} will be used in order
to appeal to some results from \cite{AtBu,BuRo} (see Remark \ref
{rem32} below).
Define
%
\begin{equation}\mathcal{K} \Df\{u \in{\R}^{\mathbf{N}}
| u = K y, y \in{\R}^{\JJ}\}.
\end{equation}
The Equivalent
Workload Formulation (EWF) and the associated control problem are
defined as follows.

\begin{defn}[{[Equivalent Workload Formulation (EWF)]}]\label{EWF}
An $\mathbf{N}$-di\-mensional adapted process $ \tilde U$, defined on some
filtered probability space $(\tilde\Omega,\allowbreak
\tilde{\mathcal{F}},
\tilde
\IP, \{\tilde
\mathcal{F}(t)\})$ which supports an $\II$-dimensional
$\{\mathcal{F}(t)\}$-Brownian motion~$\tilde X$ with drift 0 and covariance
matrix $\Sigma$ defined in (\ref{sigmadefn}), is called an
admissible control for the EWF with initial condition $w \in\mathcal
{W} $ iff
the following two properties hold $\tilde\IP$-a.s.:
%
\begin{eqnarray}\label{EWF-1}
&&\mbox{$\tilde U$ is nondecreasing}, \qquad\tilde U(0)\ge0,\qquad
\tilde U(t) \in\mathcal{K}\qquad \mbox{for all $t \ge0$,}
\nonumber\hspace*{-30pt}
\\[-8pt]
\\[-8pt]
\nonumber
&&\tilde
 W(t) \Df w + \Lambda\theta t+ \Lambda\tilde X(t) + G \tilde U(t) \in
\mathcal{W} \qquad \mbox{for all } t \ge0,\hspace*{-30pt}
\end{eqnarray}
where $\theta$ is as in (\ref{theta-defn}). We denote the class of
all such admissible controls by $\tilde\mathcal{A}_0 (w)$. The control
problem for the EWF is to
%
\begin{eqnarray}\label{cost-EWF}
\qquad\mbox{infimize } \tilde
{J}_0(w,\tilde U) \Df
\tilde\E\int_0^{\infty} e^{-\gamma t} \hat h(\tilde W(t)) \,dt +
\tilde\E
\int_{[0,{\infty})} e^{-\gamma t} p \cdot d\tilde U(t),
\end{eqnarray}
over all admissible controls $\tilde U \in\tilde\mathcal{A}_0 (w)$.
Define the value function
%
\begin{equation}
\tilde J_0^*(w) = 
\inf_{\tilde U \in\tilde\mathcal{A}_0 (w)} \tilde{J_0}(w,
\tilde U).\label{mincost-EWF}
\end{equation}
\end{defn}
From Theorem 2 of \cite{harri1} it follows that for all
$w\in
\mathcal{W}, q\in{\R}_+^{\II}$ satisfying $w=\Lambda q$,
%
\begin{equation}\label{mincost-equal}
\tilde J^*(q)=\tilde J_0^*(w).
\end{equation}
%
The following lemma will be used in order to appeal to some results
from~\cite{AtBu,BuRo}. The proof is based on arguments in \cite
{Will-Bram-2work}.
Let $\mathcal{K}_+ = \mathcal{K}\cap\R^{\mathbf{N}}_+$.
\begin{lemma}
\label{intgood}
The cones $\mathcal{K}_+$ and $G\mathcal{K}_+$ have nonempty
interiors and $\mathcal{W}^o
\cap G\mathcal{K}_+ \neq\varnothing$.
\end{lemma}
\begin{pf}
From \cite{Will-Bram-2work} (see above Corollary 7.4 therein) it
follows that $H$ has full row rank and so there is a $\BB\times\II$
matrix $H^{\dag}$ such that
$HH^{\dag} = \III$. Let $x_{b} =H^{\dag}{\mathbf{1}}_{\II}$ and let
$\varepsilon_0
\in(0, \infty)$ be
sufficiently small such that for all $\varepsilon\in(0, \varepsilon_0]$,
$x^{\varepsilon}_b = x^*_b + \varepsilon x_b > 0$. Let $\alpha_0 = \alpha+
\varepsilon
_0{\mathbf{1}}_{\II}$ and $x_0 = [x^{\varepsilon_0}_b, {\mathbf{0}}]' \in\R^{\JJ}$.
Then $Rx_0 = \alpha_0$. We will now argue that
$\vartheta= \Lambda\alpha_0 \in\mathcal{W}^o \cap(G\mathcal{K}_+)^o$.
Since the rows of $\Lambda$ are linearly independent, we can find an
$\II\times\LL$ matrix $\Lambda^\dag$ such that $\Lambda\Lambda
^\dag
= \III$.
Fix $\delta= \frac{\varepsilon}{2|\Lambda^\dag|}$. Then, whenever
$\tilde
\vartheta\in\R^{\LL}$, $|\tilde\vartheta| < \delta$, we have
$\alpha_0 + \Lambda^\dag\tilde\vartheta\in\R^{\III}_+$ and so
$\vartheta+ \tilde\vartheta= \Lambda(\alpha_0 + \Lambda^\dag
\tilde
\vartheta) \in\mathcal{W}$.
This shows that $\vartheta\in\mathcal{W}^o$. Next note that
$\vartheta=
\Lambda\alpha_0 = \Lambda Rx_0 = GKx_0.$
Since $Kx_0 = Bx^{\varepsilon_0}_b$ and $x^{\varepsilon_0}_b > 0$, we have
that $Kx_0
\in\mathcal{K}_+$ and so $\vartheta= GKx_0 \in G\mathcal{K}_+$.
Since $x^{\varepsilon_0}_b > 0$, we can find $\varepsilon_1 \in(0, \infty
)$ such
that whenever $\tilde x_b \in\R^{\BB}$ is such that
$|\tilde x_b| \le\varepsilon_1$, $x^{\varepsilon_0}_b + \tilde x_b > 0$.
Now fix
$\delta_1 = \frac{\varepsilon_1}{|H^{\dag}||\Lambda^\dag|}$. Then,
for any
$\tilde\vartheta\in
\R^{\LL}$ with $|\tilde\vartheta| \le\delta_1$ and \mbox{$\tilde x =
[H^\dag
\Lambda^\dag\tilde\vartheta, {\mathbf{0}}]' \in\R^{\JJ}$,}
\[
\vartheta+ \tilde\vartheta= \Lambda(\alpha_0+ \Lambda^\dag\tilde
\vartheta) = \Lambda( Rx_0 + HH^\dag\Lambda^\dag\tilde\vartheta)
= \Lambda(Rx_0 +R \tilde x) = G\bigl(K(x_0+\tilde x)\bigr).
\]
Since $|H^\dag\Lambda^\dag\tilde\vartheta| \le\varepsilon_1$, we have
$(x_0+\tilde x)_j > 0$ for all $j = 1, \ldots, \BB$. Also,
$ (x_0+\tilde x)_j = 0$ for all $j = \BB+1, \ldots, \JJ$.
Thus, $K(x_0+ \tilde x) \in\mathcal{K}_+$ and, therefore, $\vartheta+
\tilde
\vartheta\in G\mathcal{K}_+$. It follows that $\vartheta\in
(G\mathcal{K}_+)^o$.

Finally, we show that $\mathcal{K}_+^o \neq\varnothing$.
Since $Ax^*= {\mathbf{1}}_{\KK}$, every row of $B$ must contain at least
one strictly positive entry. Thus, $Bx^{\varepsilon_0}_b > 0$. Choose
$\varepsilon_2
\in(0, \infty)$ sufficiently small such that
$Bx^{\varepsilon_0}_b - \varepsilon_2 N {\mathbf{1}}_{\JJ-\BB} > 0$. Let $\bar
x = [x^{\varepsilon
_0}_b, -\varepsilon_2{\mathbf{1}}_{\JJ-\BB}]' \in\R^{\JJ}$. We now argue that
$\bar u = K \bar x \in(\mathcal{K}_+)^o$. Note that by construction
$\bar u >
0$. Thus, we can find $\varepsilon_3 > 0$ such that $\bar u + \tilde u
\ge0$
whenever $\tilde u \in\R^{\mathbf{N}}$ satisfies $|\tilde u| \le
\varepsilon_3$.
Also, since $K$ has full row rank (see Corollary 6.2 of \cite
{Will-Bram-2work}), we can find a $\JJ\times\mathbf{N}$ matrix
$K^\dag$
such that
$KK^\dag= \III$. Thus, $\bar u+ \tilde u = K(\bar x + K^\dag\tilde u)$
and, consequently, $\bar u + \tilde u \in\mathcal{K}_+$. The result follows.
\end{pf}

Note that the vector $x_0$ constructed in the proof of the lemma above
has the property that $Rx_0 > 0$ and $Kx_0 \ge0$. Thus, we have shown
the following:

\begin{cor}
\label{tset}
The set $\mathbb{T} = \{y \in\R^{\JJ}\dvtx  Ky \ge0, Ry > 0 \}$ is nonempty.
\end{cor}

The above result will be used in the construction of a suitable near
optimal control policy for the BCP [see below \eqref{111}].

\begin{remark}\label{rem32}
We will make use of some results from \cite{AtBu} and \cite{BuRo} that
concern a general
family of singular control problems with state constraints. We note
below some properties of the model
studied in the current paper that ensure that the assumptions of \cite
{AtBu} and \cite{BuRo} are satisfied:
\begin{longlist}[(a)]
\item[(a)] $G$ has full row rank. This follows from the observation
that $K, \Lambda$ and~$R$ have full row ranks and, therefore,
%
\[
\operatorname{rank}(G) = \operatorname{rank}(GK) = \operatorname{rank}(\Lambda R) = \operatorname
{rank}(\Lambda) = \LL.
\]
\item[(b)] $\mathcal{W}^o \cap(G\mathcal{K}_+)^o \neq\varnothing$
and $\mathcal{K}_+$
has a
nonempty interior (see Lemma \ref{intgood}).
\item[(c)] $(Gu) \cdot{\mathbf{1}}_{\LL} \ge|Gu|$, $u \cdot{\mathbf{1}}_{\mathbf{N}
} \ge|u|$ for all
$u \in\mathcal{K}_+$ and $w \cdot{\mathbf{1}}_{\LL} \ge|w|$ for all
$w\in\mathcal{W}$.
This is an immediate consequence of the fact that the entries of $G$
and $\Lambda$ are nonnegative [see above \eqref{mrgk}].
\item[(d)] Since $\Lambda$ has full row rank and, by Assumption \ref
{non-deg},
$\Sigma$ is positive definite, we have that $\Lambda\Sigma\Lambda'$
is positive definite.
\end{longlist}
The above properties along with Assumption \ref{G-column}, \eqref
{cond511} and \eqref{cond512} ensure that
Assumptions of \cite{AtBu} and \cite{BuRo} are satisfied in our
setting. In particular,
Assumption~(2.1)--(2.2) and (2.8)--(2.10) of \cite{AtBu} hold in view of
properties (b), (c) and (d) and equations \eqref{cond511} and \eqref
{cond512}. Similarly, Assumptions (1), (5)
and~2.2 of \cite{BuRo} hold in our setting [from property (c), \eqref
{cond511} and Assumption~\ref{G-column}, resp.].
Henceforth, when appealing to results from \cite{AtBu} and \cite
{BuRo}, we will
not make an explicit reference to these conditions.
\end{remark}

Recall $\tilde\zeta(t)$ and the map $\bar\Gamma$
introduced above (\ref{234b}) and \eqref{abn756}, respectively. The
following is a key step in the construction of
a near-optimal control with desirable properties.\vadjust{\goodbreak}

\begin{theorem}\label{main-bcp} Fix $q \in\R_+^{\II}$.
For each $\varepsilon\in(0, \infty)$, there exists
\mbox{$\tilde{Y}^{(1)} \in\tilde\mathcal{A}(q)$}, given on some system
$\Phi$, that
is $\varepsilon$-optimal and has the following properties:
%
\begin{equation}\label{ab315old} \tilde{Y}^{(1)} = \tilde{Y}_0^{(1)} +
\bar\Gamma\bigl(\tilde\zeta+
R\tilde{Y}_0^{(1)}\bigr),
\end{equation}
where $\tilde Y_0^{(1)}$ is an adapted
process with sample paths in $\mathcal{D}^{\JJ}$ satisfying the following:
For some $T, \eta, M \in(0, \infty)$, $p_0, j_0 \in\mathbb{N}$,
with $\theta= T/p_0$ and \mbox{$\kappa= \theta/j_0$},
\begin{enumerate}[(iii)]
\item[(i)] $\tilde{Y}_0^{(1)}(t) = \tilde{Y}_0^{(1)}(n\theta)$
for $t \in[n\theta, (n+1)\theta), n=0,1,\ldots, p_0-1$ and
$\tilde{Y}_0^{(1)}(t) = \tilde{Y}_0^{(1)}(p_0\theta)$ for $t \ge T =
p_0\theta$.

\item[(ii)] Letting $\mathcal{S}_M^{\eta}= \{b\eta\dvtx  b \in\mathbb
{Z}^{\JJ},
|b|\eta\le M, Kb \ge0\}$,
\[
\partial
\tilde{Y}_0^{(1)}(n) \Df\tilde{Y}_0^{(1)}(n\theta)-\tilde
{Y}_0^{(1)}\bigl((n-1)\theta\bigr) \in\mathcal{S}_M^{\eta},
\]
for $n=1,\ldots, p_0$
and $\partial\tilde{Y}_0^{(1)}(0) \Df\tilde{Y}_0^{(1)}(0) = 0$.

\item[(iii)] There is
an i.i.d sequence of Uniform (over $[0,1]$) random variables~$\{\tilde
{\mathcal{U}}_n\}$, that is, independent of $\tilde\zeta$, and
for each $n = 1, \ldots, p_0$
a measurable map $\varpi_n\dvtx  \R^{\II nj_0} \times[0, 1] \to
\mathcal{S}_M^{\eta}$, $n = 1, \ldots, p_0$, such that the map $x
\mapsto\varpi
_n(x, t)$ is
continuous, for a.e. $t$ in $[0,
1]$, and
%
\begin{equation}\label{ab502}
\partial\tilde{Y}_0^{(1)}(n) = \varpi_n(\mathcal{X}^\kappa(n),
\tilde
{\mathcal{U}}_n),
\end{equation}
where $\mathcal{X}^\kappa(n)=\{\tilde{X}(\ell\kappa)\dvtx  \ell
=1,\ldots,
nj_0\}$,
$n=0, 1, \ldots, p_0$.
\end{enumerate}
\end{theorem}

Proof of Theorem \ref{main-bcp} is given in Section \ref{pfmainbcp}.

\begin{remark}
\label{on3}
The above theorem provides an $\varepsilon$-optimal control $\tilde Y^{(1)}$,
which is the ``constrained''-version of a piecewise constant
process $\tilde Y^{(1)}_0$. The value of $\tilde Y^{(1)}_0$ changes
only at
time-points that are integer multiples of $\theta$ and is constant for
$t > T= p_0\theta$. Also, the changes in (the value of) the process
occur in jumps with sizes that are integer multiples of some $\eta>0$
and are bounded by $M$.
The third property in the theorem plays an important role in the weak
convergence proof [Theorem~\ref{mainweak}, see, e.g.,
\eqref{ab459n}] and says that the jump-sizes of this piecewise
constant process are determined by the Brownian motion~$\tilde X$ sampled
at discrete instants $\{\kappa, 2 \kappa,\ldots\}$
and the independent random variable $\tilde{\mathcal{U}}_n$;
furthermore, the dependence on~$\tilde X$ is continuous. The continuous
dependence
is ensured using a mollification argument [see below~\eqref
{bcp-main-18}] that has previously been used in \cite{KuMa}.
\end{remark}

The following lemma is a straightforward consequence of the
Lipschitz property of the Skorohod map,
the linearity of the cost and the state dynamics.
Proof is given in the \hyperref[appen]{Appendix}.\vadjust{\goodbreak}
\begin{lemma}
\label{ab210}
There is a $c_1\in(0, \infty)$ such that, if $q \in\R_+^{\II}, T
\in(0, \infty)$ and~$\tilde Y^1, \tilde Y^2 \in\tilde\mathcal{A}(q)$ defined on a common filtered
probability space
are such that
\[
\tilde Y^i(T+\cdot)- \tilde Y^i(T)= \bar\Gamma\bigl(\tilde Q^i(T) +
\tilde\zeta
(T+\cdot)-\tilde\zeta(T)\bigr),\qquad  i = 1,2,
\]
where $\tilde Q^i$ is defined by the right-hand side of \eqref{ab440} by
replacing~$\tilde Y$ there by~$\tilde Y^i$,
then
\[
|\tilde J(q, Y^1) - \tilde J(q, Y^2)| \le c_1\E
|Y^1 - Y^2|_{\infty, T}.
\]
\end{lemma}

Define $\vartheta\dvtx \R_+^{\II} \times\R^{\JJ} \to\R^{\JJ}$ as
\[
\vartheta(q_0, y) = y + \bar\Gamma(q_0+ Ry \curvi)(1),\qquad q_0 \in\R
_+^{\II},   y \in\R^{\JJ}.
\]
Note that $\vartheta$ is a Lipschitz map: that is, for some
$\vartheta_{\mathrm{lip}} \in(0, \infty)$, we have for $(q_0,y),
(\tilde
q_0, \tilde y) \in\R_+^{\II} \times\R^{\JJ}$,
%
\begin{equation}\label{110}
|\vartheta(q_0, y) - \vartheta(\tilde q_0, \tilde y)| \le\vartheta
_{\mathrm{lip}}
(|q_0 - \tilde q_0| + |y-\tilde y|).
\end{equation}
Also, since $\vartheta(q_0, 0) = 0$, we have for $(q_0,y) \in\R
_+^{\II
} \times\R^{\JJ}$,
%
\begin{equation}\label{111}
|\vartheta(q_0, y)| \le\vartheta_{\mathrm{lip}}|y|.
\end{equation}
We now present the near-optimal control that will be used in the proof
of Theorem~\ref{jrtoj}.
Recall the set $\mathbb{T}$ introduced in Corollary \ref{tset}.
Fix $\varepsilon> 0$, and a unit vector $y^* \in\mathbb{T}$ and define
$c_2\in[1, \infty)$ as
\[
c_2=\max\bigl\{2c_1(p_0+1)\bigl(1+L|\operatorname{diag}(x^*)C'|\bigr)|Ry^*|, 1
\bigr\}.
\]
Let $\tilde{Y}_0^{(1)}$, $\tilde{Y}^{(1)}$ be as in Theorem \ref
{main-bcp} with $\varepsilon= \varepsilon/2$.
Let $\varepsilon_0 = \varepsilon/c_2$ and define $\vartheta_{\varepsilon
_0}(x, y)=
\vartheta(x, y) + \varepsilon_0y^* $. Define control process $\tilde Y
\in\tilde
\mathcal{A}(q)$ with the corresponding state process
$\tilde Q$ [defined by the right-hand side of \eqref{ab440}] by the following
equations. For $n=0, 1, \ldots, p_0$,
%
\begin{equation}
\label{ab610}
\tilde Y(n\theta)-\tilde Y(n\theta-) = \vartheta_{\varepsilon
_0}(\tilde Q(n\theta-),
\partial\tilde Y_0^{(1)}(n)),
\end{equation}
and
%
\begin{equation}
\label{ab616}
\qquad\quad\tilde Y(t+n\theta)-\tilde Y(n\theta) = \bar\Gamma\bigl(\tilde Q(n\theta
) + \tilde\zeta
(\cdot+ n\theta) - \tilde\zeta(n\theta)\bigr)(t),
\qquad t \in[0, \theta),
\end{equation}
with the conventions that for $n=p_0$, $[0, \theta)$ is replaced by
$[0, \infty)$ and for $n=0$,
$\tilde Q(n\theta-) = q$. The control $\tilde Y$ evolves in a similar manner
to $\tilde Y^{(1)}$ at all time points excepting $n\theta$, $n=0, 1,
\ldots
, p_0$. Since $y^* \in\mathbb{T}$, for every $q_0 \in\R_+^{\II}$ and
$y \in\R^{\JJ}$,
$q_0 + R \vartheta_{\varepsilon_0}(q_0,y) > 0$. This, along with the
definition of the map~$\bar\Gamma$ (see below
Assumption \ref{ab1012}), ensures that $\tilde Q$ is nonnegative over
the time intervals $(n\theta, (n+1)\theta)$ and
$\tilde Q(n\theta) > 0$, for $n = 0, 1, \ldots, p_0$. Furthermore,~\eqref{abn756}
and the property
$Ky^* > 0$ (see definition of $\mathbb{T}$ in Corollary \ref{tset})
ensure that~$K \tilde Y$ is nondecreasing\vadjust{\goodbreak}
and nonnegative. Thus, the process defined by relations \eqref{ab610}
and \eqref{ab616} is indeed an element
of $\tilde\mathcal{A}(q)$. 
The strict positivity of $\tilde Q$ at time instants $n \theta$, $n
\le
p_0$ will be exploited in the weak convergence analysis of Section \ref
{secproof}
[see equation \eqref{ab1432} and also below~\eqref{442half}].

\begin{theorem}
\label{newmain518}
The process $\tilde Y$ defined above is $\varepsilon$-optimal for the BCP with
initial value $q$.
\end{theorem}
\begin{pf}
Since $\tilde Y^{(1)}$ is $\varepsilon/2$ optimal, in view of Lemma \ref{ab210},
it suffices to show that
%
\begin{equation}
\bigl|\tilde Y^{(1)} - \tilde Y\bigr|_{\infty, T} \le\frac
{\varepsilon}{2c_1}. \label
{ab909}
\end{equation}
For this we will introduce a collection of $\JJ$-dimensional processes
$\tilde Y_{(n)}$, $n = 0, 1, \ldots, p_0+1$,
such that $\tilde Y_{(0)} = \tilde Y^{(1)}$ and $\tilde Y_{(p_0+1)} =
\tilde Y^{(1)}$. These processes are only used in the current proof and
do not appear elsewhere in this work. 

Define, recursively, for $n=0,1,\ldots, (p_0+1)$, processes $\tilde
Y_{(n)}$ with corresponding state processes $\tilde Q_{(n)}$, as follows:
\[
\bigl(\tilde Q_{(0)},\tilde Y_{(0)}\bigr) = \bigl(\tilde Q^{(1)},\tilde Y^{(1)}\bigr),
\]
and for $n \ge0$, $t\ge0$,
\begin{eqnarray*}
\tilde Q_{(n+1)}(t)&= &\tilde Q_{(n)}(t) 1_{[0,n\theta)}(t) + \Gamma
\bigl(\tilde
H^{(n)}\bigr)(t-n\theta)1_{[n\theta, \infty)}(t),
\\
\tilde Y_{(n+1)}(t)&=& \tilde Y_{(n)}(t) 1_{[0,n\theta)}(t) + \bigl[\tilde
Y_{(n)}(n\theta-) +
\vartheta_{\varepsilon_0}\bigl(\tilde Q_{(n)}(n\theta-),\partial\tilde
Y_0(n)\bigr) \\
&&\hspace*{166pt}{}+ \bar
\Gamma\bigl(\tilde H^{(n)}\bigr)(t-n\theta)\bigr]1_{[n\theta, \infty)}(t),
\end{eqnarray*}
where for all $ t \ge0$,
\begin{eqnarray*}
\tilde H^{(n)}(t) &= &\tilde Q_{(n)}(n\theta-) + R\vartheta_{\varepsilon
_0}\bigl(\tilde
Q_{(n)}(n\theta-),\partial\tilde Y_0(n)\bigr)
+ \tilde\zeta(t+n\theta)-\tilde\zeta(n\theta)
\\
&&{}+ R[\tilde Y(t+n\theta) - \tilde Y(n\theta)].
\end{eqnarray*}
Note that for $t \in[n\theta,\infty)$,
\[
\tilde Y_{(n)}(t)=\tilde Y_{(n)}(n\theta-)+ \vartheta\bigl(\tilde
Q_{(n)}(n\theta
-),\partial\tilde Y_0(n)\bigr) +
\bar\Gamma\bigl(H^{(n)}\bigr)(t-n\theta),
\]
where for all $t \ge0$,
\begin{eqnarray*}
H^{(n)}(t) &=& \tilde Q_{(n)}(n\theta-) + R\vartheta\bigl(\tilde
Q_{(n)}(n\theta
-),\partial\tilde Y_0(n)\bigr)+
\tilde\zeta(t+n\theta)-\tilde\zeta(n\theta)
\\
&&{}+ R[\tilde Y(t+n\theta) - \tilde Y(n\theta)].
\end{eqnarray*}
Using the Lipschitz property of $\bar\Gamma$, it follows that, for $t
\ge0$,
\[
\bigl|\tilde Y_{(n+1)}(t)-\tilde Y_{(n)}(t)\bigr| \le\varepsilon_0\bigl(1+
L|\operatorname{diag}
(x^*)C'|\bigr)|Ry^*|.
\]
Thus, for $t\ge0$,
\[
\bigl|\tilde Y_{(p_0+1)}(t)-\tilde Y_{(0)}(t)\bigr| \le(p_0+1)\varepsilon_0\bigl(1+
L|\operatorname{diag}
(x^*)C'|\bigr)|Ry^*| \le\frac{\varepsilon}{2c_1}.
\]
The result follows on noting that $\tilde Y_{(0)} = \tilde Y^{(1)}$ and
$\tilde
Y_{(p_0+1)} = \tilde Y$.
\end{pf}

\subsection{\texorpdfstring{Proof of Theorem \protect\ref{main-bcp}}{Proof of Theorem 3.5}} \label{pfmainbcp}

Throughout this section we fix $q \in\R_+^{\II}$ and $\varepsilon
\in(0,
\infty)$.
We begin with some preparatory results. Let $\tilde\mathcal{A}_0$ be
the class
of all $\JJ$-dimensional adapted processes $Y$ given on some filtered
probability space such that $U=KY$
is nondecreasing, $U(0) \ge0$ and $q+RY(0) \ge0$. Note that $\tilde
\mathcal{A}
(q) \subset\tilde\mathcal{A}_0$. Also, a given $\tilde Y \in\tilde
\mathcal{A}
_0$ is in
$\tilde\mathcal{A}(q)$ if and only if
\eqref{ab440} is satisfied.

Given an adapted process $Y_0$, on some system $\Phi$, with sample
paths in~$\mathcal{D}^{\JJ}$, and satisfying
$q+RY(0) \ge0$,
we will
denote the process $Y$, defined by
%
\begin{equation}\label{ab315} Y = Y_0 +
\bar\Gamma(\tilde\zeta+
RY_0),
\end{equation}
as $\Upsilon(Y_0)$. We claim that
%
\begin{equation}\label{abnov7}
\mbox{if } Y_0\in\tilde\mathcal{A}_0,
\mbox{ then } \tilde Y =
\Upsilon(\tilde Y_0) \in\tilde\mathcal{A}(q).
\end{equation}
Indeed, $\tilde Q = \tilde\zeta+ R\tilde Y = \Gamma(\tilde\zeta+
R\tilde Y_0) \ge0$.
Also, $K\tilde Y = K \tilde Y_0 + K \bar\Gamma(\tilde\zeta+ R
\tilde Y_0)$. Since
$Y_0 \in\tilde\mathcal{A}_0$, $KY_0$ is nondecreasing and $KY_0(0)
\ge0$.
Also, for $x \in\mathcal{D}_+^{\II}$,
$K \bar\Gamma(x) = K\operatorname{diag}(x^*) C' \hat\Gamma(x)$,
which is a
nonnegative and nondecreasing function
since $\hat\Gamma(x)$ has these properties and the matrix
\[
K\operatorname{diag}(x^*)C' = \left[
\matrix{
Bx^*_n &0\cr
0 & 0}
\right] C'
\]
has nonnegative entries.
Combining these observations, we see that the process $\tilde Y =
\Upsilon(\tilde
Y_0)$ satisfies
\eqref{ab440} and \eqref{ab934}. The claim follows. %

Next, from the Lipschitz property of $\Gamma$ it follows that there
is a $\bar L \in(1, \infty)$ such that, if $\tilde Y_0^{(i)} \in
\tilde\mathcal{A}
_0$, $i=1,2$, then for all $T > 0$,
%
\begin{equation}
\label{ab161}
\bigl|\Upsilon\bigl(\tilde Y_0^{(1)}\bigr) - \Upsilon\bigl(\tilde Y_0^{(2)}\bigr)\bigr|_{\infty,T}
\le\bar L
\bigl|\tilde Y_0^{(1)}-\tilde Y_0^{(2)}\bigr|_{\infty,T}.
\end{equation}
In what follows, we will denote
$\sigma\{\tilde X_s\dvtx  0 \le s \le t\}$ by $\mathcal{F}_t^{\tilde X}$.

\begin{theorem}\label{stepwise-bcp}
Let $Y \in\tilde\mathcal{A}(q)$ be a $\{\mathcal{F}_t^{\tilde X}\}
$-adapted process
with a.s. continuous paths.
Suppose further that for some $m>0$,
%
%
\begin{equation}
\label{stepwise-bcp-1} \E[|Y|_{\infty,t}^m] <\infty\qquad  \mbox{for all } t>0.
\end{equation}
Then for any $\varepsilon_1, T \in(0, \infty)$, there are $\eta, M \in(0,
\infty)$, $p_0 \in\N$ and a $\tilde Y^{(1)} \in\tilde\mathcal{A}(q)$
such that $\tilde Y^{(1)} = \Upsilon(\tilde Y_0^{(1)})$ for some
$\tilde Y_0^{(1)}
\in\tilde\mathcal{A}_0$, that is, $\{\mathcal{F}_t^{\tilde X}\}
$-adapted and
satisfies~\textup{(i)}
and \textup{(ii)} of Theorem \ref{main-bcp} with $\theta= T/p_0$ and
%
\begin{equation}\label{ab102}
\E\bigl[\bigl|Y-\tilde Y_0^{(1)}\bigr|_{\infty,
T}^m\bigr]<\varepsilon
_1.
\end{equation}
\end{theorem}

\begin{pf}
The construction of $\tilde Y^{(1)}$ proceeds by defining,
successively, simpler approximations of $Y$, denoted as $Y^{(1)},
Y^{(2)}, Y^{(3)}, Y^{(4)}$. The process $Y^{(1)}$ is given in terms of
a sequence $\{Y_n\}$ of $\JJ$-dimensional processes,
whereas the processes $Y^{(3)}$ and $Y^{(4)}$ are given in terms of one
parameter families of $\JJ$-dimensional processes $\{Y(\theta, \cdot),
\theta> 0 \}$, $\{Y^*(M, \cdot), M > 0 \}$, respectively. All the
processes $Y^{(i)}$, $i=1,2,3,4$, and
$\{Y_n\}$, $\{Y(\theta, \cdot)\}$, $\{Y^*(M, \cdot)\}$ are only used in
this proof and do not appear elsewhere in the paper. %

Fix $\varepsilon_1, T \in(0, \infty)$.
Define $Y_n(t) = n\int_{(t-{1}/{n})^+}^t Y(s)\,ds, n\geq1, t \geq
0$. Note that for all
$t, t' \in[0,T]$,
\[
|Y_n(t)-Y_n(t')| \le2n |t-t'| |Y|_{\infty, T}.
\]
Hence, by \bbref{stepwise-bcp-1}, we have, for each $n$,
%
%
\begin{equation}
\label{step-pf-1} \E\biggl[ \sup_{t,t'\in[0,T]}
\biggl|\frac{Y_n(t)-Y_n(t')}{t-t'} \biggr|^m\biggr] < \infty.
\end{equation}
Note that for $t \in[1/n, T]$, $|Y_n(t) - Y(t)| \le\varpi
^T_{Y}(1/n)$ and
for $t \in[0, 1/n]$, $|Y_n(t) - Y(t)| \le2|Y|_{\infty, 1/n}$.
Since $Y$ is continuous and $Y(0) = 0$, $\varpi^T_{Y}(1/n) +
2|Y|_{\infty, 1/n} \to0$ a.s.
Combining this with \eqref{stepwise-bcp-1} and the estimate $ |Y_n -
Y|_{\infty, T} \le2 |Y|_{\infty, T}$,
we now have that, for some $n_0 \in\mathbb{N}$, $Y^{(1)} \Df Y_{n_0}$
satisfies
%
\begin{equation}
\label{step-pf-2} \E\bigl|Y^{(1)} - Y\bigr|_{\infty, T}^m \le\tilde\varepsilon_0,
\end{equation}
with $\tilde\varepsilon_0 = \varepsilon_1/4$. Also,
%
%
\begin{equation}
\label{step-pf-3}
\E\biggl[\sup_{t,t'\in[0,T]}\biggl|\frac{Y^1(t)-Y^1(t')}{t-t'}
\bigg|^m\biggr] \Df C_1 <\infty.
\end{equation}
Note that
%
%
\begin{equation}
\label{step-pf-4} Y^{(1)}(0)=0 \quad\mbox{and}\quad
Y^{(1)}\in\tilde\mathcal{A}_0.
\end{equation}
Given $p_0\in{\N}$ and $\theta= T/p_0$, define
\[
Y(\theta, t) =
Y^1 \biggl(\biggl\lfloor\frac{t}{\theta}\biggr\rfloor\theta
\biggr)1_{[0,T)}(t) + Y^1(T)1_{[T,\infty)}(t).
\]
Fix $p_0$ large enough so that $\theta< (\frac{\tilde\varepsilon
_0}{C_1})^{1/m}$ and set $Y^{(2)}(t) = Y(\theta, t)$.
Then, from~\eqref{step-pf-3} we have
%
%
\begin{eqnarray}\label{ab1206}
&&\E\bigl[\bigl|Y^{(2)}(t)-Y^{(1)}(t)\bigr|^m_{\infty,T}\bigr]
\nonumber\\
&&\qquad= \E\Bigl[\max_{n=0,\ldots,p_0-1}\Bigl\{\sup_{t\in[n\theta,
(n+1)\theta)}
\bigl|Y^{(1)}(n\theta)-Y^{(1)}(t)\bigr|^m\Bigr\}\Bigr]
\nonumber
\\[-8pt]
\\[-8pt]
\nonumber
&&\qquad\leq \theta^m
\E\biggl[\sup_{n=0,1,\ldots,p_0-1}\biggl\{\sup_{t,t'\in[n\theta,
(n+1)\theta)}
\frac{|Y^{(1)}(t)-Y^{(1)}(t')|}{|t-t'|}\biggr\}^m\biggr]\\
&&\qquad\leq \theta^m C_1 <  \tilde\varepsilon_0.\nonumber
\end{eqnarray}
From \eqref{step-pf-4} we have
%
\begin{equation}
Y^{(2)}(0)=0 \quad\mbox{and}\quad Y^{(2)} \in\tilde\mathcal{A}_0.
\label
{step-pf-6}
\end{equation}
For $x \in\R$, let $\lceil x\rceil$ denote the smallest integer upper
bound for $x$.
For \mbox{$x \in\R^{\JJ}$}, let $\lceil x\rceil= (\lceil x_1\rceil,
\ldots,
\lceil x_{\JJ}\rceil)'$.
Fix $\eta\le\frac{{\tilde\varepsilon_0}^{1/m}}{p_0\sqrt{\JJ}}$ and, with
convention $Y^{(2)}(-\theta) = 0$,
define for $t \ge0$
\[
Y^{(3)}(t) = \sum_{n=0}^{\lfloor{t\wedge T}/{\theta}\rfloor}
\bigl\lceil
\partial Y^{(2)}(n)/\eta\bigr\rceil\eta,
\]
where for $y \in\mathcal{D}^{\JJ}$, $\partial y(n)$ denotes
$y(n\theta) -
y((n-1)\theta)$.
Note that for $t \le T$,
\[
Y^{(2)}(t) =\sum_{n=0}^{\lfloor{t\wedge T}/{\theta}\rfloor}
\partial Y^{(2)}(n).
\]
Observing that for $x \in\R^{\JJ}$, $|x - \lfloor x/\eta\rfloor
\eta
| \le\eta\sqrt{\JJ}$ and recalling that $T=p_0\theta$, we have that
%
\begin{equation}
\E\bigl[\bigl|Y^{(3)}-Y^{(2)}\bigr|^m_{\infty,T}\bigr] \leq
\bigl(p_0\eta\sqrt{\JJ}\bigr)^m \le\tilde\varepsilon_0. \label{ab1208}
\end{equation}
Note that if $y \in\R^{\JJ}$ satisfies $Ky \geq0$, then $y_j \le0$
for all $j = \BB+1 , \ldots,\JJ$
and, consequently, for such $j$, $\lceil y_j\rceil\le0$. Combining
this with the fact that $A$ has nonnegative entries, we see that
$K\lceil y \rceil\geq0$.
From this observation, along with \eqref{step-pf-6}, we have
%
\begin{equation}
Y^{(3)}(0)=0 \quad\mbox{and}\quad Y^{(3)}\in\tilde\mathcal{A}_0.
\label
{step-pf-6ins}
\end{equation}
The process $Y^{(3)}$ constructed above is constant on $[n\theta,
(n+1)\theta)$ and the jumps $\partial Y^{(3)}(n)$ take value in the
lattice $\{k\eta\dvtx k\in{\IZ}\}$, for $n = 0, \ldots, p_0$.
Also, $Y^{(3)}(t) = Y^{(3)}(\theta p_0) = Y^{(3)}(T)$ for $t \ge T$.

For fixed $M \in(0, \infty)$, define
\[
Y^*(M,t) = \sum_{n=0}^{\lfloor{t}/{\theta}\rfloor}
\partial Y^{(3)}(n) I_{\{|\partial Y^{(3)}(n)|\leq M\}},\qquad t \ge0.
\]
Then there exists $C_2 \in(0, \infty)$ such that, for all $M > 0$,
%
\begin{equation}\label{step-pf-75}
\E\bigl[\bigl|Y^*(M, \cdot)-Y^{(3)}\bigr|^m_{\infty, T}\bigr] \leq
C_2\sum_{n=0}^{p_0}\E\bigl[ \bigl|\partial Y^{(3)}(n)\bigr|^m I_{(|\partial
Y^{(3)}(n)|>M)}\bigr].\hspace*{-30pt}
\end{equation}
Also, for some $C_3 \in(0, \infty)$, we have from \eqref{stepwise-bcp-1},
\eqref{step-pf-2}, \eqref{ab1206} and \eqref{ab1208} that, for $n = 0,
1,\ldots, p_0$,
\[
\E\bigl[\bigl|\partial Y^{(3)}(n)\bigr|^m\bigr] \le C_3(\E[|Y|^m_{\infty, T}] + 1) <
\infty.
\]
Fix $M > 0$ such that the right-hand side of \eqref{step-pf-75} is bounded
by $\tilde\varepsilon_0$.
Setting $Y^{(4)} = Y^*(M, \cdot)$, we now have that
%
\begin{equation}
\E\bigl[\bigl|Y^{(4)}-Y^{(3)}\bigr|^m_{\infty,T}\bigr] \leq
\tilde\varepsilon_0. \label{ab1221}
\end{equation}
Also,
%
\begin{equation}\label{step-pf-79}
Y^{(4)}(0)=0 \quad\mbox{and}\quad Y^{(4)}
\in\tilde\mathcal{A}_0.
\end{equation}
Combining \eqref{stepwise-bcp-1},
\eqref{step-pf-2}, \eqref{ab1206}, \eqref{ab1208} and \eqref{ab1221},
we now have that $\tilde Y_0^{(1)} = Y^{(4)}$ satisfies
\eqref{ab102} as well as (i) and (ii) of Theorem \ref{main-bcp}.
This completes the proof.
\end{pf}


\begin{lemma}\label{cont-bcp}
For each $\varepsilon_1 > 0$ there exists an $\varepsilon_1$-optimal $Y \in
\tilde\mathcal{A}
(q)$, which is $\{\mathcal{F}_t^{\tilde X}\}$-adapted, continuous
a.s., and satisfies
%
\begin{equation}
\label{ab127}
\limsup_{T \to\infty} e^{-\gamma T} \E|Y|_{\infty, T}^m =0  \qquad\mbox{for every } m > 0.
\end{equation}
\end{lemma}
\begin{pf}Fix $\varepsilon_1 > 0$ and let $w = \Lambda q$. Applying Theorem
2.1(iv) of \cite{AtBu},
we have that $\tilde J_0^*(w) = \inf\tilde J_0(w, U)$, where the
infimum is
taken over all $\{\mathcal{F}^{\tilde X}_t\}$-adapted controls $U$. Hence,
using \eqref{mincost-equal}, we conclude that
there is an $\mathbf{N}$-dimensional $\{\mathcal{F}_t^{\tilde X}\}
$-adapted process
$U$ for which \eqref{EWF-1} holds and
%
\begin{equation}
\label{bd622}
\tilde J^*(q) = \tilde J_0^*(w) \ge\tilde J_0(w, U) - \varepsilon_1.
\end{equation}
From
Lemma 4.7 of \cite{AtBu} and following the construction of Proposition
3.3 of~\cite{BuRo} [cf. (12) and (14) of that paper], we can assume without
loss of generality that $U$ has continuous sample paths
and for all $m > 0$,
\[
\limsup_{T \to\infty} e^{-\gamma T}\E|GU|_{\infty, T}^m = 0.
\]
Hence, using properties of the $G$ matrix (see Assumption \ref
{G-column} ), we have that
%
\begin{eqnarray}
\label{cond612}
\limsup_{T \to\infty} e^{-\gamma T}\E|U(T)|^m \le c^{-m}
\limsup_{T \to\infty} e^{-\gamma T}\E|GU|_{\infty, T}^m = 0
\nonumber
\\[-8pt]
\\[-8pt]
\eqntext{\mbox{for all } m > 0.}
\end{eqnarray}
We will now use a construction given in the proof of Theorem 1 of \cite
{harri1}. This construction shows that
there is a $\JJ\times\mathbf{N}$ matrix $F_1$ and a $\JJ\times\II$
matrix $F_2$
such that letting
%
\begin{equation}
\label{cond615}
Y(t) = F_1 U(t) + F_2 \bigl(q^*(\tilde W(t)) - q - \tilde X(t)\bigr),\qquad t > 0,
\end{equation}
we have that $Y \in\tilde\mathcal{A}(q)$ and $\tilde J_0(w, U) =
\tilde
J(q, Y).$
We refer the reader to equations~(35) and (36) of \cite{harri1} for
definitions and constructions of these matrices.
From~\eqref{bd622} we now have that $Y$ is an $\varepsilon_1$-optimal control,
has continuous sample paths a.s. and is $\{\mathcal{F}_t^{X}\}$-adapted.
Finally from \eqref{cond615}, we have that for some $C_2 \in(0,
\infty)$,
\[
\E|Y|_{\infty, T}^m \le C_2 \bigl(1 + \E|U(T)|^m + T^m + \E|q^*(\tilde
W)|_{\infty, T}^{m}\bigr).
\]
By combining \eqref{32half} and \eqref{EWF-1}, the fourth term on the
right-hand side can be bounded above
by $C_3(1 + T^m + \E|U(T)|^m)$ for some $C_3>0$.
The result then follows on using
\eqref{cond612}.
\end{pf}

The following construction will be used in the proof of Theorem \ref{main-bcp}.
%
\begin{lemma}
\label{inter159}
Fix $Y \in\tilde\mathcal{A}(q)$ such that $\tilde J(q, Y) < \infty$ and
\eqref
{ab127} holds. For $T > 0$, let $\tilde Y^T_0(t) = Y(t \wedge T)$, $t >
0$, and
$Y^T = \Upsilon(\tilde Y^T_0)$. Then given $\varepsilon_1 > 0$, there
exists $T \in
(0, \infty)$ such that
$|\tilde J(q, Y) - \tilde J(q, Y^T)| < \varepsilon_1$ and \eqref{ab127}
holds with
$Y$ replaced by $Y^T$.
\end{lemma}
\begin{pf}
Since $\tilde J(q, Y) < \infty$, we have that
\[
L(T, Y) \Df\E\int_T^{\infty} e^{-\gamma t} h \cdot\tilde Q(t) \,dt +
\E
\int_{(T,\infty)} e^{-\gamma t} p \cdot dU(t)
\to0\qquad \mbox{as } T \to\infty,
\]
where $\tilde Q$ is the state process corresponding to $Y$ and $U = KY$.
Choose $T_1$ large enough so that
%
\begin{equation}\label{lab16}
L(T, Y) < \varepsilon_1/2
\qquad\mbox{for } T \ge T_1.
\end{equation}
Using the Lipschitz property of $\Gamma$ and $\hat\Gamma$ (see
Assumption \ref{ab1012}), we can find
$C_1 \in(0, \infty)$ such that, for all $T > 0$,
%
\begin{eqnarray}
\label{ab302}
&&|\tilde Q^T(t) - \tilde Q(T)| + |Y^T(t) - Y(T)|
\nonumber
\\[-8pt]
\\[-8pt]
\nonumber
&&\qquad\le C_1 \sup_{T \le s
\le t}
|\tilde\zeta(s) - \tilde\zeta(T)|,\qquad
t \ge T,
\end{eqnarray}
where $\tilde Q^T$ is the state process corresponding to $Y^T$. Thus, for
some $C_2 \in(0, \infty)$,
\[
\E\int_T^{\infty} e^{-\gamma t} h \cdot\tilde Q^T(t) \,dt \le\frac
{\bar
h}{\gamma} e^{-\gamma T} \E|\tilde Q(T)|
+ C_2 \int_T^{\infty} e^{-\gamma t} (1+t) \,dt,
\]
where $\bar h = \max_{i \in\II} h_i$. Using \eqref{ab127}, we can now
choose $T_2$ large enough so that
%
\begin{equation}\label{lab17} \E\int_T^{\infty} e^{-\gamma t} h
\cdot\tilde Q^T(t)
\,dt \le\varepsilon_1/4\qquad \mbox{for all } T > T_2.
\end{equation}
Next, letting $U^T = KY^T$, we have from \eqref{ab302} that for some
$C_3 \in(0, \infty)$,
\[
\E|U^T(t) - U(T)| \le C_3 \bigl(1 + (t-T)\bigr)\qquad \mbox{for } 0 < T < t.
\]
Integration by parts now yields that for some $T_3 > 0$,
%
\begin{equation}\label{lab18}\E\int_{(T, \infty)} e^{-\gamma T} p
\cdot dU^T(t)
\le\varepsilon_1/4\qquad \mbox{for } T \ge T_3.
\end{equation}
Combining the estimates in \eqref{lab16}, \eqref{lab17} and \eqref
{lab18}, we now have that
for all $T \ge\max\{T_1, T_2, T_3\}$,
\[
|\tilde J(q, Y) - \tilde J(q, Y^T)| \le L(T, Y) + L(T, Y^T) \le
\varepsilon_1/2 +
\varepsilon_1/4 + \varepsilon_1/4 = \varepsilon_1.\vadjust{\goodbreak}
\]
Finally, the fact that \eqref{ab127} holds with $Y$ replaced by $Y^T$
is an immediate consequence of \eqref{ab302}.
\end{pf}

We can now complete the proof of Theorem \ref{main-bcp}.

\begin{pf*}{Proof of Theorem \protect\ref{main-bcp}}
Using Lemma
\ref{cont-bcp}, one can find $Y \in\mathcal{A}(q)$ which is
$\varepsilon
/5$-optimal,
has continuous paths a.s., is $\{\mathcal{F}_t^{\tilde X}\}$-adapted and
satisfies \eqref{ab127}. Using Lemma \ref{inter159},
we can find $T \in(0, \infty)$ such that $Y^T = \Upsilon( Y(\cdot
\wedge
T))$ is $2\varepsilon/5$-optimal and
\eqref{ab127} holds with $Y$ replaced by $Y^T$.
We will apply Theorem \ref{stepwise-bcp}, with $m=1$, $Y$ replaced with
$Y^T$, $\varepsilon_1$ replaced by $\varepsilon/(5c_1\bar L)$ and denote the
corresponding processes
obtained from Theorem~\ref{stepwise-bcp}, once again by $\tilde Y^{(1)}$
and $ \tilde Y_0^{(1)}$. 
In particular,
$\tilde Y^{(1)} \in\mathcal{A}(q)$ is such that $\tilde Y^{(1)} =
\Upsilon
(\tilde
Y^{(1)}_0)$, where $\tilde Y_0$ is
$\{\mathcal{F}_t^{\tilde X}\}$-adapted, satisfies (i) and (ii) of
Theorem~\ref{main-bcp}, for some $\eta, M, \theta\in(0, \infty)$ and $p_0 \in
\N
$, and~\eqref{ab102} holds
with $m=1$, $Y$ replaced by $Y^T$ and
where $\bar L$ is as in \eqref{ab161}. Then
\[
\E\bigl[\bigl|Y^T - \tilde Y^{(1)} \bigr|_{\infty, T}\bigr] \le\bar L \E\bigl[\bigl|Y^T-\tilde
Y_0^{(1)}\bigr|_{\infty,T}\bigr] \le\frac{\varepsilon}{5c_1}.
\]
Thus, from Lemma \ref{ab210}, $\tilde Y^{(1)}$
is $3\varepsilon/5$-optimal.

Processes $\tilde Y_0^{(1)}, \tilde Y^{(1)}$ as in the statement of
Theorem \ref{main-bcp} will be constructed by modifying
the processes $\tilde Y_0^{(1)}, \tilde Y^{(1)}$ above (but denoted
once more by the same symbols), by constructing successive
approximations $(Y_0^{(\kappa)}, Y^{(\kappa)})$ and $(\tilde
Y_0^{(\gamma)}, Y^{(\gamma)})$. These approximations are only used in
the current proof and do not appear elsewhere in the paper. %

Consider $\kappa> 0$ such that $\theta/\kappa\in\N$. Let $\mathcal{G}
^{(n)}_\kappa= \sigma\{ \mathcal{X}^\kappa(n)\},
\mathcal{G}^{(n)} = \sigma\{\tilde X(s)\dvtx\allowbreak  s \le n\theta\}, \mathcal
{G}= \sigma\{
\tilde X(s)\dvtx  s \ge0\}, n \in\mathbb{N},
\kappa> 0$.
Since $\mathcal{G}^{(n)}_\kappa\uparrow
\mathcal{G}^{(n)}$ as $\kappa\downarrow0$ and $\tilde Y^{(1)}_0$ is
$\{
\mathcal{F}
_t^{\tilde X}\}$-adapted, we have for each fixed $\varsigma\in
\mathcal{S}_M^{\eta}$,
%
\begin{eqnarray}\label{ab401}
\IP\bigl[\partial
\tilde Y^{(1)}_0(n )=\varsigma |  \mathcal{G}_\kappa^{(n)}\bigr]
&\rightarrow&
\IP\bigl[\partial
\tilde Y^{(1)}_0(n )=\varsigma |   \mathcal{G}^{(n)}\bigr]
\nonumber
\\[-8pt]
\\[-8pt]
\nonumber
& = & 1_{\{
\partial
\tilde
Y^{(1)}_0(n) = \varsigma\}}\qquad \mbox{a.e., as } \kappa\downarrow0.
\end{eqnarray}
Note that for fixed $\varsigma\in\mathcal{S}_M^{\eta}$ and $n = 1,
\ldots, p_0$,
\[
\IP\bigl[\partial
\tilde Y_0^{(1)}(n )=\varsigma |  \mathcal{G}_\kappa^{(n)}\bigr] =
p^{\kappa
}_{n,\varsigma}(\mathcal{X}^{\kappa}(n)),
\]
for some measurable map
$p^{\kappa}_{n,\varsigma}\dvtx  \R^{nj_0\II} \to[0, 1]$ satisfying
%
\begin{equation}
\label{ab316}
\mbox{for all } x \in\R^{\II nj_0}\qquad
\sum_{\varsigma\in\mathcal{S}_M^{\eta}} p^{\kappa}_{n,\varsigma
}(x) = 1.
\end{equation}
Using Lemma \ref{appnov} in the \hyperref[appen]{Appendix}, we can construct, by suitably
augmenting the
filtered probability space, an adapted process $Y_0^{(\kappa)}$ such that
$Y_0^{(\kappa)}(0) = 0$, $Y_0^{(\kappa)}(t) = Y_0^{(\kappa)}(\lfloor
t/\theta\rfloor\theta)$
and $\partial Y_0^{(\kappa)} (n) \Df Y_0^{(\kappa)}(n \theta) -
Y_0^{(\kappa)}((n-1)\theta)$ satisfies
%
\begin{equation}\label{ab459}\IP\bigl[\partial
 Y_0^{(\kappa)}(n )=\varsigma |  \mathcal{G}\vee\mathcal
{Y}_0^{n-1}\bigr] =
p^{\kappa
}_{n,\varsigma}(\mathcal{X}^{\kappa}(n)),\qquad
1 \le n \le p_0, \varsigma\in\mathcal{S}_M^{\eta},\hspace*{-35pt}
\end{equation}
where $\mathcal{Y}_0^n = \sigma\{\partial Y_0^{(\kappa)} (j), j \le
n \}$.
Note that if $f$ is a real bounded continuous map on $\mathcal{C}^{\II}$,
then by successive conditioning
and \eqref{ab401}, as $\kappa\downarrow0$,
\begin{eqnarray*}
\E\Biggl( f(\tilde X) \prod_{n=1}^{p_0} 1_{\{\partial
Y_0^{(\kappa)}(n ) = \varsigma_n\}}\Biggr)& =&
\E\Biggl( f(\tilde X) \prod_{n=1}^{p_0} p^{\kappa}_{n,\varsigma
_n}(\mathcal{X}
^{\kappa}(n))\Biggr)\\[-2pt]
&\to&
\E\Biggl( f(\tilde X) \prod_{n=1}^{p_0} 1_{\{\partial
\tilde Y_0^{(1)}(n ) = \varsigma_n\}}\Biggr),
\end{eqnarray*}
for any fixed $(\varsigma_1, \ldots, \varsigma_{p_0})$.
Thus, in particular,
%
\begin{equation}
\label{ab419}
\bigl(\tilde X, Y_0^{(\kappa)}\bigr) \Rightarrow\bigl(\tilde X, \tilde Y_0^{(1)}\bigr)
\qquad\mbox{as }
\kappa\downarrow0.
\end{equation}
Define $Y^{(\kappa)} = \Upsilon(Y_0^{(\kappa)})$.
Then $Y^{(\kappa)} \in\mathcal{A}(q)$ and satisfies (i) and (ii) of Theorem~\ref{main-bcp}. Also,
\eqref{ab419} along with the Lipschitz property of the Skorohod map yields
that $\tilde J(q,Y^{(\kappa)}) \to\tilde J(q,\tilde Y^{(1)})$ as
$\kappa
\downarrow0$. Fix $\kappa$ sufficiently small
so that~$Y^{(\kappa)}$ is $4\varepsilon/5$-optimal and, suppressing
$\kappa$,
denote $Y_0^{(\kappa)}$ by $\tilde Y_0$,
$Y^{(\kappa)}$ by $\tilde Y$
and the generating kernels by
$p_{n,\varsigma}$.

Finally, in order to ensure property (iii) in the theorem, we mollify
the kernels $p_{n,\varsigma}$ as follows. For
$\gamma> 0$, and $n=1,\ldots, p_0$, define
%
\begin{equation}\label{bcp-main-18}\hat p^\gamma_{n,\varsigma}(x)
\Df
\int_{\R^{nj_0\II}} p_{n,\varsigma}(x+z) \prod_{j=1}^{nj_0}(\phi
^\gamma(z_j)
\,dz_j),\qquad
x\in\R^{nj_0\II}, \varsigma\in\mathcal{S}_M^{\eta},\hspace*{-35pt}
\end{equation}
where $\phi^\gamma$ is the density function of an $\II$-dimensional
Normal random variable with
mean 0 and variance $\gamma I$. Note that the map
$x \mapsto\hat p^\gamma_{n,\varsigma}(x)$ is continuous for every
$\gamma,
n, \varsigma$ and
\eqref{ab316} is satisfied with $p^{\kappa}_{n,\varsigma}(x)$ replaced
with $\hat p^\gamma_{n,\varsigma}(x)$.
From continuity of the maps\vspace*{1pt} $\hat p^\gamma_{n,\varsigma}$, we can
find measurable
maps (suppressing dependence on $\gamma$ in notation)
$\varpi_n\dvtx  \R^{\II nj_0} \times[0, 1] \to\mathcal{S}_M^{\eta}$,
$n = 1, \ldots
,p_0$, such that
$\varpi_n(\cdot, t)$ is continuous, at every $x \in\R^{\II nj_0}$, for
a.e. $t$ in $[0, 1]$, and
if $\mathcal{U}$ is a Uniform random variable on $[0,1]$, then
\[
\IP\bigl(\varpi_n(x,\mathcal{U})=\varsigma\bigr)=\hat p^\gamma_{n,\varsigma}(x),
\qquad x\in\R
^{\II nj_0}, \varsigma\in
\mathcal{S}_M^{\eta}.
\]
Let $\{\tilde{\mathcal{U}}_n\}$ be
an i.i.d sequence of Uniform random variables, that is, independent of~$\tilde X$.
Now construct $\tilde Y_0^{(\gamma)}$ such that $\tilde Y_0^{(\gamma)}(0)=0$,
\[
\partial\tilde{Y}_0^{(\gamma)}(n) = \tilde{Y}_0^{(\gamma)}(n\theta
)-\tilde
{Y}_0^{(\gamma)}\bigl((n-1)\theta\bigr)=
\varpi_n(\mathcal{X}^\kappa(n), \tilde{\mathcal{U}}_n),\qquad n \ge1
\]
and $\tilde{Y}_0^{(\gamma)}(t)=\tilde{Y}_0^{(\gamma)}(\theta
\lfloor t/\theta
\rfloor)$. Note that
\[
\IP\bigl[\partial
\tilde Y_0^{(\gamma)}(n )=\varsigma |  \mathcal{G}\vee\mathcal
{Y}_0^{n-1}\bigr] =
\hat
p^{\gamma}_{n,\varsigma}(\mathcal{X}^{\kappa}(n)),
\qquad n \ge1, \varsigma\in\mathcal{S}_M^{\eta}.\vadjust{\goodbreak}
\]
Since $\hat p_{n,\varsigma}^{\gamma} \to p_{n,\varsigma}$ pointwise, as
$\gamma\to0$, we have
for every real bounded map~$f$ on $\mathcal{C}^{\II}$,
\begin{eqnarray*}
\E\Biggl( f(\tilde X) \prod_{n=1}^{p_0} 1_{\{\partial
\tilde Y_0^{(\gamma)}(n ) = \varsigma_n\}}\Biggr) &=&
\E\Biggl( f(\tilde X) \prod_{n=1}^{p_0} \hat p^{\gamma}_{n,\varsigma
_n}(\mathcal{X}^{\kappa}(n))\Biggr)\\
&\to&
\E\Biggl( f(\tilde X) \prod_{n=1}^{p_0} p_{n,\varsigma_n}(\mathcal
{X}^{\kappa
}(n))\Biggr)
\\
&=& \E\Biggl( f(\tilde X) \prod_{n=1}^{p_0} 1_{\{\partial
\tilde Y_0(n ) = \varsigma_n\}}\Biggr),
\end{eqnarray*}
for every $(\varsigma_1, \ldots,\varsigma_{p_0})$, as $\gamma\to0$.
Thus,
$
(\tilde X, \tilde Y_0^{(\gamma)}) \Rightarrow(\tilde X, \tilde Y_0)$,
as $\gamma\to0$.
Let $\tilde Y^{(\gamma)} = \Upsilon(\tilde Y_0^{(\gamma)})$.
Then $\tilde Y^{(\gamma)} \in\mathcal{A}(q)$ and
the above weak convergence and, once more, the Lipschitz property of
the Skorohod map yields
that $\tilde J(q,\tilde Y^{(\gamma)}) \to\tilde J(q,\tilde Y)$ as
$\gamma\downarrow
0$. Recall that $\tilde Y$ is $4\varepsilon/5$-optimal.
We now choose $\gamma$ sufficiently small so that $\tilde Y^{(\gamma
)}$ is
$\varepsilon$-optimal. By construction,
$\tilde Y_0^{(1)} \Df\tilde Y_0^{(\gamma)}$ and $\tilde Y^{(1)} \Df
\tilde Y^{(\gamma
)}$ satisfy all the properties stated in the theorem.
\end{pf*}

\section{Asymptotically near-optimal controls for SPN}\label{secproof}

The goal of this section is to prove Theorem \ref{jrtoj}. Fix $q \in
\R
_+^{\II}$ and $\varepsilon\in(0,1)$. Let $\tilde Y \in\tilde\mathcal{A}(q)$
be the $\varepsilon
$-optimal control introduced
above Theorem \ref{newmain518}. Fix $q^r \in\N_0^{\II}$, $r > 0$ such
that $\hat q^r \to q$, as $r \to\infty$. Section \ref{subconstruct} below
gives the construction of
the sequence of policies $\{T^r\}$, $T^r \in\mathcal{A}^r(q^r)$, such that
$J^r(q^r, T^r) \to\tilde J(q, \tilde Y)$, yielding the proof of
Theorem \ref{jrtoj}.
The latter convergence of costs is proved in Section \ref{secconvprf}.
The main ingredient in this proof is
Theorem \ref{mainweak} whose proof is given in Section \ref{weakcgce}.
For the rest of this section $\tilde Y$ as in Theorem \ref{newmain518}
and parameters $T, \eta, M, p_0, j_0, \theta, \kappa$ that specify
$\tilde
Y$ shall be fixed. In addition,
let $\rho\in(0, \infty)$ and $r_0 \ge1$ be such that
%
\begin{equation}\label{rho-choice}
\rho\Bigl[\min_{x_j^* \neq0} x_j^*\Bigr] >
M(\vartheta_{\mathrm{lip}}+1)
\quad \mbox{and} \quad  r_0 \theta> \rho,
\end{equation}
where $\vartheta_{\mathrm{lip}}$ is as in \eqref{110}.
\subsection{Construction of the policy sequence} \label{subconstruct}
$\!\!\!$We will only specify a \mbox{$T^r\!\in\!\mathcal{A}^r(q^r)$} for $r\ge r_0$ and so
henceforth, without loss of generality, we assume
$r > r_0$.

In this section, since $r \ge r_0$ will be fixed,
the superscript $r$ will frequently be suppressed from the notation.
The following additional notation will be used.
For $n \le p_0$, define $a(n) = nr^2\theta$, $b(n) = nr^2\theta+
r\rho$,
$\mathcal{I}(n) = [a(n), a(n+1))$, $\mathcal{I}_1(n) = [a(n), b(n))$
and $\mathcal{I}_2(n) = [b(n), a(n+1))$, where we set
$a(p_0+1) = \infty$. Note $\mathcal{I}(n) = \mathcal{I}_1(n) \cup
\mathcal{I}_2(n)$, $n
\le p_0$.

Recall $\curvm$ introduced in Assumption \ref{ldp}. Fix $\curvk\in(0,
1)$ such that
%
\begin{equation}\label{lab23}
 \curvk(1 +\curvm) - 2 \Df\upsilon> 1.\vadjust{\goodbreak}
\end{equation}
Fix $d_1 \in(0, \infty)$ such that $d_1 > \sup_r\{\beta_j^r\dvtx  j \in
\mathbb{J}\} + 1$.
Define
%
\begin{equation}
\label{ab411}
\Delta^r \Df r^{\curvk},\qquad \Theta^r(s) \Df d_1r^{\curvk} 1_{(\cup
_{n\le p_0}\mathcal{I}_2(n))}(s)\qquad \mbox{for } s \ge0.
\end{equation}
%
Also
define $\curvp^r\dvtx  [0, \infty) \to[0, \infty)$ as
\[
\curvp^r(s) \Df\cases{
s, &\quad $\mbox{if } s \in\mathcal{I}_1(n), n = 0, 1, \ldots,
p_0,$\vspace*{3pt}\cr
b(n) + \displaystyle\biggl\lfloor\frac{s-b(n)}{\Delta} \biggr\rfloor\Delta, & \quad$\mbox{if } s
\in\mathcal{I}_2(n), n = 0, 1, \ldots, p_0.$}
\]
%
Thus, if $s \in\mathcal{I}_2(n)$ for some $n$, $\Theta^r(s) =
d_1\Delta^r$
and $\curvp^r(s)$ equals the left end point of the $\Delta$-subinterval
in which
$s$ falls.
Otherwise,
if $s \in\mathcal{I}_1(n)$ for some~$n$, $\Theta^r(s)=0$ and $\curvp^r(s)
= s$.


Recall the probability space $(\Omega, \mathcal{F}, \IP)$ introduced
in Section
\ref{secsetup} which supports all the random
variables and stochastic processes introduced therein. Let $\{\mathcal
{U}_i\dvtx  i
\in\mathbb{N}\}$ be a sequence of Uniform random
variables on $[0, 1]$ on this probability space (constructed by
augmenting the space if needed), independent of the $\sigma$-field
$\bar\mathcal{F}^r$ defined in \eqref{fr-sigmafield-defn}.
This sequence will be used in the construction of
the control policy.

%
\begin{figure}

\includegraphics{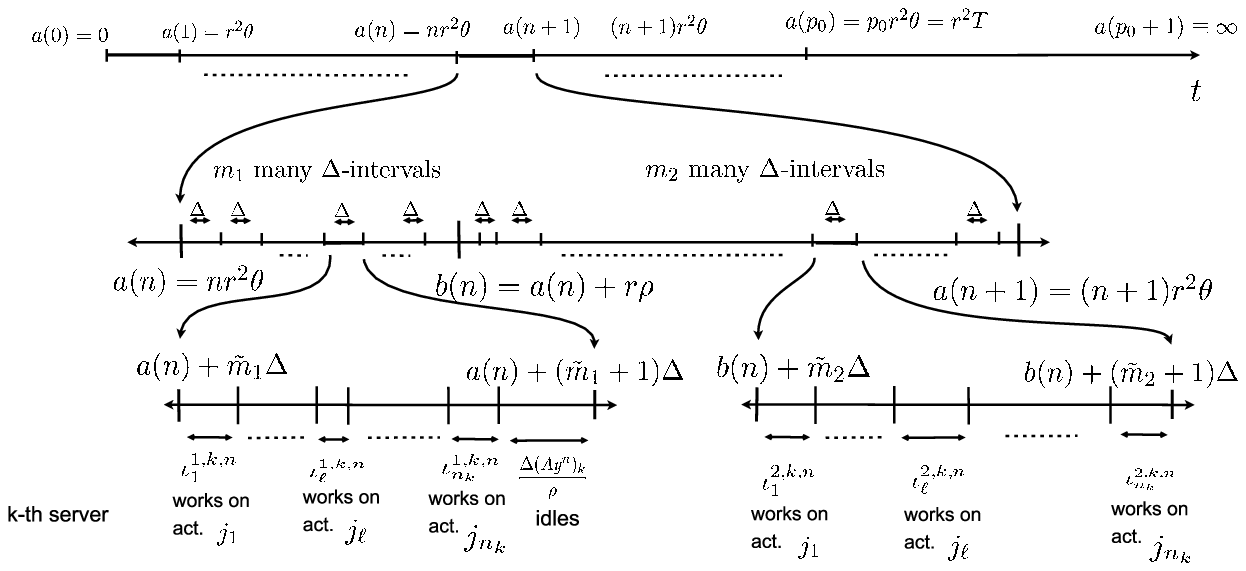}

\caption{The figure above shows the behavior, under $T^{(1)}$, of the
$k$th server, that is, responsible for the $n_k$ activities ($j_1, \ldots
, j_{n_k}$). 
The actual policy $T$ is given as a certain modification of $T^{(1)}$
that ensures feasibility and a strict positivity property.}
\label{policyfig}
\vspace*{6pt}
\end{figure}

The policy $T \equiv T^r$ is constructed recursively over time
intervals $\mathcal{I}(n), n = 0, \ldots, p_0$, as follows. We will describe
the effect of the policy on the $k$th server (for each $k \in\mathbb{K}$)
at every time-instant $s \ge0$ (see Figure \ref{policyfig}). Recall the set $\mathbb{J}(k)$ introduced
in Section \ref{secsetup}. Let for $k \in\mathbb{K}$,
$j_1 < j_2< \cdots<j_{n_k}$ be the ordered elements of $\mathbb{J}(k)$ (the
set of activities that the $k$th server can perform).

\textit{Step 1}: [$T^r(s)$ when $s \in\mathcal{I}(0)$].
Let $m_0 = \lfloor\frac{r^2\theta}{\Delta} \rfloor$.
Recall that $\tilde Y(0) = \varepsilon_0 y^*$.
Write $\nu^0 = \varepsilon_0 y^*$. For $k \in\KK$, let
%
\begin{equation}\label{abn532a}
\nabla_\ell^{1,k,0} \Df\biggl(x_{j_\ell
}^* - \frac
{\nu^0_{j_\ell}}{\rho}\biggr)
\Delta,\quad \nabla_\ell^{2, k,0} \Df x_{j_\ell}^* \Delta,\quad
\ell
= 1, \ldots, n_k.
\end{equation}
Note that, by our choice of $\rho$, $\nabla_\ell^{1,k,0} > 0$ if $j_l$
is a basic activity. Also, if
$j_l$ is nonbasic, then $x^*_{j_l} = 0$, but since $K\nu^0 \ge0$,
we have $\nu^0_{j_l} \le0$.
Thus, $\nabla_\ell^{1,k,0} \ge0$ for every
$l,k$. Also, from Assumption \ref{assum-HT} [see \eqref{xstar-defn}]
and recalling that $A\nu^0 \ge0$, we get
%
\begin{equation}\label{ab1743}
\sum_{l=1}^{n_k}\nabla_\ell^{1,k,0}
= \Delta\biggl(1
- \frac{1}{\rho} (A\nu^0)_k\biggr) \le\Delta,\quad
\sum_{l=1}^{n_k}\nabla_\ell^{2,k,0} = \Delta,\quad k
\in\mathbb{K}.\hspace*{-30pt}
\end{equation}
Let $m_1 \equiv m_1^r \Df\lfloor\frac{r\rho}{\Delta}\rfloor$ and
$m_2 \equiv m_2^r \Df\lfloor\frac{r^2\theta- r\rho}{\Delta
}\rfloor$.
Define, for each $k \in\mathbb{K}$, and $l= 1, \ldots, n_k,$
\[
\dot{T}^{(1)}_{j_l}(s) = \cases{
\displaystyle\sum_{m=0}^{m_1 -1} 1_{\mathcal{E}_l^{m, 0}}(s), & \quad$s \in\mathcal
{I}_1(0),$\vspace*{2pt}\cr
\displaystyle\sum_{m=0}^{m_2 -1} 1_{\hat\mathcal{E}_l^{m, 0}}(s), & \quad$s \in
\mathcal{I}_2(0),$}
\]
where for $\tilde m_i = 0, \ldots,m_i-1$, $i=1,2$, $l = 1, \ldots, n_k$,
%
\begin{eqnarray}
\mathcal{E}_l^{\tilde m_1,0} &=& \Biggl[a(0) + \tilde m_1\Delta+
\sum_{i=0}^{l-1} \nabla_i^{1,k,0}, a(0) + \tilde m_1\Delta+
\sum_{i=0}^{l} \nabla_i^{1,k,0}\Biggr),\hspace*{-55pt}\label{ab8a}\\
\hat\mathcal{E}_l^{\tilde m_2, 0} &=& \Biggl[b(0)
+ \tilde m_2\Delta+ \sum_{i=0}^{l-1} \nabla_i^{2,k,0}, b(0) +
\tilde m_2\Delta+
\sum_{i=0}^{l} \nabla_i^{2,k,0}\Biggr)\hspace*{-55pt} \label{ab8b}
\end{eqnarray}
and, by convention, $\nabla_0^{i,k,0} = 0$, $i=1,2$.
Since $\{\mathbb{J}(k), k \in\mathbb{K}\}$ gives a partition of
$\mathbb{J}$, the
above defines
the $\JJ$-dimensional process $\dot{T}^{(1)}(s)$ for $s \in[a(0),
a(1))$. We set
\[
T^{(1)}(s) = \int_0^s \dot{T}^{(1)}(s) \,ds,\qquad s \in\mathcal{I}(0).
\]
Next, define
$\mathcal{V}(s) \equiv\mathcal{V}^r(s) = (Q^r(s), X^r(s), T^r(s))$
for $s \in
\mathcal{I}
(0)$ by the system of equations below: 
%
\begin{eqnarray}\label{ab543}
 Q^r(s) &= &q^r + X^r(s) + r \bigl(\theta_1^rs - (C - P')\operatorname
{diag}(r^2\theta_2^r)
\bar{T}^r(s/r^2)\bigr)\nonumber\\
&&{} + R\bigl(x^*s - T^r(s)\bigr),
\nonumber
\\[-8pt]
\\[-8pt]
\nonumber
X^r(s) &=& r \hat{X}^r(s/r^2),\qquad \hat X^r \mbox{ defined by \eqref
{xr-defn}, } \\
T^r_j(s) &=& \int_{[0, s]} 1_{\{Q^r_{\sigma_1(j)}(u-) > 0\}} 1_{\{
Q^r_{\sigma_1(j)}(\curvp^r(u))
> \Theta^r(u)\}} \,dT^{r,(1)}_j(u),\qquad j \in\mathbb{J},\hspace*{-35pt}
\nonumber
\end{eqnarray}
where (as introduced in Section \ref{secsetup}), $\sigma_1(j)$ denotes
the index of the buffer that the $j$th activity is associated with.
The above construction can be interpreted as follows.
The policy $T$ ``attempts'' to implement $T^{(1)}$ over~$\mathcal{I}_1(0)$,
unless the corresponding buffer is empty (in which case, it idles).
Over $\mathcal{I}_2(0)$, $T$ has a similar behavior, but, in addition, it
idles when the corresponding buffer does not have at least $\Theta$---many
jobs at the beginning of the $\Delta$-subintervals.
Processing jobs only when there is a ``safety-stock'' at each buffer at
the beginning of the $\Delta$-subinterval ensures that with ``high probability''
either all the activities associated with the given buffer receive
nominal time effort over the interval or all of them receive zero
effort. This, in particular, makes sure that the idling processes
associated with the policy are consistent with the reflection terms for
the Skorohod map associated with
the constraint matrix~$D$ [see \eqref{ab1249}]. This property will be
exploited in the weak convergence arguments of Section \ref{secconvprf}
[see, e.g., arguments below~\eqref{ab1042}].

Let $I^r$ be defined by \eqref{i-defn}.
By construction, $T^r(0) =0$, $I^r(0) = 0$ and
%
\begin{eqnarray}
\label{ab549}
Q^r(s) \ge0,\qquad T^r(s) \ge0,\qquad I^r(s) \ge0
\nonumber
\\[-8pt]
\\[-8pt]
\eqntext{\mbox{for } s\in\mathcal{I}
(0) \mbox{ and }
T^r, I^r \mbox{ are nondecreasing on } \mathcal{I}(0).}
\end{eqnarray}
This completes the construction of the policy and the associated
processes on $\mathcal{I}(0)$.



\textit{Step 2}: [$T^r(s)$ for $s \in\mathcal{I}(n)$, for $1\le n \le
p_0$].
Suppose now that the process~$\mathcal{V}(s)$ has been defined
for $s \in[a(0), a(n))$, where $1\le n \le p_0$. We now describe the
construction over the interval
$\mathcal{I}(n) =[a(n), a(n + 1))$. Let $\hat\chi^{\kappa, r}(n) =
\{ \hat X^r(l\kappa)\dvtx  l = 1, \ldots, n j_0\}$.
Recall that $nj_0\kappa= n\theta$. Define
%
\begin{equation}\label{ab201}\bar\nu^{r, n} \Df\varpi_{n}(\hat
\chi^{\kappa,
r}(n), \mathcal{U}_{n}),\qquad
\nu^{r, n} \Df\vartheta_{\varepsilon_0}(\hat Q^r(n\theta-), \bar\nu
^{r, n}),
\end{equation}
where, as in Section \ref{secsetup}, $\hat Q^r(t) = Q^r(r^2t)/r$, $t
\ge0$.
We will suppress $\kappa$ and~$r$ from the notation and write
$(\hat\chi^{\kappa, r}(n),\bar\nu^{r, n}, \nu^{r, n})
\equiv(\hat\chi(n),\bar\nu^{n}, \nu^{n})$. For each $k \in
\mathbb{K}$, define
%
\begin{equation}\label{abn532b} \nabla_\ell^{1,k,n} \Df\biggl(x_{j_\ell
}^* - \frac
{\nu^n_{j_\ell}}{\rho}\biggr)
\Delta,\quad \nabla_\ell^{2, k,n} \Df x_{j_\ell}^* \Delta,\quad
\ell
= 1, \ldots, n_k.
\end{equation}
As before, \eqref{ab1743} (with $\nabla_\ell^{j,k,0}$ replaced by
$\nabla_\ell^{j,k,n}$) is satisfied.
Define, for each $k \in\mathbb{K}$, and $l= 1, \ldots, n_k,$
%
\begin{equation}
\dot{T}^1_{j_l}(s) = \cases{
\displaystyle\sum_{m=0}^{m_1 -1} 1_{\mathcal{E}_l^{m, n}}(s), & \quad$s \in\mathcal
{I}_1(n),$\vspace*{2pt}\cr
\displaystyle\sum_{m=0}^{m_2 -1} 1_{\hat\mathcal{E}_l^{m, n}}(s), &\quad $ s \in
\mathcal{I}_2(n),$}\label{abnov8}
\end{equation}
where for $\tilde m_i = 0, \ldots,m_i-1$, $i=1,2$, $l = 1, \ldots, n_k$,
$\mathcal{E}_l^{\tilde m_1,n}, \hat\mathcal{E}_l^{\tilde m_2, n}$
are defined by
the right-hand side of \eqref{ab8a} and \eqref{ab8b} respectively,\vadjust{\goodbreak} with
$(a(0), b(0), \nabla_i^{1,k,0},\allowbreak \nabla_i^{2,k,0})$ replaced by $(a(n),
b(n), \nabla_i^{1,k,n}, \nabla_i^{2,k,n})$
We set
\[
T^1(s) = T(a(n)) + \int_{a(n)}^s \dot{T}^1(s) \,ds,\qquad s \in\mathcal{I}(n)
\]
and define $\mathcal{V}(s)$ for $s \in\mathcal{I}(n)$ by the system of
equations in
\eqref{ab543}.

The above recursive procedure gives a construction for the process
$\mathcal{V}
^r(s) = (Q^r(s), X^r(s), T^r(s))$
for all $s \in[0, \infty)$.

The policy constructed above clearly satisfies parts (i), (ii) and
(iii) of Definition~\ref{t-adm-defn}. In fact, it also satisfies part
(iv) of the
definition and, consequently, we have the following result. The proof
is given in the \hyperref[appen]{Appendix}.

\begin{prop}
\label{tisadmisnew}
For
all $r \ge r_0$, $T^r \in\mathcal{A}^r(q^r)$.
\end{prop}

\subsection{Convergence of costs} \label{secconvprf}
In this section we prove Theorem \ref{jrtoj} by showing that, with $\{
T^r\}$ as in Section \ref{subconstruct} and $\tilde Y$ as introduced above
Theorem \ref{newmain518}, $J^r(q^r, T^r) \to\tilde J(q, \tilde Y)$,
as $r
\to\infty$. We begin with an elementary lemma.

\begin{lemma}
\label{invpr}
Let $g^n \in\mathcal{D}^{\II}$, $n\ge1$, $g \in\mathcal{C}^{\II}$
be such that $g^n \to g$ u.o.c. as \mbox{$n \to\infty$}. Let $\varepsilon_n \in(0,
\infty)$, $n\ge1$,
be such that $\varepsilon_n \to0$ as $n \to\infty$. Let $B^n, B$ be
$\II
\times\II$ matrices such that $B^n \to B$ as $n\to\infty$.
Suppose $f^n, h^n \in\mathcal{D}^{\II}$ and $\gamma^n \in\mathcal
{D}^1$ satisfy for
all $t \ge0$:
\begin{longlist}[(iii)]
\item[(i)] $f^n(t) \ge0$, $\gamma^n(t) \ge0$,

\item[(ii)]
$f^n(t) = g^n(t) + 
B^n h^n(t), h^n_i(t) = \int_{[0,t]} 1_{\{f^n_i(\gamma^n(s)) \le
\varepsilon_n\}}\,ds,
i \in\mathbb{I}$
and

\item[(iii)] $|\gamma^n(t) - t| \le\varepsilon_n$.
\end{longlist}
Then $(f^n, g^n, h^n, \gamma^n)$ is precompact in $\mathcal{D}
^{3\II
+1}$ and any limit point $(f,g,h,\gamma)$ satisfies for all $t\ge0$,
$\gamma(t)=t$; $f(t)= g(t)+Bh(t)$; $h(t)= \int_{[0,t]}\tilde h(s) \,ds$,
where $\tilde h\dvtx  [0, \infty) \to[0,1]^{\II}$ is a measurable map
such that
\[
\int_{[0,t]} 1_{\{f_i(s) > 0\}} \tilde h_i(s) \,ds = 0 \qquad\mbox{for all
} t
\ge0.
\]
\end{lemma}

The proof of the above lemma is given in the \hyperref[appen]{Appendix}. Recall the
definitions of various scaled processes given in \eqref
{fl-scaled}--\eqref{w-relation}, and that $\curvi(t) = t$ for $t \ge0$.
In addition, we define $\bar T^{r,1}(t) = T^{r,1}(r^2t)/r^2$,
$\hat T^{r,1}(t) = T^{r,1}(r^2t)/r$, $r> 0$, $t \ge0$.

\begin{prop}
\label{tisadmis}
As $r \to\infty$, $\bar T^r \to x^*\curvi$,
u.o.c. in probability.
\end{prop}
\begin{pf}
From the definition of $T^{r,1}$ [see \eqref{abnov8}], and the
observation that
the interval $\hat\mathcal{E}_l^{m, n}$ has length $x^*_{j_l}\Delta$ [see
\eqref{ab1743} and \eqref{ab8b}], we have
that over each interval $[m\Delta, (m+1)\Delta)$, that is, contained in
$\mathcal{I}_2(n)$ for some $n \le p_0$, the Lebesgue measure of the time
instants $s$ such that $\dot{T}_j^1(s) =1$ equals $x^*_j \Delta$, $j
\in\mathbb{J}$.
This is equivalent to the statement that
%
\begin{eqnarray}
\label{ab248}
\int_{[m\Delta, (m+1)\Delta)} d\bigl(T^{r,1}_j(s)-x^*_js\bigr)=0
\nonumber
\\[-8pt]
\\[-8pt]
\eqntext{ \mbox{whenever
} \bigl[m\Delta, (m+1)\Delta\bigr) \subset
\mathcal{I}_2(n) \mbox{ for some } n \le p_0 .}
\end{eqnarray}
Also, noting that for $j \in\mathbb{J}$, $0 \le\dot{T}^{r,1}_j \le1$
and $x^*_j \le1$, we have that for some $C_1>0$, 
%
\begin{eqnarray}
\label{ab248a}
\int_{[m\Delta, (m+1)\Delta)} d\bigl(T^{r,(1)}_j(s)+x^*_js\bigr)<C_1\Delta
\nonumber
\\[-8pt]
\\[-8pt]
\eqntext{\mbox
{whenever } \bigl[m\Delta, (m+1)\Delta\bigr) \subset
\mathcal{I}(n) \mbox{ for some } n \le p_0 .}
\end{eqnarray}
Fix $j \in\mathbb{J}$ and $t > 0$ such that $r^2t \in[a(n), a(n+1))$
for some $n \le p_0$. Then
\begin{eqnarray*}
&& r^2|x^*_jt - \bar T^{r,1}_j(t)|\\
&&\qquad \le\sum_{n_0=0}^{n}\sum_{\ell
=1}^2\sum
_{m\dvtx  [m\Delta, (m+1)\Delta) \in\mathcal{I}_\ell(n_0)}
\int_{[m\Delta\wedge t, (m+1)\Delta\wedge t)} d\bigl(T^{r,1}_j(s)-x^*_js\bigr).
\end{eqnarray*}
Using \eqref{ab248} and \eqref{ab248a} and the fact that the number of
$\Delta$ intervals in $\mathcal{I}_1(n)$ is bounded by $r\rho/\Delta$,
we see that
\[
r^2|x_j^* t - \bar T^{r,1}_j(t)| \le\biggl(\frac{r\rho(p_0+1)}{\Delta}
+1\biggr)C_1\Delta.
\]
Thus, for some $\varrho\in(0, \infty)$, 
%
\begin{equation}
\label{ab1254}
\sup_{0 \le t < \infty} |x^*t - \bar T^{r,1}(t)| \le\varrho/r
\qquad\mbox
{for all } r \ge r_0.
\end{equation}
Also, since $\bar T^r_j(t) \le t$ for all $t \ge0$, $j \in\mathbb
{J}$, we
get from
\eqref{xr-defn} and standard estimates for renewal processes (see, e.g.,
Lemma 3.5 of \cite{BudGho2}) that $\bar X^r \Df\hat X^r/r$ converges
to $0$ u.o.c. in probability, as $r \to\infty$. Combining this with
Assumption \ref{assum-HT} and \eqref{ab1254}, we have
%
\begin{equation}\label{ab218n} \bar\zeta^r \Df\hat\zeta^r/r \mbox{ converges to
} 0,
\mbox{ u.o.c., in probability. }
\end{equation}
Next, define $\bar\curvp^r(s) = \curvp^r(r^2s)/r^2$, $\bar\Theta
^r(s)=\Theta^r(r^2s)/r^2$ and
%
\begin{eqnarray}\label{lab33}
\mathcal{S}_j^r(t) &=& \bigl\{s \in[0, t]\dvtx
\bar
Q^r_{\sigma
_1(j)}(s)=0, \mbox{or }
\bar Q^r_{\sigma_1(j)}(\bar\curvp^r(s))\le
\bar\Theta^r(s) \bigr\}\nonumber\\
&= &\bigl\{s \in[0, t]\dvtx  \bar Q^r_{\sigma_1(j)}(\bar\curvp^r(s))\le
\bar\Theta^r(s) \bigr\}
\nonumber
\\[-8pt]
\\[-8pt]
\nonumber
&&{} \cup\bigl\{s \in[0, t]\dvtx  \bar Q^r_{\sigma_1(j)}(s)=0,
\bar Q^r_{\sigma_1(j)}(\bar\curvp^r(s)) >
\bar\Theta^r(s) \bigr\} \\
&\Df& \mathcal{S}_j^{r, 1}(t) \cup\mathcal{S}_j^{r, 2}(t).\nonumber
\end{eqnarray}
We will fix $t>0$ for the rest of the proof and suppress $t$ from the
notation when writing $\mathcal{S}_j^{r, i}(t)$, unless there is scope for
confusion. Using the above display and~\eqref{ab543}, we have that
%
\begin{equation}
\label{ab104}
\bar T^{r,1}_j(t) - \bar T^r_j(t) = \int_{[0, t]} 1_{\mathcal{S}_j^r}(s)
\,d\bar T^{r,1}_j(s),\qquad j \in\mathbb{J}, t \ge0.
\end{equation}
Using the fact that $\bar Q^r(\cdot) = \frac{\hat Q^r(\cdot)}{r}$ along
with \eqref{q-relation} 
and \eqref{ab104}, we can write
%
\begin{equation}\label{413half} \bar Q^r(t) = \bar\zeta^r(t) + R\bar
H^r(t) + R
\bar L^r(t),\qquad t \ge0,
\end{equation}
where for $j \in\mathbb{J}$,
%
\begin{eqnarray}\label{ab210n}\bar H^r_j(t) &=& \bigl(x^*_jt
- \bar T^{r,1}(t)\bigr) + \int_{[0, t]} 1_{\mathcal{S}_j^{r,1}}(s)\, d\bigl(\bar
T^{r,1}_j(s) - x^*_js\bigr)
\nonumber
\\[-8pt]
\\[-8pt]
\nonumber
&&{}+ \int_{[0, t]} 1_{\mathcal{S}_j^{r,2}}(s)
\,d\bar T^{r,1}_j(s)
\end{eqnarray}
and
\[
\bar L^r_j(t) = x^*_j\int_{[0, t]} 1_{\mathcal{S}_j^{r,1}}(s) \,ds =
(\operatorname{diag}(x^*) C'\bar\ell^r)_j(t),
\]
where
\[
\bar\ell^r_i(t) = \int_{[0, t]} 1_{\{\bar Q^r_{i}(\bar\curvp
^r(s))\le
\bar\Theta^r(s)\}} \,ds,\qquad t \ge0, i \in\mathbb{I}.
\]
Since $D= R \operatorname{diag}(x^*)C'$ [see \eqref{R-defn} and
\eqref{refmat}],
we have
%
\begin{equation}\label{ab1249} \bar Q^r = \bar\zeta^r + R\bar H^r +
D \bar\ell
^r.
\end{equation}
Next, letting
$\hat m_0^r = \lfloor\frac{r^2t}{\Delta^r}\rfloor+ 1$, we have from
the choice of $d_1$ [see above \eqref{ab411}] that
%
\begin{eqnarray}\label{lab71}
&&\IP(\mathcal{S}_j^{r, 2} \neq\varnothing)\nonumber\\
&&\qquad\le \sum_{k=0}^{\hat m_0^r} \sum_{i \in\mathbb{I}} \IP\bigl(Q_i^r(k\Delta
^r) >
d_1r^{\curvk}, Q_i^r(u) = 0\nonumber\\
&&\qquad\qquad\mbox{ for some } u \in\bigl[k\Delta^r, (k+1)\Delta^r\bigr)\bigr)\nonumber\\
&&\qquad\le \II\sum_{k=0}^{\hat m_0^r} \sum_{j \in\mathbb{J}}
\IP\bigl(S_j^r\bigl(T_j^r(k\Delta^r)+\Delta^r\bigr) - S_j^r(T_j^r(k\Delta^r)) \ge
d_1r^{\curvk} \bigr)\\
&&\qquad\le {\hat m_0^r}\II\JJ\frac{\varsigma(1)}{r^{\curvk\curvm}}
\nonumber\\
&&\qquad\le C_2 (t+1) r^{2-\curvk(1+\curvm)},\nonumber
\end{eqnarray}
for some $C_2 \in(0, \infty)$, where the next to last inequality makes
use of Assumption~\ref{ldp}. Recalling [from \eqref{lab23}] that
$\upsilon= \curvk(1+\curvm)-2 > 1$, we get that
%
\begin{equation}
\label{ab1230}
\IP(\mathcal{S}_j^{r, 2} \neq\varnothing) \le C_2(t+1) r^{-\upsilon}
\to0\qquad
\mbox{as }
r \to\infty.
\end{equation}
We
note that the above convergence only requires that $\upsilon>0$. The
property $\upsilon> 1$ will, however, be needed in the
proof of Proposition \ref{tight424} [see \eqref{ab632}]. Next, using
\eqref{ab248} and \eqref{ab248a}, we have 
for some $\varrho_1 \in(0, \infty)$,
%
\begin{equation}\label{ab1240}
\sup_{0 \le u < \infty}\biggl|\int_{[0, u]} 1_{\mathcal{S}_j^{r,1}}
\,d\bigl(\bar T^{r,1}_j(s) - x^*_js\bigr)\biggr| \le\frac{\varrho_1}{r} \to0
\qquad\mbox{as } r \to\infty,
\end{equation}
for all $j \in\mathbb{J}$. The above inequality follows from the fact that
the integral can be written as the sum of integrals over $\Delta
$-subintervals: When the subinterval is within some $\mathcal
{I}_2(n)$, the
integral is zero [using the definition of $\curvp^r(s)$ in such
intervals and \eqref{ab248}] and when the subinterval is within some
$\mathcal{I}_1(n)$
[the number of such intervals is $(p_0+1)m_1^r$ which can be bounded by
$C_3 \frac{r\rho}{\Delta}$ for some $C_3$], the integral
is bounded by $C_1\Delta/r^2$
from \eqref{ab248a}.

Now, combining \eqref{ab1254}, \eqref{ab1230} and \eqref{ab1240},
we have that, for each $j \in\mathbb{J}$, \mbox{$\bar H^r_j \to0$,}
u.o.c. in probability, as $r \to\infty$.
From \eqref{ab1249}, \eqref{ab218n}, Lemma \ref{invpr} and unique
solvability of the Skorohod problem for $(0, D)$, we now have that
$(\bar Q^r, \bar\ell^r) \to(0,0)$,
u.o.c. in probability, as $r \to\infty$. Thus, $\bar L^r$ converges to
$0$ as well.
The result now follows on noting that
$\bar T^r = x^*\curvi-\bar H^r - \bar L^r $.
\end{pf}

The following proposition gives a key estimate in the proof of Theorem~\ref{jrtoj}.

\begin{prop}\label{tight424}
For some $c_2 \in(0, \infty)$ and $\bar\upsilon\in(1, \infty)$,
\[
\E(|\hat Q^r|^{\bar\upsilon}_{\infty, t} + |\hat\zeta^r|^{\bar
\upsilon}_{\infty, t} +
|\hat Y^r|^{\bar\upsilon}_{\infty, t} ) \le c_2 (1 + t^3)
\qquad\mbox{for all } r \ge r_0, t \ge0.
\]
\end{prop}
\begin{pf}
Using standard moment estimates for renewal processes (cf. Lem\-ma~3.5 of
\cite{BudGho2}),
one can find $C_1 \in(0, \infty)$ such that
%
\begin{equation}
\label{ab425}
\E|\hat\zeta^r|^2_{\infty, t} \le C_1 (1+t^2) \qquad\mbox{for all
} r \ge r_0, t \ge0.
\end{equation}
From \eqref{413half} and \eqref{ab1249}, we have
\[
\hat Q^r(t) = \hat\zeta^r(t) + R\hat H^r(t) + R \hat L^r(t)
=\hat\zeta^r(t) + R\hat H^r(t) + D\hat\ell^r(t) ,\qquad t \ge0,
\]
where $\hat H^r = r \bar H^r$, $\hat L^r = r \bar L^r$ and
$\hat\ell^r = r \bar\ell^r$. We rewrite the above display as
\[
\hat Q^r(\bar\curvp^r(t)) = [\hat Q^r(\bar\curvp^r(t))-\hat
Q^r(t) ]+ \hat\zeta^r(t) + R\hat H^r(t) +
D\hat\ell^r(t) ,\qquad t \ge0.
\]
From Theorem 5.1
of \cite{will-invpr}, for some $C_2 \in(0, \infty)$,
%
\begin{eqnarray}\label{ab610n}
\qquad &&|\hat\ell^r|_{\infty, t} + |\hat Q^r|_{\infty, t}
\nonumber
\\[-8pt]
\\[-8pt]
\nonumber
&&\qquad\le C_2 \biggl(\hat q^r + |\hat\zeta^r|_{\infty, t} + |\hat
H^r|_{\infty, t} +
|\hat Q^r(\bar\curvp^r(\cdot))-\hat Q^r|_{\infty,t}+\frac
{d_1r^{\curvk}}{r}\biggr).
\end{eqnarray}
Also, from \eqref{ab210n}, \eqref{ab1254} and \eqref{ab1240}, for all
$t \ge0$ and $r \ge r_0$,
%
\begin{equation}
\biggl|\hat H^r_j(t) - r \int_{[0, t]} 1_{\mathcal{S}_j^{r,2}}(s)
\,d\bar T^{r,1}_j(s) \biggr| \le\varrho+ \varrho_1. \label
{ab612}
\end{equation}
Next, using \eqref{ab1230}, we get, for some $C_3 \in(0, \infty)$,
%
\begin{eqnarray}\label{ab632}
\E\biggl(r \sup_{0 \le u \le t}\int_{[0, u]}
1_{\mathcal{S}
_j^{r,2}}(s)
\,d\bar T^{r,1}_j(s) \biggr)^{\upsilon\wedge2} &\le& r^{\upsilon}C_3(t^3+1)
r^{-\upsilon}
\nonumber
\\[-8pt]
\\[-8pt]
\nonumber
 &\le&
C_3(t^3+1).
\end{eqnarray}
Finally, for some $C_4 \in(1, \infty)$, for all $a \ge1$,
%
\begin{eqnarray}\label{lab72}
&&\IP\bigl(|\hat Q^r(\bar\curvp^r(\cdot))-\hat Q^r(\cdot
)|_{\infty,t} \ge a\bigr)
\nonumber
\\[-8pt]
\\[-8pt]
\nonumber
&&\qquad \le \frac{r^2t}{\Delta^r}\biggl(\sum_{i \in\mathbb{I}}\IP
\biggl(A_i^r(\Delta
^r) \ge\frac{ar}{C_4}\biggr)
+\sum_{j \in\mathbb{J}}\IP\biggl(S_j^r(\Delta^r)
\ge\frac{ar}{C_4}\biggr)\biggr).\hspace*{-30pt}
\end{eqnarray}
Using moment estimates for renewal process once more (Lemma 3.5 of
\cite
{BudGho2}), we can find $C_5 \in(0, \infty)$
such that
\begin{eqnarray*}
\IP\bigl(|\hat Q^r(\bar\curvp^r(\cdot))-\hat Q^r(\cdot)|_{\infty,t}
\ge a\bigr)
&\le&
\frac{r^2t}{\Delta^r}\frac{C_5\Delta^r}{(ar-C_5\Delta^r)^2}\\
&\le& \frac{1}{a^2} \frac{C_5 r^2t}{(r-{C_5\Delta^r}/{a})^2}.
\end{eqnarray*}
Thus, there is an $r_1 \in(r_0, \infty)$ and $C_6 \in(0, \infty)$
such that, for all $r \ge r_1$,
%
\begin{equation}\label{ab248n}
\IP\bigl(|\hat Q^r(\bar\curvp^r(\cdot))-\hat Q^r(\cdot)|_{\infty,t}
\ge
a\bigr) \le\frac{1}{a^2}C_6(t+1).
\end{equation}
This shows that for some $\upsilon_1 \in(1, \infty)$ and $C_7\in(0,
\infty)$,
%
\begin{equation}
\label{ab629}
\E|\hat Q^r(\bar\curvp^r(\cdot))-\hat Q^r(\cdot)|^{\upsilon
_1}_{\infty
,t} \le C_7(1+t),\qquad t \ge0.
\end{equation}
The result now follows on using \eqref{ab425}, \eqref{ab612}, \eqref
{ab632} and \eqref{ab629} in \eqref{ab610n}
and observing that $\hat Y^r = \hat H^r + \operatorname
{diag}(x^*)C'\hat\ell^r$.
\end{pf}

In preparation for the proof of Theorem \ref{jrtoj}, we introduce the
following notation.
For $n = 0, \ldots, p_0-1$, we define processes $\curvq^{r,n},\curvz
^{r,n}$ with paths
in~$\mathcal{D}^{\II}_{\theta}$ and $\mathcal{C}^{\JJ}_{\theta}$,
respectively, as
\[
(\curvq^{r,n}(t), \curvz^{r,n}(t))\hspace*{-0.5pt} = \hspace*{-0.5pt}\cases{
\bigl(\hat Q^r(n\theta+ \rho/r), 0 \bigr), & \quad$t \in[0, \rho/r),$\vspace*{2pt}\cr
\bigl(\hat Q^r(t+ n\theta), \hat Y^r(t+n\theta) - \hat Y^r(n\theta
+ \rho/r)\bigr),& \quad$t \in[\rho/r,\theta].$}
\]
We denote by $\curvq^{r,p_0}, \curvz^{r,p_0}$ the processes with paths
in $\mathcal{D}^{\II}$
and $\mathcal{C}^{\JJ}$, respectively,
defined by the right-hand side in the display\vadjust{\goodbreak} above, by replacing $n$ by
$p_0$ and
$[\rho/r,\theta]$
with $[\rho/r,\infty)$.
Then
\[
\curvq^r \Df(\curvq^{r, *}, \curvq^{r, p_0}) \in\mathcal
{D}_{\theta
}^{p_0\II
}\times\mathcal{D}^{\II}
\qquad \mbox{a.s.},
\]
where $\curvq^{r,*} = (\curvq^{r, 0}, \ldots,\curvq^{r, p_0-1})$.
Similarly, define
the process $\curvz^r$ with paths in $\mathcal{C}_{\theta}^{p_0\JJ
}\times
\mathcal{C}
^{\JJ}$. Recall $\bar\nu^{r, n},
\nu^{r, n}$ introduced in
\eqref{ab201}.
Denote $\bar\nu^r = (\bar\nu^{r, 0}, \bar\nu^{r,1}, \ldots,\allowbreak\bar\nu^{r, p_0})$
and $\nu^r = ( \nu^{r, 0}, \nu^{r,1}, \ldots,\nu^{r,
p_0})$, where we set
$\bar\nu^{r, 0} = 0$ and $\nu^{r, 0} =\varepsilon_0 y^*$.

Next, for $n = 0, \ldots, p_0-1$, define processes $\curvq
^{(n)},\curvz
^{(n)}$ with paths
in $\mathcal{C}^{\II}_{\theta}$ and~$\mathcal{C}^{\JJ}_{\theta}$,
respectively, as
\[
\bigl(\curvq^{(n)}(t),\curvz^{(n)}(t)\bigr) = \cases{
\bigl(\tilde Q(t+n\theta), \tilde Y(t+n\theta) - \tilde Y(n\theta)\bigr),& \quad$t
\in[0,
\theta),$\vspace*{2pt}\cr
\bigl(\tilde Q\bigl((n+1)\theta-\bigr),\tilde Y\bigl((n+1)\theta-\bigr) - \tilde
Y(n\theta)\bigr),
& \quad$t
=\theta.$}
\]
Also, define $\curvq^{(p_0)}, \curvz^{(p_0)}$ by the first line of the
above display
by replacing~$\theta$ by $\infty$.
Then $\curvq\Df(\curvq^{*}, \curvq^{ (p_0)}) \in\mathcal
{C}_{\theta
}^{p_0\II
}\times\mathcal{C}^{\II}$,
a.s., where $\curvq^{*} = (\curvq^{(0)}, \ldots,\curvq^{(p_0-1)})$.
Similarly, define
the process $\curvz$ with paths in $\mathcal{C}_{\theta}^{p_0\JJ
}\times
\mathcal{C}
^{\JJ}$. Also, let for $n = 1, \ldots, p_0$, $\bar\nu^{(n)} =
\partial\tilde Y_0^{(1)}(n)=\tilde Y_0^{(1)}(n\theta) -
\tilde Y_0^{(1)}((n-1)\theta)$, where $\tilde Y_0^{(1)}$ is as
above~\eqref{ab610}, and $\nu^{(n)} = \tilde Y(n\theta)-\tilde Y(n\theta-)$. Then
%
\begin{equation}\label{ab735}\bar\nu^{(n)} \Df\varpi_{n}(\chi
^{\kappa}(n), \tilde
\mathcal{U}_{n}),\qquad
\nu^{(n)} \Df\vartheta_{\varepsilon_0}\bigl(\tilde Q(n\theta-), \bar\nu
^{(n)}\bigr).
\end{equation}

Define
$\bar\nu= (\bar\nu^{(0)}, \bar\nu^{(1)}, \ldots,\bar
\nu^{(p_0)})$
and
$\nu= (\nu^{(0)}, \nu^{(1)}, \ldots,\nu^{(p_0)})$, where
$\bar\nu^{(0)} = 0$ and $\nu^{(0)} = \varepsilon_0y^*$. Then $\bar
\nu^r, \bar\nu\in
(\mathcal{S}_M^{\eta})^{\otimes(p_0+1)}$ and $\nu^r,\nu\in\R
^{\JJ(p_0+1)}$.
Next, let
\[
\nu_0^{r,n} = \hat Y^r(n\theta+ \rho/r),\qquad \nu_0^{(n)} =
\tilde Y(n\theta),\qquad
n = 0, 1, \ldots, p_0 .
\]
Then $\nu_0^r \Df(\nu_0^{r,0}, \ldots,\nu_0^{r,p_0});
\nu_0 \Df(\nu_0^{(0)}, \ldots,\nu_0^{(p_0)}) \in\R^{\JJ
(p_0+1)}.$ Let
\[
\Xi= \mathcal{D}^{\II} \times(\mathcal{S}_M^{\eta})^{\otimes(p_0+1)}
\times\bigl(\R^{\JJ(p_0+1)}\bigr)
\times\bigl(\R^{\JJ(p_0+1)}\bigr) \times(\mathcal{D}_{\theta}^{p_0\II}
\times
\mathcal{D}
^{\II})
\times(\mathcal{C}_{\theta}^{p_0\JJ} \times\mathcal{C}^{\JJ}).
\]
Note that $\mathcal{J}^r \Df(\hat\zeta^r, \bar\nu^r, \nu^r, \nu
_0^r, \curvq^r, \curvz^r)$, $r \ge1$ and
$\mathcal{J}= (\tilde\zeta, \bar\nu, \nu, \nu_0, \curvq, \curvz
)$ are
$\Xi$-valued random variables.
The following is the main step in the proof of Theorem~\ref{jrtoj}.
\begin{theorem}
\label{mainweak}
As $r \to\infty$,
$ \mathcal{J}^r \Rightarrow\mathcal{J}.$
\end{theorem}

Proof of the above theorem is given in the next subsection. Using
Theorem~\ref{mainweak}, the proof of Theorem~\ref{jrtoj} is now
completed as follows.

\begin{pf*}{Proof of Theorem \protect\ref{jrtoj}}
From proposition and integration by parts,
%
\begin{eqnarray}
\label{ab804}
J^r(q^r, T^r) &=&
\sum_{n=0}^{p_0} \biggl[ \E\int_{[b^r(n)/r^2, a^r(n+1)/r^2)}
e^{-\gamma t} \bigl(h\cdot\hat Q^r(t) + \gamma p\cdot\hat U^r(t)\bigr) \,dt
\nonumber\\[-2pt]
&&\hspace*{18pt}{}+ \E\int_{[a^r(n)/r^2, b^r(n)/r^2)}
e^{-\gamma t} \bigl(h\cdot\hat Q^r(t) + \gamma p\cdot\hat U^r(t)\bigr) \,dt
\biggr] \\[-2pt]
&=&
\sum_{n=0}^{p_0} \E\int_{[b^r(n)/r^2, a^r(n+1)/r^2)}
e^{-\gamma t} \bigl(h\cdot\hat Q^r(t) + \gamma p\cdot\hat U^r(t)\bigr) \,dt
+\varepsilon_r,\nonumber\vadjust{\goodbreak}
\end{eqnarray}
where, using Proposition \ref{tight424} and the observation that
$a^r(n+1)/r^2- b^r(n)/\allowbreak r^2 \le\rho/r \to0$,
we have that $\varepsilon_r \to0$ as $r\to\infty$. %
From Theorem \ref{mainweak}, as $r\to\infty$,
$\curvq^r \Rightarrow\curvq$. Combining this with Proposition \ref
{tight424},
we get for every $n = 0, \ldots,p_0$,
%
\begin{eqnarray}\label{ab753}
&&\lim_{r\to\infty} \E\int_{[b^r(n)/r^2, a^r(n+1)/r^2)}
e^{-\gamma t} h\cdot\hat Q^r(t) \,dt\nonumber\\
&&\qquad= \lim_{r\to\infty} \E\int_{[n\theta+\rho/r, (n+1)\theta)}
e^{-\gamma t} h\cdot\hat Q^r(t) \,dt\nonumber\\
&&\qquad= \lim_{r\to\infty} \E\int_{[n\theta, (n+1)\theta)}
e^{-\gamma t} h\cdot\curvq^{r,n}(t-n\theta) \,dt\\
&&\qquad= \E\int_{[n\theta, (n+1)\theta)}
e^{-\gamma t} h\cdot\curvq^{(n)}(t-n\theta) \,dt\nonumber\\
&&\qquad= \E\int_{[n\theta, (n+1)\theta)}
e^{-\gamma t} h\cdot\tilde Q(t) \,dt,\nonumber
\end{eqnarray}
where, by convention, $[n\theta, (n+1)\theta)= [p_0\theta, \infty)$
when $n = p_0$.
Next, for $t \in[n\theta+\rho/r, (n+1)\theta)$, $n = 0, \ldots, p_0$,
%
\begin{eqnarray}\label{ab1945}
\qquad \gamma p \cdot\hat U^r(t)
&=& \gamma p \cdot K\bigl(\hat Y^r(t) - \hat Y^r(n\theta+\rho/r)\bigr) +
\gamma p \cdot K\hat Y^r(n\theta+\rho/r)
\nonumber
\\[-8pt]
\\[-8pt]
\nonumber
&=& \gamma p \cdot K \curvz^{r,n}(t-n\theta) + \gamma p \cdot K\nu
_0^{r,n}.
\end{eqnarray}
From Theorem \ref{mainweak}, as $r \to\infty$,
\[
\gamma p \cdot K \curvz^{r,n}(\cdot) + \gamma p \cdot K\nu_0^{r,n}
\Rightarrow\gamma p \cdot K \curvz^{(n)}(\cdot) + \gamma p \cdot
K\nu_0^{(n)}
\]
in $\mathcal{C}_{\theta}^1$.
Combining this with \eqref{ab1945} and Proposition \ref{tight424} we
now get similarly to~\eqref{ab753}, for $n = 0, \ldots,p_0$,
\begin{eqnarray*}\lim_{r \to\infty}
&&\E\int_{[b^r(n)/r^2, a^r(n+1)/r^2)}
\gamma e^{-\gamma t} p\cdot\hat U^r(t) \,dt\\
&&\qquad=\lim_{r \to\infty} E \int_{[n\theta, (n+1)\theta)}
\gamma e^{-\gamma t} p \cdot K \bigl( \curvz^{r,n}(t-n\theta) +\nu
_0^{r,n}\bigr) \,dt\\
&&\qquad=E \int_{[n\theta, (n+1)\theta)}
\gamma e^{-\gamma t} p \cdot K \bigl(\curvz^{(n)}(t-n\theta) + \nu
_0^{(n)}\bigr) \,dt.
\end{eqnarray*}
Note that for $t \in[n\theta, (n+1)\theta)$,
\[
\curvz^{(n)}(t-n\theta)+\nu_0^{(n)}= \tilde Y(n\theta) + \tilde
Y(t) - \tilde
Y(n\theta) = \tilde Y(t).
\]
Thus, the expression on the right-hand side of the above display equals
\[
\E\int_{[n\theta, (n+1)\theta)}
\gamma e^{-\gamma t} p\cdot\tilde U(t) \,dt.
\]
The result now follows on using this observation along with \eqref
{ab753} in \eqref{ab804}.
\end{pf*}

\subsection{\texorpdfstring{Proof of Theorem \protect\ref{mainweak}}{Proof of Theorem 4.5}}
\label{weakcgce}
For $j \in\mathbb{J}$, $n = 0, 1, \ldots, p_0$, and $\omega\in\Omega
$, define
%
\begin{equation}\label{lab107}\check S_j^{r,n}(\omega) = \{s \in[0,
\rho]\dvtx
Q^r_{\sigma_1(j)}(nr^2\theta+ rs, \omega) = 0\}.
\end{equation}
From
the definition of $\vartheta_{\varepsilon_0}$, it follows that for some
$\gamma_0 > 0$,
\[
\inf_{r \ge r_0} \min_{n=0, \ldots, p_0; j \in\mathbb{J}} (R\nu
^{r,n})_j \ge\gamma_0\qquad \mbox{a.e.}
\]
As a consequence of this observation, we have the following result.
The proof is given in Section~\ref{prop46}.
\begin{prop}\label{lemma45a}
For some $\{\rho_r\} \subset[0, \rho]$ such that $\rho_r \to0$ as $r
\to\infty$, we have
%
\begin{equation}
\label{claim1022}
\IP(\Psi^r_n) \to1 \qquad \mbox{as } r
\to\infty, \mbox{ for all } n =0, 1, \ldots, p_0,
\end{equation}
where $\Psi^r_n = \{\omega\in\Omega\dvtx  (\bigcup_{j\in\mathbb
{J}}\check
S_j^{r,n}(\omega)) \cap[\rho_r, \rho] =\varnothing\}$.
\end{prop}
For $n=0, 1, \ldots, p_0$, let $\bar\nu^r[n]=(\bar\nu
^{r,0},\ldots,\bar\nu^{r,n})$. We define
$\nu^r[n], \nu_0^r[n], \curvq^r[n],\allowbreak  \curvz^r[n]$ and their
limiting analogues $\bar\nu[n],\nu[n], \nu_0[n], \curvq[n],
\curvz[n]$ in a similar fashion.
Set
\begin{eqnarray*}
\mathcal{J}^r[n] &= &(\hat\zeta^r, \bar\nu^r[n], \nu^r[n], \nu
_0^r[n], \curvq^r[n], \curvz^r[n]),\\
\mathcal{J}[n] &=& (\hat\zeta, \bar\nu[n], \nu[n], \nu_0[n],
\curvq
[n], \curvz[n]).
\end{eqnarray*}
Then $\mathcal{J}^r[n], \mathcal{J}[n]$ are $\Xi[n]$-valued random
variables, with
\[
\Xi[n] = \mathcal{D}^{\II} \times(\mathcal{S}_M^{\eta})^{\otimes(n+1)}
\times\bigl(\R^{\JJ(n+1)}\bigr)
\times\bigl(\R^{\JJ(n+1)}\bigr) \times\mathcal{D}_{\theta}^{(n+1)\II}
\times\mathcal{C}_{\theta}^{(n+1)\JJ},
\]
where we follow the usual convention for $n=p_0$.
To prove the
theorem, we need to show that $\mathcal{J}^r[p_0]\Rightarrow\mathcal
{J}[p_0]$. In
the lemma below we will in fact show, recursively in $n$, that
$\mathcal{J}^r[n]\Rightarrow\mathcal{J}[n]$ as $r \to\infty$, for
each $n=0, 1,
\ldots, p_0$, which will complete the proof of Theorem \ref{mainweak}.
\begin{lemma}
\label{riszero} For each $n = 0, 1, \ldots, p_0$, $\mathcal{J}
^r[n]\Rightarrow
\mathcal{J}[n]$, as $r \to\infty$.
\end{lemma}
\begin{pf}
The proof will follow the following two steps:
\begin{longlist}[(ii)]
\item[(i)] As $r \to\infty$, $\mathcal{J}^r[0]\Rightarrow\mathcal{J}[0]$.

\item[(ii)] Suppose that $\mathcal{J}^r[k]\Rightarrow\mathcal{J}[k]$ as $r
\to\infty$ for
$k=0, 1,\ldots,n$, for some $n < p_0$. Then, as $r \to\infty$,
$\mathcal{J}^r[n+1]\Rightarrow\mathcal{J}[n+1]$.
\end{longlist}

Consider (i). Define scaled processes $\check
Q^r(t) = Q^r(rt)/r$, $\check Y^r(t)= Y^r(rt)/r$. Processes
$\check X^r$, $\check T^r$, $\check T^{r,1}$, $\check\zeta^r$ are defined
similarly.
By the functional central limit theorem for renewal
processes and Proposition \ref{tisadmis}, it follows that (cf. Lemma
3.3 of \cite{BudGho2})
%
\begin{equation}\label{ab1122}\hat\zeta^r
\Rightarrow\tilde\zeta.
\end{equation}
Also, convergence of $(\bar\nu^r[0],\nu
^r[0])$ follows trivially since $\bar\nu^r[0] =\bar
\nu[0]=0$ and $\nu^r[0] =\nu[0]=\varepsilon_0 y^*$. Next,
consider\vadjust{\goodbreak}
$\nu_0^r[0] = \hat Y^r(\rho/r)$. From the definition of the scaled
processes defined above \eqref{lab107}, we have that
%
\begin{equation}\label{abnew558}
\check\zeta^r(t) = \check q^r + \check X^r(t) + \frac{1}{r}[\theta_1^r
t - (C-P')
\operatorname{diag}(\theta_2^r) \check
T^r(t)],
\end{equation}
and
%
\begin{equation}\label{abnew558b}\check Q^r(t) = \check\zeta^r(t) +
R \check
Y^r(t), \qquad t \in[0, \rho].
\end{equation}
Also, observe that $\check Y^r$ can be written as
%
\begin{equation}
\label{abnew558c}\check Y^r(t) = \biggl(\biggl(x^*-\frac{\varepsilon_0y^*}{\rho}\biggr)t -
\check T^{r,1}(t)\biggr) +
\frac{1}{\rho} \varepsilon_0y^*t + \check N^r(t),
\end{equation}
where, with $\check S_j^{r,0}$ defined in \eqref{lab107},
%
\begin{equation}
\label{abnew558d}\check N_j^r(t) = \int_{[0,t]}1_{\check S_j^{r,0}}(s)
\,d \check T^{r,1}_j(s).
\end{equation}
Next, note that, for a suitable $C_1 \in(0, \infty)$,
%
\begin{equation}\biggl|\check T^{r,1} - \biggl(x^*-\frac{\varepsilon_0y^*}{\rho}\biggr)
\curvi
\biggr|_{\infty, \rho} \le C_1 \frac{\Delta^r}{r}.\label{ab1026}
\end{equation}
Also,
$|\check N^r|_{\infty,\rho} \le\rho_r + \rho1_{(\Psi^r_n)^c}.$
Thus, from Proposition \ref{lemma45a}, $|\check N^r|_{\infty,\rho}$
converges to $0$ in probability as $r \to\infty$, which shows that
%
\begin{eqnarray}\label{442half}
\biggl|\check Y^r - \frac{\varepsilon
_0}{\rho} y^*
\biggr|_{\infty, \rho}= \biggl|\biggl(x^*-\frac{\varepsilon_0}{\rho} y^*\biggr)\curvi-
\check
T^r\biggr|_{\infty, \rho}
\to0
\nonumber
\\[-8pt]
\\[-8pt]
\eqntext{ \mbox{in probability, as } r \to\infty.}
\end{eqnarray}
The above convergence is the key reason for introducing the
modification of~$\tilde Y^{(1)}$, through the vector $y^*$, described
above Theorem
\ref{newmain518}.

Next, standard moment bounds for renewal processes (see, e.g.,
Lemma~3.5 of~\cite{BudGho2}) yield that $|\check X^r|_{\infty, \rho}$
converges to
zero in probability as $r \to\infty$. Combining these observations, we
get from \eqref{abnew558} and \eqref{abnew558b} that
$(\check Q^r,\check Y^r)$ converge, uniformly over $[0, \rho]$, in
probability, to $(q+\frac{\varepsilon_0}{\rho}Ry^*\curvi,
\frac{\varepsilon_0y^*}{\rho}\curvi)$.
In particular, this shows that
%
\begin{eqnarray}\label{ab1042}
&&
(\curvq^{r,0}(0), \nu_0^{r,0})=(\check
Q^r(\rho),
\check Y^r(\rho))
\nonumber
\\[-8pt]
\\[-8pt]
\nonumber
&&\quad \Rightarrow\quad(q+\varepsilon_0Ry^*, \varepsilon_0y^*)
= \bigl(\curvq^{(0)}(0), \nu_0^{(0)}\bigr).
\end{eqnarray}
Finally, we prove the convergence of $(\curvq^{r,0}, \curvz^{r,0})$ to
$(\curvq^{(0)}, \curvz^{(0)})$.
We will apply Theorem 4.1 of \cite{will-invpr}. Note that
%
\begin{eqnarray}\label{lab99}\curvq^{r,0}(t) = \curvq^{r,0}(0) +
w^{r,0}(t) +
R\curvz^{r,0}(t), \qquad t \in[0, \theta],
\end{eqnarray}
where $w^{r,0}$ is a $\mathcal{D}_{\theta}^{\II}$-valued random variable
defined as $w^{r,0}(t) = (\hat\zeta^r(t)-\hat\zeta^r(\rho
/r))1_{[\rho
/r, \infty)}(t)$.
From \eqref{ab1042} and \eqref{ab1122}
%
\begin{equation}
\label{ab430} \curvq^{r,0}(0)\to\curvq^{(0)}(0)\quad \mbox{and}\quad
w^{r,0} \Rightarrow\tilde\zeta \qquad\mbox{as } r \to\infty.\vadjust{\goodbreak}
\end{equation}
Next, for $t \in[\rho/r, \theta)$,
write
%
\begin{equation}
\label{ab441half}
\curvz^{r,0}(t) = \hat H^{r,0}(t) + \hat L^{r,0}(t),
\end{equation}
where for $j \in\mathbb{J}$,
\begin{eqnarray*}
\hat H^{r,0}_j(t) &=& r\int_{[\rho/r,t]}\bigl(1-1_{\mathcal{S}
_j^{r,1}}(s)\bigr)\, d\bigl(x^*_js -\bar T^{r,1}_j(s) \bigr)
+ r \int_{[\rho/r,t]}1_{\mathcal{S}_j^{r,2}}(s) \,d\bar T^{r,1}_j(s),\\
\hat L^{r,0}_j(t) &=& rx^*_j\int_{[\rho/r,t]}1_{\mathcal{S}_j^{r,1}}(s)\,ds,
\end{eqnarray*}
with $\mathcal{S}_j^{r,i}$ defined in \eqref{lab33}.
Using calculations similar to those in the proof of Proposition \ref
{tisadmis} 
[see \eqref{lab71}], we get
%
\begin{equation}
\sup_{\rho/r \le t \le\theta}\biggl| r\int_{[\rho/r, t]} 1_{\mathcal{S}
_j^{r,2}} (s) \bar T^{r,1}_j(s)\biggr|
\to0\qquad \mbox{in probability, as } r \to\infty,\hspace*{-35pt}
\end{equation}
for all $j \in\mathbb{J}$.

Also, from \eqref{ab248} it follows that
\[
\sup_{\rho/r \le t \le\theta}\biggl| r\int_{[\rho/r, t]} \bigl(1 -
1_{\mathcal{S}_j^{r,1}} (s)\bigr)\, d\bigl(x^*_js- \bar T^{r,1}_j(s)\bigr)\biggr| \le
\frac{\Delta^r}{r}.
\]
Combining the above estimates,
%
\begin{equation}\label{ab1211} \sup_{\rho/r \le t \le\theta}|\hat
H_j^{r,0}(t)|
\to0\qquad \mbox{in probability, as } r \to\infty.
\end{equation}
Also,
$\hat L^{r,0}(t)=(\operatorname{diag}(x^*))C'\hat\ell^{r,0}(t),$
where for $i \in\mathbb{I}$ and $t \in[\rho/r,\theta]$,
%
\begin{eqnarray}\label{ab1452n}
\hat\ell_i^{r,0}(t)&=& r \int_{[\rho/r,t]} 1_{\{\hat Q^r_i(\bar
\curvp
^r(s)) \le r\bar\Theta^r(s)\}} \,ds
\nonumber
\\[-8pt]
\\[-8pt]
\nonumber
&=&
r \int_{[\rho/r,t]} 1_{\{\curvq^{r,0}_i(\bar\curvp^r(s)) \le r\bar
\Theta^r(s)\}} \,ds.
\end{eqnarray}
Recall that $\curvz^{r,0}(t) = 0$ for $t \in[0, \rho/r]$. Hence,
setting $\hat H^{r,0}(t)= \hat\ell^{r,0}(t)=0$ for $t \in[0, \rho
/r]$, we have
from \eqref{lab99} and \eqref{ab441half}
%
\begin{equation}\label{ab444}\qquad
\curvq^{r,0}(t) = \curvq^{r,0}(0) +
w^{r,0}(t) +
R\hat H^{r,0}(t) + D\hat\ell^{r,0}(t), \quad t \in[0, \theta].
\end{equation}
From \eqref{ab430} and \eqref{ab1211} we now have that, as $r \to
\infty$,
%
\begin{equation}\label{ab445}
q^{r,0}(0) + w^{r,0} + R\hat H^{r,0} \Rightarrow q ^{(0)}+ \tilde\zeta
\end{equation}
in $\mathcal{D}_{\theta}^{\II}$.
Using the definition of $\curvp^r$,
Assumption \ref{ldp} and elementary properties of renewal processes
[see similar arguments in \eqref{lab71} and \eqref{lab72}], we have
that for some $C_2$, as $r \to\infty$,
\[
\IP\Bigl(\sup_{s \in[\rho/r,\theta]}|\hat Q^r(\bar\curvp^r(s)) - \hat
Q^r(s)| > \varepsilon\Bigr)
\le C_2\frac{r^2\theta}{\Delta^r} \frac{1}{(\Delta^r)^{\curvm}} = C_2
\frac{r^2\theta}{r^{\curvk(\curvm+1)}} \to0.
\]
This shows that
%
\begin{equation}
\label{ab426}\qquad
\sup_{s \in[0,\theta]}|\curvq^{r,0}(s) - \curvq^{r,0}(\bar\curvp
^r(s))| \to0\quad \mbox{in probability, as } r \to\infty.
\end{equation}
Using Theorem 4.1 of \cite{will-invpr} along with \eqref{ab1452n},
\eqref{ab444}, \eqref{ab445} and \eqref{ab426}, we now have that
\[
(\hat\zeta^r, \curvq^{r,0}, \hat\ell^{r,0}) \Rightarrow\bigl(\tilde
\zeta, \Gamma\bigl(q^{(0)}(0) +\tilde\zeta\bigr),
\hat\Gamma\bigl(q^{(0)}(0) +\tilde\zeta\bigr)\bigr)\qquad \mbox{as } r \to\infty,
\]
as $\mathcal{D}^{\II}\times\mathcal{D}_{\theta}^{\II}\times
\mathcal{D}_{\theta
}^{\II}$
valued random variables.
Since
\[
\Gamma\bigl(q^{(0)}(0) +\tilde\zeta\bigr) = \curvq^{(0)} \quad\mbox{and}\quad
\operatorname{diag}(x^*) C'\hat\Gamma\bigl(q^{(0)}(0) +\tilde\zeta\bigr) =
\curvz^{(0)},
\]
we get from the above display that $(\curvq^{r,0}, \curvz^{r,0})
\Rightarrow(\curvq^{(0)}, \curvz^{(0)})$ as $r \rightarrow\infty$.
Combining this with \eqref{ab1042} and observations below \eqref{ab1122},
we have $\mathcal{J}^r[0] \Rightarrow\mathcal{J}[0]$,
which completes the proof of (i).

We now prove (ii). We can write
\begin{eqnarray*}\label{conti}\mathcal{J}^r[n+1] &= &(\mathcal
{J}^r[n], (\bar\nu
^{r,n+1}, \nu^{r,n+1}, \nu_0^{r,n+1}, \curvq^{r,n+1},
\curvz^{r,n+1}))\\
\mathcal{J}[n+1] &=& \bigl(\mathcal{J}^r[n],\bigl(\bar\nu^{(n+1)}, \nu^{(n+1)},
\nu_0^{(n+1)},
\curvq^{(n+1)}, \curvz^{(n+1)}\bigr)\bigr).
\end{eqnarray*}
By assumption,
%
\begin{equation}\label{ab505}\mathcal{J}^r[n]
\Rightarrow\mathcal{J}[n]
\end{equation}
and, thus, in particular, \eqref{ab1122} holds.
This shows that $\hat\chi^{\kappa,r}(n+1) \Rightarrow
\chi^{\kappa}(n+1)$ and as a consequence, using continuity properties of
$\varpi_{n+1}$,
%
\begin{eqnarray}\label{ab459n}\qquad
\bar\nu^{r,n+1} &\Df&\varpi_{n+1}\bigl(\hat\chi^{\kappa,r}(n+1),
\mathcal{U}_{n+1}\bigr)
\Rightarrow\varpi_{n+1}\bigl(\chi^{\kappa}(n+1), \tilde\mathcal
{U}_{n+1}\bigr)
\nonumber
\\[-8pt]
\\[-8pt]
\nonumber
& =&
\bar\nu
^{(n+1)}.
\end{eqnarray}
In fact, this shows the joint convergence: $(\mathcal{J}^r[n], \bar
\nu
^{r,n+1}) \Rightarrow(\mathcal{J}[n], \bar\nu^{n+1})$.
In particular, we have
\begin{eqnarray*}
&&\bigl(\bar\nu^{r,n+1}, \hat Q^r\bigl((n+1)\theta\bigr)\bigr)=(\bar\nu^{r,n+1},\curvq
^{r,n}(\theta))\\
&&\quad\Rightarrow\quad\bigl(\bar\nu^{(n+1)},\curvq^{(n)}(\theta)\bigr)
=\bigl(\bar\nu
^{(n+1)}, \tilde Q((n+1)\theta-)\bigr).
\end{eqnarray*}
For the remaining proof, to keep the presentation simple, we will not
explicitly note the joint convergence of all
the processes being considered.
From continuity of the map $\vartheta_{\varepsilon_0}$, we now have that
\[
\nu^{r,n+1} = \vartheta_{\varepsilon_0}\bigl(\hat Q^r\bigl((n+1)\theta\bigr), \bar
\nu
^{r,n+1}\bigr) \quad\Rightarrow\quad
\vartheta_{\varepsilon_0}\bigl(\tilde Q\bigl((n+1)\theta-\bigr), \bar\nu^{(n+1)}\bigr)
=\nu^{(n+1)}.
\]
Next, we consider the weak convergence of $\nu_0^{r,n+1}$ to $\nu
_0^{(n+1)}$. The proof is similar to the case $n+1=0$
treated in the first part of the lemma [cf. below~\eqref{ab1122}] and
so only a sketch will be provided.
Note that
\begin{eqnarray*}
\nu_0^{r,n+1} &=& \nu_0^{r,n} + \curvz^{r,n}(\theta) + \bigl(\hat
Y^r\bigl((n+1)\theta+ \rho/r\bigr) - \hat Y^r\bigl((n+1)\theta\bigr)\bigr)\\
\curvq^{r,n+1}(0) &=& \curvq^{r,n}(\theta) + \bigl(\hat Q^r\bigl((n+1)\theta+
\rho/r\bigr) - \hat Q^r\bigl((n+1)\theta\bigr)\bigr).
\end{eqnarray*}
Weak convergence of
$(\curvq^{r,n}(\theta), \nu_0^{r,n} + \curvz^{r,n}(\theta))$ to
$(\curvq^{(n)}(\theta), \nu_0^{(n)} + \curvz^{(n)}(\theta))$ is
a~consequence of \eqref{ab505}.
Next, abusing notation introduced above \eqref{ab1122}, define for $t
\in[0, \rho]$,
\begin{eqnarray*}
\check Q^r(t) &=& r^{-1}Q^r\bigl(r^2\theta(n+1) + rt\bigr), \qquad \check q^r = \check
Q^r(0),\\
\check Y^r(t) &= & r^{-1}\bigl(Y^r\bigl(r^2\theta(n+1) + rt\bigr) - Y^r\bigl(r^2\theta
(n+1)\bigr)\bigr).
\end{eqnarray*}
Processes $\check X^r, \check T^r, \check T^{r,1}, \check\zeta^r$ are
defined similarly to $\check Y^r$.
Then, equations \eqref{abnew558} and \eqref{abnew558b} are satisfied with
these new definitions.
Hence, using arguments similar to the ones used in the proof of \eqref
{ab1042} (in particular, making use of Proposition \ref{lemma45a}), we
have that
$(\check Y^r, \check Q^r)$ converges in distribution to
\[
\biggl( \frac{\nu^{(n+1)}}{\rho}\curvi, \tilde Q\bigl((n+1)\theta-\bigr) + \frac
{1}{\rho}R\nu^{(n+1)}\curvi\biggr),
\]
as $r \to\infty$. Combining the above observations, we have, as $r
\to
\infty$,
%
\begin{eqnarray}
\label{ab604}
&&(\curvq^{r,n+1}(0), \nu_0^{r,n+1}) = (\check Q^r(\rho), \check
Y^r(\rho))
\nonumber
\\[-8pt]
\\[-8pt]
\nonumber
&&\quad\Rightarrow\quad\bigl(\tilde Q\bigl((n+1)\theta\bigr), \nu_0^{(n)}+ \curvz
^{(n)}(\theta)+\nu^{(n+1)}\bigr)=\bigl(\curvq^{(n+1)}(0),
\nu_0^{(n+1)}\bigr).\nonumber\hspace*{-30pt}
\end{eqnarray}
Finally, we consider weak convergence of $(\curvq^{r,n+1}, \curvz
^{r,n+1})$ to
$(\curvq^{(n+1)}, \curvz^{(n+1)})$.
Similar to \eqref{lab99}, we have
\[
\curvq^{r,n+1}(t) = \curvq^{r,n+1}(0) + w^{r,n+1}(t) + R\curvz
^{r,n+1}(t),\qquad t \in[0, \theta],
\]
where $w^{r,n+1}$ is a $\mathcal{D}_{\theta}^{\II}$-valued random variable
defined as
\[
w^{r,n+1}(t) = \bigl(\hat\zeta^r\bigl(t+(n+1)\theta\bigr)-\hat\zeta
^r\bigl((n+1)\theta+ \rho/r\bigr)\bigr)1_{[\rho/r, \infty)}(t).
\]
Using \eqref{ab604} and \eqref{ab1122}, as $r \to\infty$,
%
\begin{equation}\qquad
\label{ab430new} \curvq^{r,n+1}(0) + w^{r,n+1} \Rightarrow\curvq^{(n+1)}(0)
+ \tilde{\zeta}\bigl((n+1)\theta+\cdot\bigr) -\tilde\zeta
\bigl((n+1)\theta\bigr) .
\end{equation}
Weak convergence of $(\curvq^{r,n+1}, \curvz^{r,n+1})$ to
$(\curvq^{(n+1)}, \curvz^{(n+1)})$ now follows exactly as below
\eqref
{ab430}. Combining the above weak convergence properties,
we now have $\mathcal{J}^r[n+1] \Rightarrow\mathcal{J}[n+1]$ and the
result follows.
\end{pf}

\subsection{\texorpdfstring{Proof of Proposition \protect\ref{lemma45a}}{Proof of Proposition 4.6}} \label{prop46}
We will only consider the case $n=0$. The general case is treated similarly.
Let $M_r = \lfloor\frac{r\rho}{\Delta^r}\rfloor$. From
Assumption \ref{ldp}, for each $\delta> 0$,
one can find $C_1(\delta)$ such that, for $i \in\mathbb{I}$, $j \in
\mathbb{J}
$, $r \ge1$ and $k \le M_r$,
%
\begin{eqnarray}\label{lab198}
&\displaystyle\IP\bigl(|E_i^r\bigl((k\,{+}\,1)\Delta^r\bigr)\,{-}\,E_i^r(k\Delta^r)\,{-}\,
    \alpha_i^r \Delta^r\bigr|\,{\ge}\,\delta\Delta^r\bigr)\,{\le}\,
    \frac{C_1(\delta)}{r^{\curvk\curvm}},\hspace*{-30pt}&\nonumber\\
&\displaystyle\IP\bigl(\bigl|S_j^r\bigl(T_j^{r,1}\bigl((k\,{+}\,1)\Delta^r\bigr)\bigr)\,{-}\,
    S_j^r(T_j^{r,1}(k\Delta^r))\,{-}\,\beta_j^r\tau_j^{r,k}\bigr|\,{\ge}\,\delta\Delta^r\bigr)
    \,{\le}\, \frac{C_1(\delta)}{r^{\curvk\curvm}},\hspace*{-30pt}&\nonumber\\[-8pt]\\[-8pt]
\qquad&\displaystyle\IP\bigl(\bigl|\Phi_i^{j,r}\bigl(S_j^r\bigl(T_j^{r,1}\bigl((k\,{+}\,
    1)\Delta^r\bigr)\bigr)\bigr)\,{-}\,\Phi_i^{j,r}(S_j^r(T_j^{r,1}(k\Delta^r)))\,{-}\,
    p^j_i\beta_j^r \tau_j^{r,k}\bigr|\,{\ge}\,\delta\Delta^r\bigr)\hspace*{-30pt}&\nonumber\\
&\displaystyle\hspace*{-242pt}\,{\le}\,
    \frac{C_1(\delta)}{r^{\curvk\curvm}},\hspace*{-30pt}&\nonumber
\end{eqnarray}
where $\tau_j^{r,k}=T_j^{r,1}((k+1)\Delta^r)-T_j^{r,1}(k\Delta^r) =
(x^*_j-\frac{\varepsilon_0 y^*_j}{\rho})\Delta^r.$
Denote the union, over all $i,j$, of events on the left-hand side of the
three displays in \eqref{lab198}, by $H_k^r$. Then, the above
estimates, along with
\eqref{lab23}, yield
%
\begin{equation}
\label{ab1140}
\IP\Biggl(\bigcup_{k=0}^{M_r} H_k^r\Biggr) \to0 \qquad\mbox{as } r \to\infty.
\end{equation}
Define for $k=0,1,\ldots, M_r$,
\begin{eqnarray*}
Q^r_0\bigl((k+1)\Delta^r\bigr) &\Df& Q^r(k\Delta^r) + E^r\bigl((k+1)\Delta^r\bigr) -
E^r((k\Delta^r))\\[-2pt]
&&{}- C \bigl(S^r\bigl(T^{r,1}\bigl((k+1)\Delta^r\bigr)\bigr)- S^r(T^{r,1}(k\Delta^r))\bigr)
\\[-2pt]
&&{}+ \Phi^r\bigl(S^r\bigl(T^{r,1}\bigl((k+1)\Delta^r\bigr)\bigr)\bigr)-
\Phi^{r}(S^r(T^{r,1}(k\Delta^r))).
\end{eqnarray*}
Note that, on the set $(H_k^r)^c$, we have for some $C_2>0$,
\begin{eqnarray*}
&&Q^r_0\bigl((k+1)\Delta^r\bigr)\\[-2pt]
&&\qquad\ge Q^r(k\Delta^r) + \alpha^r\Delta^r -
(C-P')\operatorname{diag}(\beta^r)\biggl(x^*- \frac{\varepsilon_0y^*}{\rho
}\biggr) \Delta^r -
C_2\delta
\Delta^r{\mathbf{1}}_{\II} .
\end{eqnarray*}
Using Assumption \ref{assum-limit-param}, we now have that for some
$C_3$, on the set $(H_k^r)^c$, for all $r \ge r_0$,
\[
Q^r_0\bigl((k+1)\Delta^r\bigr) \ge Q^r(k\Delta^r) +\frac{\varepsilon_0}{\rho}
Ry^*\Delta
^r - \biggl(C_2\delta\Delta^r + C_3 \frac{\Delta^r}{r}\biggr){\mathbf{1}}_{\II}.
\]
Recall that $Ry^* > \gamma_0 {\mathbf{1}}_{\II}$. Fix $\delta$ small enough
so that for some $\varepsilon_1 > 0$ and $r_1 > r_0$,
%
\begin{equation}\label{ab1432} \frac{\varepsilon_0}{\rho} \gamma_0 -
\biggl(C_2\delta+ \frac
{C_3}{r}\biggr) \ge\varepsilon_1 \qquad\mbox{for all }
r \ge r_1.
\end{equation}
Then, for every $k =0, 1, \ldots, M_r$,
%
\begin{eqnarray}
\label{ab1202}
\mbox{on the set } (H_k^r)^c,\quad Q^r\bigl((k+1)\Delta^r\bigr) &\ge&
Q^r_0\bigl((k+1)\Delta^r\bigr)
\nonumber
\\[-9.5pt]
\\[-9.5pt]
\nonumber
&\ge &Q^r(k\Delta^r) + \varepsilon_1 {\mathbf{1}}_{\II}
\Delta
^r,\quad
r \ge r_1.\hspace*{-30pt}
\end{eqnarray}
Recall $d_1$ introduced above \eqref{ab411}. Let $m_0$ be large enough
so that $\varepsilon_1 m_0
> d_1$.
Then using \eqref{ab1202}, we get that 
%
\begin{eqnarray}
\label{ab1202a}
\mbox{on the set } \bigcap_{k=0}^{M_r}(H_k^r)^c,\quad
Q^r(k\Delta^r)
\ge d_1 {\mathbf{1}}_{\II} \Delta^r,\\[-2pt]
\eqntext{ \mbox{for all } k =m_0, \ldots, M_r,
r \ge r_1.}
\end{eqnarray}
Next, let
%
\begin{eqnarray}\label{frk}
F_k^r & = & \Bigl\{\omega\dvtx  \inf_{i \in
\mathbb{I}} Q_i^r(k\Delta
^r, \omega) \ge d_1\Delta^r\Bigr\}\nonumber\\
&&{} \cap
\{\omega\dvtx  Q_i^r(t, \omega) = 0 \mbox{ for some } i \in\mathbb{I},
t \in[k\Delta^r,
(k+1)\Delta^r]\}\\
& \Df& G_k^r \cap B_k^r .\nonumber\vadjust{\goodbreak}
\end{eqnarray}
Using estimates
below \eqref{ab1249}, we see
that
%
\begin{equation}\label{ab1208n} \IP\Biggl(\bigcup_{k=0}^{M_r} F_k^r\Biggr) \to0
\qquad\mbox{as } r \to\infty.\vspace*{-2pt}
\end{equation}
%
Also, from \eqref{ab1202a},
%
\begin{eqnarray}\label{ab565half} \liminf_{r\to\infty} \IP
(G_{m_0}^r)& \ge&\liminf
_{r\to\infty}
\IP\Biggl( G_{m_0}^r \cap\biggl[ \bigcap_{k=0}^{M_r}(H_k^r)^c
\biggr]\Biggr)
\nonumber
\\[-10pt]
\\[-10pt]
\nonumber
&= &\liminf_{r\to\infty}
\IP\Biggl( \bigcap_{k=0}^{M_r}(H_k^r)^c\Biggr) = 1.\vspace*{-2pt}
\end{eqnarray}
Next, for $ r \ge r_1$,
%
\begin{equation}
\label{ab1235n}
\IP\Biggl( \bigcup_{k=m_0}^{M_r}B_k^r \Biggr) \le\IP\Biggl( \bigcup
_{k=m_0}^{M_r}(B_k^r \cap G^r_{m_0})\Biggr)
+ \IP((G^r_{m_0})^c).\vspace*{-2pt}
\end{equation}
Also,
%
\begin{eqnarray}\label{new1}\IP\Biggl( \bigcup_{k=m_0}^{M_r}(B_k^r \cap
G^r_{m_0})\Biggr) &=&
\IP\Biggl(\Biggl [\bigcup_{k=m_0}^{M_r}(B_k^r \cap G^r_{m_0})\Biggr]
\cap\Biggl[ \bigcup_{k=0}^{M_r}H_k^r\Biggr]\Biggr)\nonumber\\[-3pt]
&&{}+\IP\Biggl( \Biggl[\bigcup_{k=m_0}^{M_r}(B_k^r \cap G^r_{m_0})\Biggr]
\cap\Biggl[ \bigcap_{k=0}^{M_r}(H_k^r)^c\Biggr]\Biggr)\\[-3pt]
&\le&
\IP\Biggl( \bigcup_{k=0}^{M_r}H_k^r\Biggr)+\IP\Biggl( \bigcup
_{k=m_0}^{M_r}(B_k^r \cap G^r_{k})\Biggr),\nonumber\vspace*{-2pt}
\end{eqnarray}
where the last inequality is a consequence of the fact that on $\bigcap
_{k=0}^{M_r}(H_k^r)^c$, $G^r_k \subseteq G^r_{k+1}$ for $k\ge m_0$.
From \eqref{ab1140} and \eqref{ab1208n} the above expression is seen to
approach zero as $r \to\infty$.
Using this observation and \eqref{ab565half} in \eqref{ab1235n}, we now
see that $\IP( \bigcup_{k=m_0}^{M_r}B_k^r ) \to0$
as $r \to\infty$.
Finally, recalling the definition of~$B_k^r$, we have
\[
\IP\bigl(Q_i^r(s) = 0, \mbox{ for some } i \in\mathbb{I}\mbox{ and } s
\in
[m_0\Delta^r, r\rho]\bigr) = \IP\Biggl(\bigcup_{k=m_0}^{M^r}B_k^r\Biggr).\vspace*{-2pt}
\]
The proposition now follows on setting $\rho_r = \frac{m_0 \Delta
^r}{r}$. \vspace*{-3pt}

\begin{appendix}
\section*{Appendix}\vspace*{-3pt}\label{appen}

\renewcommand{\thelemm}{\Alph{section}.\arabic{lemm}}
\begin{lemm} \label{appnov}
Let $\{\tilde Y_n\}_{n \ge1}$ be a sequence of random variables, with
values in a finite set $\mathbb{S}$, given on a probability
space
$(\Omega, \mathcal{F}, \IP)$. Let $\mathcal{G}$ be a sub-$\sigma$
field of $\mathcal{F}$.
Suppose that $\{\mathcal{G}_n\}_{n \ge1}$
is a sequence of sub-$\sigma$ fields of $\mathcal{G}$ and $\{X_n\}
_{n\ge1}$ a
sequence of $\{\mathcal{G}_n\}$-adapted, $\R^d$-valued random
variables such that
\[
\IP(\tilde Y_n = \zeta\vert\mathcal{G}_n) = p_{n,\zeta}(X_n),\qquad n \ge1,
\zeta
\in\mathbb{S},\vspace*{-2pt}\vadjust{\goodbreak}
\]
where $p_{n,\zeta}\dvtx  \R^d \to[0,1]$ are measurable maps such that
$\Sigma_{\zeta\in\mathbb{S}} p_{n,\zeta}(x)=1$ for all $x \in\R^d$,
$n \ge1$.
Then there is a sequence of $\mathbb{S}$-valued random variables $\{
Y_n\}$ defined on
an augmentation of $(\Omega, \mathcal{F}, \IP)$ such that
\[
\IP(Y_n = \zeta\vert\mathcal{G}\vee\mathcal{Y}_0^{n-1}) =
p_{n,\zeta}(X_n),\qquad n
\ge1, \zeta\in\mathbb{S},
\]
where $\mathcal{Y}_0^{n-1} = \sigma\{Y_1, \ldots, Y_{n-1} \}$.
\end{lemm}

\begin{pf}
By suitably augmenting the space, we can assume that the probability
space $(\Omega, \mathcal{F}, \IP)$ supports an i.i.d. sequence $\{
U_n\}_{n
\ge
1}$ of Uniform $[0,1]$ random variables, independent of $\mathcal{G}$.
Let, for $n \ge1, \zeta\in\mathbb{S}$, $a_{n,\zeta}, b_{n, \zeta} \dvtx
\R^d \to[0,1]$ be measurable maps, such that
for all $x \in\R^d$:
\begin{longlist}[(iii)]
\item[(i)] $b_{n,\zeta}(x) - a_{n,\zeta}(x) = p_{n,\zeta}(x)$, $n
\ge1$, $\zeta\in\mathbb{S}$.
\item[(ii)] $[a_{n,\zeta}(x), b_{n,\zeta}(x)) \cap[a_{n,\zeta'}(x),
b_{n,\zeta'}(x)) = \varnothing$,
$\zeta, \zeta' \in\mathbb{S}$, $\zeta\neq\zeta'$, $n \ge1$.
\item[(iii)] $\bigcup_{\zeta\in\mathbb{S}} [a_{n,\zeta}(x),
b_{n,\zeta
}(x)) = [0, 1)$.
The result follows on defining
\[
Y_n = \sum_{\zeta\in\mathbb{S}} \zeta1_{[a_{n,\zeta}(X_n),
b_{n,\zeta}(X_n))}(U_n),\qquad n \ge1.
\]
\end{longlist}
\upqed\end{pf}

\subsection{\texorpdfstring{Proof of Lemma \protect\ref{ab210}}{Proof of Lemma 3.7}}
From \eqref{ab440} we have that for some $C_1 \in(0, \infty)$,
\setcounter{equation}{0}
\begin{equation}
|\tilde Q^1 - \tilde Q^2|_{\infty, T} \le C_1 |\tilde Y^1 - \tilde
Y^2|_{\infty, T}. \label{novab1}
\end{equation}
Thus, for some $C_2 \in(0, \infty)$,
%
\begin{equation}
\label{novab2}
\biggl| \tilde\E\int_0^T e^{-\gamma t} h \cdot\tilde Q^1(t) \,dt -
\tilde\E\int_0^T e^{-\gamma t} h \cdot\tilde Q^2(t) \,dt \biggr|
\le C_2 \tilde\E|\tilde Y^1 - \tilde Y^2|_{\infty, T}.\hspace*{-35pt}
\end{equation}
Next, for $t \ge0$ and $i=1,2$,
%
\begin{equation}\label{novab25}
\tilde Q^i(t+T) = \Gamma\bigl(\tilde Q^i(T) + \tilde\zeta(T+ \cdot) -
\tilde\zeta(T)\bigr).
\end{equation}
Using Assumption \ref{ab1012} and \eqref{novab1}, we now have that for
all $t \ge0$,
%
\begin{equation}\label{novab3}
|\tilde Q^1(t+T) - \tilde Q^2(t+T)| \le L |\tilde Q^1(T) - \tilde
Q^2(T)| \le LC_1 |\tilde Y^1 - \tilde Y^2|_{\infty, T}.\hspace*{-35pt}
\end{equation}
This shows that, for some $C_3 \in(0, \infty)$,
%
\begin{equation}
\label{novab4}
\biggl| \tilde\E\int_T^{\infty} e^{-\gamma t} h \cdot\tilde Q^1(t)
\,dt - \tilde\E\int_T^{\infty} e^{-\gamma t} h \cdot\tilde Q^2(t) \,dt
\biggr|
\le C_3 \tilde\E|\tilde Y^1 - \tilde Y^2|_{\infty, T}.\hspace*{-35pt}
\end{equation}
Next, for some $C_4 \in(0, \infty)$,
\begin{eqnarray}\label{novab45}
&&\biggl| \int_{[0,T]} e^{-\gamma t} p \cdot d \tilde U^1(t) - \int
_{[0,T]} e^{-\gamma t} p \cdot d \tilde U^2(t) \biggr|\nonumber\\[-2pt]
&&\qquad\le |p| \biggl[ |\tilde U^1(0) - \tilde U^2(0)| + e^{-\gamma T} |\tilde
U^1(T) - \tilde U^2(T)|
\nonumber
\\[-10pt]
\\[-10pt]
\nonumber
&&\hspace*{78pt}\quad\qquad{}+ \gamma\int_{[0,T]} |\tilde U^1(t) - \tilde U^2(t)| \,dt \biggr] \nonumber
\\[-2pt]
&&\qquad\le C_4 |\tilde Y^1 - \tilde Y^2|_{\infty, T}.\nonumber
\end{eqnarray}
Next, note that for $S > T$, $i=1,2$,
%
\begin{eqnarray}
\label{novab5}
\qquad&&\int_{(T,S]} e^{-\gamma t} p \cdot d \tilde U^i(t)
\nonumber
\\[-8pt]
\\[-8pt]
\nonumber
&&\qquad = \gamma\int_T^S
e^{-\gamma t} p \cdot[ \tilde U^i(t) -
\tilde U^i(T)] \,dt + e^{-\gamma S} p \cdot[ \tilde U^i(S) -
\tilde U^i(T)].
\end{eqnarray}
Also, using Assumption \ref{ab1012} and \eqref{ab934}, for some $C_5
\in(0, \infty)$,
\[
|\tilde U^i(t) -\tilde U^i(T)| \le C_5 |\tilde\zeta(T+ \cdot) -
\tilde
\zeta(T)|_{\infty, t-T},\qquad t \ge T,
\]
which shows that, for $i=1,2$,
\[
\tilde\E|e^{-\gamma S} p \cdot[ \tilde U^i(S) -
\tilde U^i(T)]| \to0\qquad \mbox{as } S \to\infty.
\]
Combining this observation with \eqref{novab5}, we now have, on sending
$S \to\infty$, that for $i=1,2$,
\[
\tilde E\int_{(T,\infty)} e^{-\gamma t} p \cdot d \tilde U^i(t) =
\gamma\tilde\E
\int_T^{\infty} e^{-\gamma t} p \cdot[ \tilde U^i(t) -
\tilde U^i(T)] \,dt.
\]
Thus, for some $C_6 \in(0, \infty)$,
\begin{eqnarray*}
&& \tilde\E
\biggl|\int_{(T,\infty)} e^{-\gamma t} p \cdot d \tilde U^1(t) - \int
_{(T,\infty)} e^{-\gamma t} p \cdot d \tilde U^2(t)\biggr| \\
&&\qquad \le\gamma|p| \int_{(T,\infty)} e^{-\gamma t} |
[\tilde U^1(t) -
\tilde U^1(T)] - [\tilde U^2(t) -
\tilde U^2(T)]| \,dt \\
&&\qquad \le C_6 \tilde\E\int_{(T,\infty)} e^{-\gamma t}
|\tilde Q^1(T) - \tilde Q^2(T)| \,dt,
\end{eqnarray*}
where the last equality follows on using Assumption \ref{ab1012} and
\eqref{novab25}.
Combining this with \eqref{novab1}, we now have that
\[
\tilde\E
\biggl|\int_{(T,\infty)} e^{-\gamma t} p \cdot d \tilde U^1(t) - \int
_{(T,\infty)} e^{-\gamma t} p \cdot d \tilde U^2(t)\biggr|
\le C_7 \tilde\E|\tilde Y^1 - \tilde Y^2|_{\infty, T}.
\]
The result now follows on combining the above estimate with \eqref
{novab2}, \eqref{novab4} and~\eqref{novab45}.

\subsection{\texorpdfstring{Proof of Proposition \protect\ref{tisadmisnew}}{Proof of Proposition 4.1}}
It is immediate from the construction that $T^r$ satisfies (i)--(iii) of
Definition
\ref{t-adm-defn}. We now verify that, with $\mathcal{G}= \sigma\{
\mathcal{U}_i, i
\ge1\}$,
%
\begin{equation}\label{admis1232} T^r \mbox{ satisfies (iv)}.
\end{equation}
The proof of \eqref{admis1232}
is similar to that of Theorem 5.4 in \cite{BudGho2}, which shows
that if a policy satisfies certain natural conditions (see Assumptions
5.1, 5.2, 5.3 therein),
then it is admissible (in the sense of Definition \ref{t-adm-defn} of
the current paper). The policy~$T^r$
constructed in Section \ref{subconstruct} does not exactly satisfy
conditions in Section 5 of \cite{BudGho2}, but it has similar properties.
Since most of the arguments are similar to \cite{BudGho2}, we only provide
a sketch, emphasizing only the changes that are needed. For the
convenience of the reader,
we use similar notation as in \cite{BudGho2}. Also, we suppress the
superscript $r$ from the notation.
Recall from (\ref{ab543}) that
%
%
\begin{eqnarray}\label{policy-12345}T_j(t) = \int_{[0, t]} 1_{\{
Q_{\sigma_1(j)}(u) > 0\}}
1_{\{Q_{\sigma_1(j)}(\curvp(u)) > \Theta(u)\}} {\dot
{T}^{(1)}_j}(u) \,du,\\
\eqntext{ j \in\mathbb{J}, t \ge0.}
\end{eqnarray}
In particular, Assumption 5.1 of \cite{BudGho2} is satisfied. In view
of \eqref{stpos1029},
the integrand above has countably many points (a.s.) where the value
of ${\dot{T}^{(1)}_j}$ changes from~0 to~1 (or vice versa).
Denote these points by $\{\bar{\Upsilon}_{\ell}^1\}_{\ell\in\N}$.
Set $\bar{\Upsilon}_{0}^1 =
0$. We refer to these points as the ``break-points'' of $T$.
Break-points are boundaries of the intervals of the form $\{[n r^2
\theta, (n+1) r^2 \theta)\}_{n}$ or
those of subintervals of length $\nabla^{1,k,n}_i$ or $\nabla
^{2,k,n}_i$ for some $k, i, n$
[see \eqref{abn532b}, see also \eqref{abn532a} for $\nabla^{1,k}_i$ or
$\nabla^{2,k}_i$] that are used to define the policy $T$. Next, define
$\{\bar{\Upsilon}_{\ell}^0\}_{\ell\in\N_0}$ as the countable set
of (random)
``event-points'' as defined in \cite{BudGho2} (denoted there as $\{
\Upsilon_{\ell}\}_{\ell\in\N_0}$).
These are the points where either an arrival of a job or service
completion of a job takes
place anywhere in the network. Combining
the event-points and the break-points, we get the set of
``change-points'' of the
policy $T$ denoted by $\{\Upsilon_{\ell}\}$:
\[
\{\Upsilon_{\ell} \} = \{\bar{\Upsilon}_{\ell}^0 \}\cup\{\bar
{\Upsilon}_{\ell}^1 \}.
\]
We will assume that the sequence $\{\Upsilon_{\ell}\}$ [resp. $\{\bar
{\Upsilon}_{\ell}^0\}$, $\{\bar{\Upsilon}_{\ell}^1\}$] is indexed
such that $\Upsilon_{\ell}$ [resp. $\bar{\Upsilon}_{\ell}^0$,
$\bar{\Upsilon}_{\ell}^1$] is a strictly
increasing sequence in $\ell$.

As noted earlier, \cite{BudGho2} uses the notation $\{\Upsilon_{\ell
} \}$,
instead of $\{\bar{\Upsilon}_{\ell}^0 \}$,
for event points. We have made this change of notation since $\{
\Upsilon_{\ell} \}$ here plays an identical role
as that of event-points in the proof of \cite{BudGho2}.
In particular, it is easily seen that Assumption 5.2 of \cite{BudGho2}
holds with this new definition
of $\{\Upsilon_{\ell} \}$. We will next verify Assumption 5.3 (a
nonanticipativity condition) of \cite{BudGho2}
in
Lemma \ref{appen-lemma} below.

For $i\in\mathbb{I}$ and $\ell\in
\N_0$, let $u_i^{\ell} \Df\xi_i(E_i(\Upsilon_{\ell}) + 1) -
\Upsilon_{\ell}$.
Thus,
$u_i^{\ell}$ is the
residual (exogenous) arrival time at the $i$th buffer at time $\Upsilon
_{\ell}$, unless
an arrival of the $i$th class occurred at time
$\Upsilon_{\ell}$, in which case it equals
$u_i^{\ell}=u_i(E_i(\Upsilon_{\ell})+1)$. Similarly, for $j\in
\mathbb{J}$, $\ell
\in\N_0$,
define $v_j^{\ell} \Df\eta_j(S_j(\Upsilon_{\ell}) +
1) - \Upsilon_{\ell}$.
Next, write $T(t) = \int_0^t \dot T(s) \,ds$, where $\dot T$ is right continuous.
For $i\in\mathbb{I}$, set $Q_{i,0}=0$, and for $\ell\ge1$,
$Q_{i,\ell}\Df
Q_i(\Upsilon_{\ell})$.
Also, for $j\in\mathbb{J}$, and
$\ell\ge
0$, let $\dot{T}_j^{\ell} \Df\dot{T}_j(\Upsilon_{\ell})$. Let
$\dot{T}_j^{-1}\Df0$.
Finally, define for $\ell\ge0$,
%
\begin{equation}
\label{chdef538}
\hspace*{31pt}\chi^{\ell}\Df\bigl\{(\Upsilon_{\ell'}, u_i^{\ell'}, v_j^{\ell'},
Q_i^{\ell'}, \dot{T}_j^{\ell'-1}\dvtx
i \in\mathbb{I}, j\in\mathbb{J}, \{ \mathcal{U}_i\dvtx  i \in\mathbb
{N}\} )\dvtx
\ell
'=0,\ldots, \ell
\bigr\}.\hspace*{-10pt}
\end{equation}
The definition of $\chi^{\ell}$ above is similar to that in
\cite{BudGho2}, with the
exception of the sequence $\{ \mathcal{U}_i\dvtx  i \in\mathbb{N}\}$. This
enlargement of the collection $\chi^{\ell}$
is needed due to the randomization step, involving the sequence $\{
\mathcal{U}
_i\}$, in the construction of the policy [see \eqref{ab201}].
In \cite{BudGho2}, part (iv) of the admissibility requirement (for the
smaller class of policies considered there)
was in fact shown with respect to a smaller filtration,
namely, $\bar\mathcal{F}^r((m,n))$. Here, using the above
enlargement, we
will show that part (iv) holds (for the policy in Section~\ref{subconstruct})
with $\mathcal{F}^r((m,n)) = \bar\mathcal{F}^r((m,n)) \vee\{
\mathcal{U}_i\}$. In
Lemma \ref{appen-lemma} below, we prove that~$\dot{T}(\Upsilon_{\ell
})$ is a
measurable
function of $\chi^{\ell}$, for all $\ell\in\N_0$. This shows that $T$
satisfies Assumptions 5.1--5.3 of \cite{BudGho2} with the modified
definition of $\Upsilon_{\ell}$ and $\chi^{\ell}$. Now part (iv) of the
admissibility requirement [i.e., \eqref{admis1232}] follows exactly as
the proof of Theorem
5.4 of \cite{BudGho2}. This completes the proof of the proposition.

\renewcommand{\thelemm}{\Alph{section}.\arabic{lemm}}
\setcounter{lemm}{1}
\begin{lemm}\label{appen-lemma} $\dot{T}(\Upsilon_{\ell})$ is a measurable
function of $\chi^{\ell}$, for all $\ell\in\N_0$.
\end{lemm}
\begin{pf}
Let for $m\in\N_0$, $L_m = (\bar{\Upsilon}_{m+1}^1-\bar{\Upsilon
}_{m}^1)$ denote the length of
the $m$th break-point interval.
Define $\kappa^0=0$ and $\kappa^{\ell}=\max\{m\ge0\dvtx  \bar{\Upsilon
}_{m}^1 \le\Upsilon_{\ell}\}$ for $\ell\in\N_0$. Hence, $\kappa
^{\ell}$ denotes the
number of break-points that preceded the $\ell$th change-point, and
$\bar{\Upsilon}_{\kappa^{\ell}}^1$ is the ``last'' break-point
before the $\ell$th
change-point $\Upsilon_{\ell}$ (note that $\kappa^{\ell} \le\ell
$) for all
$\ell\in\N_0$. Also, define for $\ell\in\N_0$, $\Delta^{\ell} =
(\bar{\Upsilon}_{\kappa^{\ell}}^1+ L_{\kappa^{\ell}}-\Upsilon
_{\ell})$ as the
``residual'' time for the next break-point after $\Upsilon_{\ell}$. In
particular, $\Delta_{\ell}=0$ implies that $\Upsilon_{\ell}$ itself
is a break-point.
By definition of $T$ [see \eqref{ab201}], it follows that $\kappa
^{\ell
}$, $\bar{\Upsilon}_{\kappa^{\ell}}^1$, $L_{\kappa^{\ell}}$, and,
hence, $\Delta^{\ell}$ are all measurable functions of $\chi^{\ell}$
for ${\ell} \in\N_0$.
Summarizing this, we get
%
\begin{equation}\label{ch-pt-1}\qquad
\kappa^{\ell}, \Delta^{\ell}, \bar
{\Upsilon}_{\kappa^{\ell}}^1
\mbox{ are measurable functions of } \chi^{\ell}\qquad \mbox{for }
{\ell
} \in\N_0.
\end{equation}
Using notation from \cite{BudGho2}, let $\mathcal{J}_{i} = \{j \in
\mathbb
{J}\dvtx  \sigma
_1(j) = i\}$ be the set of all activities that are associated
with the buffer $i$ and, for $a \in\{0,1\}^{\JJ}$, $\mathcal
{J}_{i}(a)$ be as
defined by equation~(5.2) of \cite{BudGho2}. Then
$\mathcal{J}_{i}(\dot{T}(t))$ denotes all activities in $\mathcal
{J}_{i}$ that are active
at time $t \ge0$,
under $T$.
Clearly, for $\ell\in\N_0$,
%
\begin{equation}\label{appen-time-relation} \Upsilon_{\ell}
=\Upsilon_{\ell-1} +
\min_{i\in\mathbb{I}} \min\{\Delta^{\ell-1}, u_i^{\ell-1},
v_j^{\ell-1}\dvtx  j \in\mathcal{J}_{i}(\dot{T}^{\ell-1})\} .
\end{equation}
For $i \in\mathbb{I}$, let $\mathcal{I}_{i}^{\ell}$ be the indicator
function of the event that at the change-point $\Upsilon_{\ell}$ an arrival
or service completion occurs at buffer $i$. More precisely, for
$i \in\mathbb{I}$ and $\ell\ge0$,
%
\begin{equation}
\qquad\quad\mathcal{I}_{i}^{\ell} =
\cases{1, &\quad$\mbox{if }
\min\{u_i^{\ell-1}, v_j^{\ell-1}\dvtx  j \in
\mathcal{J}_{i}(\dot{T}^{\ell-1})\}$\cr
& \quad\qquad$\displaystyle=\mathop{\min}_{i'\in\mathbb{I}}
\min\{\Delta^{\ell-1}, u_{i'}^{\ell-1}, v_j^{\ell-1}\dvtx  j \in
\mathcal{J}_{i'}(\dot{T}^{\ell-1})\},$\cr
 0, 
&\quad$\mbox{otherwise. }$}\label{ch-pt-2}
\end{equation}
From (\ref{ch-pt-1}) and (\ref{chdef538}), it
follows that
%
\begin{equation}
\label{new545}\qquad
\mbox{both } \mathcal{I}_{i}^{\ell} ,
\Upsilon_{\ell} \mbox{ are measurable functions of }
\chi^{\ell},\qquad \ell\in\N_0.
\end{equation}
Using \eqref{ab201} and the construction below it, along with \eqref
{new545}, it is easily checked that
$\dot{T}^{(1)}(\Upsilon_{\ell})$ is a measurable function of $\chi
^{\ell}$.
Next,
for $j \in\mathbb{J}$,
%
\begin{equation}\label{abn617} \dot{T}_j(\Upsilon_{\ell}) =
1_{\{Q_{\sigma_1(j)}(\Upsilon_{\ell}) > 0\}} 1_{\{Q_{\sigma
_1(j)}(\curvp
(\Upsilon_{\ell})) > \Theta(\Upsilon_{\ell})\}}
{\dot{T}^{(1)}_j}(\Upsilon_{\ell}).
\end{equation}
From \eqref{new545} and (\ref{chdef538}), $\Theta(\Upsilon_{\ell
})$ and
$Q_{\sigma_1(j)}(\Upsilon_{\ell})$ are $\chi^{\ell}$ measurable,
thus, so is the
first indicator in the above display. Also, since $\curvp(\Upsilon
_{\ell
})$ is either~$\Upsilon_{\ell}$ or
$\bar\Upsilon^1_{\kappa^{\ell}}$---depending on whether $\Upsilon
_{\ell}$ is in $\bigcup_{n\le p_0} \mathcal{I}_2(n)$ or not---and both
$Q_{\sigma_1(j)}(\Upsilon_{\ell})$ and~$Q_{\sigma
_1(j)}(\bar
\Upsilon^1_{\kappa^{\ell}})$
are $\chi^{\ell}$ measurable, we see that the second indicator
in~\eqref{abn617} is~$\chi^{\ell}$ measurable as well.
The lemma follows on combining the above observations.
\end{pf}

\subsection{\texorpdfstring{Proof of Lemma \protect\ref{invpr}}{Proof of Lemma 4.2}}
Since $h^n$ is equicontinuous, pre-compactness of $(f^n, g^n, h^n,
\gamma^n)$ is immediate.
Suppose now that $(f^n, g^n, h^n,\gamma^n)$ converges (in $\mathcal
{D}^{3\II
+1}$), along
some subsequence, to $(f,g,h,\gamma)$. Then $\gamma(t)=t$ for $t\ge0$
and $f,g,h \in\mathcal{C}^{\II}$. Also, for suitable measurable maps
$\tilde h_i \dvtx  [0, \infty) \to[0, 1]$, $i \in\mathbb{I}$,
\[
h_i(t) = \int_{[0,t]} \tilde h_i(s) \,ds, \qquad i \in\mathbb{I}, t\ge0.
\]
If $\psi\dvtx  [0, \infty) \to\R$
is a continuous map with compact support, then, along the above subsequence,
\[
\int_{[0, t]} \psi(f^n_i(s)) 1_{\{f^n_i(\gamma^n(s)) \le\varepsilon
_n\}} \,ds
\to
\int_{[0, t]} \psi(f_i(s)) \tilde h_i(s) \,ds,\qquad t \ge0, i \in\mathbb{I}.
\]
Suppose now that $\operatorname{supp}(\psi) \subset(\delta, \infty)$ for some
$\delta> 0$.
Then the left-hand side of the above display converges to $0$ and so for all
$t, i,$
$\int_{[0, t]} \psi(f_i(s)) \tilde h_i(s) \,ds = 0$ for such~$\psi$. Since
$\delta> 0$ is arbitrary,
we get
\[
\int_{[0, t]} 1_{\{f_i(s)=0\}} \tilde h_i(s) \,ds = h_i(t).
\]
The result follows.
\end{appendix}


\section*{Acknowledgment} We thank an anonymous referee for pointing
us to the paper \cite{K-T}.

%
%
%
%
%
%
%
%
%
%
%
%
%
%
%
%
%
%
%

%
%
%
%
%
%



\printaddresses

\end{document}